\documentclass[review,12pt]{elsarticle}

\usepackage{fancyvrb}
\usepackage{booktabs}

\usepackage{subcaption}

\usepackage{amsmath,amssymb,bm,graphicx}
\usepackage{algorithm}
\usepackage{algpseudocode}

\journal{}

\usepackage{amsthm}

\theoremstyle{definition}
\newtheorem{remark}{Remark}

\begin{document}

\begin{frontmatter}

\title{Stabilized Isogeometric Topology
Optimization on Trimmed Geometries}

\author[inst1]{Haoyu Su}

\author[inst2]{Giuliano Guarino}

\author[inst1]{Xiaodong Wei\corref{cor1}}

\cortext[cor1]{Corresponding author.}
\emailauthor{xiaodong.wei@sjtu.edu.cn}{X. Wei}

\affiliation[inst1]{
    organization={Shanghai Jiao Tong University Global College},
    city={Shanghai},
    postcode={200240},
    country={China}
}

\affiliation[inst2]{
    organization={Department of Engineering, University of Palermo},
    city={Palermo},
    postcode={90128},
    country={Italy}
}

\begin{abstract}
This paper presents a novel stabilized isogeometric topology optimization method for complex design domains represented by trimmed geometries. Minimum compliance is taken as the model problem, where a generalized Cahn-Hilliard model, a time-dependent fourth-order partial differential equation, is solved and its steady-state solutions yield the optimized structures. This approach does not require sensitivity analysis and can naturally accommodate material constraints and complex topological changes. To impose Dirichlet boundary conditions and solve the instability issue caused by trimming, Nitsche's formulation with minimal stabilization is proposed for the high-order problem. Moreover, the immersed boundary-conformal method is extended to handle complex non-design domains, allowing material-void interfaces to be treated in a unified manner. Finally, a three-dimensional biharmonic problem and several two-dimensional topology optimization examples are presented to demonstrate the effectiveness, stability, and flexibility of the proposed method in dealing with complex design domains.
\end{abstract}

\begin{keyword}
Isogeometric Topology Optimization\sep
Non-Design Domain\sep
Cahn--Hilliard\sep
Trimming\sep
Immersed Boundary\sep
Stabilized Nitsche
\end{keyword}

\end{frontmatter}

\section{Introduction}

Topology optimization (TO) seeks the optimal placement of material within a prescribed design domain to fulfill desired functions under given constraints \cite{Sigmund_TOReview,Shin_MLTOReview}. 
TO has been successfully applied in many fields, such as aerospace engineering \cite{Dong_AircraftLeadingEdge}, automotive lightweight design \cite{Lu_EVChassis}, biomedical and bio-inspired design \cite{Wang_LumbarFusion}, and additive manufacturing \cite{Zhu_AMReview}.

Conventional TO often relies on the finite element method (FEM) and thus is fundamentally disconnected with Computer-Aided Design (CAD). 
Isogeometric analysis (IGA) \cite{Hughes_IGA} has accordingly been introduced in TO to bridge the gap between CAD and analysis by using the same spline functions, such as B-splines and Non-Uniform Rational B-splines (NURBS) \cite{Piegl_NURBSBook}. Isogeometric topology optimization (ITO) has gained popularities over the last decade and has been studied with different descriptions of material evolutions, such as density-based Solid Isotropic Material with Penalization (SIMP) \cite{Gao_DensityDistributionITO,Qian_BSplineFilter}, level-set methods \cite{WangBenson_GeometricConstraintTrimmed}, Moving Morphable Components \cite{Hou_ExplicitMMCITO}, Bidirectional Evolutionary Structural Optimization \cite{Yin_IsogeometricBESO}, and phase-field methods \cite{Dede_PhaseFieldITO}. Besides being able to integrate with CAD, ITO also benefits from the smoothness of spline basis functions in terms of accuracy and structural performance \cite{Hughes_IGA,Gao_ITOReview,Gao_DensityDistributionITO}. In fact, smooth splines themselves can already play the role of filters, an important building block in FEM-based TO \cite{Sigmund_TOReview,Qian_BSplineFilter}, so ITO can usually be used without traditional filters.
 
However, the majority of ITO methods have been focused on topologically/geometrically simple design domains (e.g., rectangular domains) due to the tensor-product structure of B-splines/NURBS.
On the other hand, design domains can be complex due to the need of design iterations \cite{Dede_PhaseFieldITO}, environment constraints \cite{WangBenson_GeometricConstraintTrimmed}, and the presence of non-design domains \cite{Burman_NondesignCutTO}. They are often represented by trimmed geometries in CAD as trimming is a universal operation for modeling complex geometries. To facilitate IGA on such models, we can either adopt boundary-unfitted methods \cite{Parvizian_FCMTopologyOptimization} or reparameterize them into watertight representations \cite{Wei_ScalableReparameterizationTrimmed}. We focus on the former in this work to keep geometric representations consistent.

The main challenge of performing ITO on trimmed geometries lies in dealing with trimming in analysis. Trimming inevitably generates so-called cut elements, which are a result of arbitrary intersections between trimming surfaces/curves and a background mesh. Cut elements cause several fundamental issues that are common in boundary-unfitted methods, such as stability of discretized problems \cite{dePrenter_ImmersedReview,Buffa_MinimalStabilization}, imposition of Dirichlet boundary conditions \cite{Chouly_EssentialBCReview,Wei_IBCM}, numerical integration \cite{Marussig_TrimmingReview,Kim_TrimmingComplexTopology,Elbadry_BoundaryConformalQuadrature}, and conditioning \cite{dePrenter_FCMConditioning}.
Among these issues, particular emphasis will be placed on stability and Dirichlet boundary conditions as they are the most important issues for the proposed method.

While Dirichlet boundary conditions are often imposed by Nitsche's method due to its consistent formulation, a variety of different stabilization methods have been actively studied over the past decade. For example, the finite cell method solves the governing problem on an enlarged computational domain, where a small positive stiffness is assigned to the fictitious region \cite{Parvizian_FCMTopologyOptimization,Schillinger_FCMReview,dePrenter_FCMConditioning}. The ghost-penalty method introduces an additional term in the weak form, either by a face term that penalizes the jumps in the (high-order) normal derivatives across cut element interfaces \cite{Burman_GhostPenalty,BurmanEtAl_CutFEMGeometryPDE}, or by a volumetric term defined on cut elements that penalizes the difference between the functions and their extension counterparts from the neighboring elements \cite{BurmanHansboLarson_LockingFreeGhostPenalty}. The shifted boundary method transfers the boundary condition to an inner surrogate boundary through Taylor expansion around the true boundary, thereby circumventing the need to treat cut elements directly \cite{Main_SBMPartI,Main_SBMPartII,Atallah_SBMSecondGen}. The basis-function-removal strategy discards basis functions whose support inside the physical domain is smaller than a prescribed threshold \cite{Elfverson_BasisRemoval,Verhoosel_TrabecularBoneIGA}.

Several other methods rely on polynomial extensions from uncut/well-posed elements to construct the approximation space on cut elements. The aggregated unfitted finite element method replaces ill-posed basis functions by the polynomial extension of associated well-posed ones \cite{Badia_AggregatedUnfittedFEM}. 
Extended B-splines were introduced to modify basis functions defined on cut elements via polynomial extension of splines from neighboring intact elements \cite{Marussig_StableTrimmedGeometries}, where the extension is done in the parameter domain and it needs to be performed for all the cut elements.
On the other hand, minimal stabilization only stabilizes the problematic terms via polynomial extension, which is local to badly cut elements and can be done in the physical domain when non-trivial geometric mapping is involved (e.g., in the case of IGA) \cite{Buffa_MinimalStabilization,Antolin_OverlapMinimal,Wei_StokesOverlap}.

In this work, we propose a novel stabilized ITO method for complex design domains represented by trimmed NURBS. 
The material evolution is driven by the generalized Cahn--Hilliard equation \cite{Dede_PhaseFieldITO}, which is a time-dependent fourth-order partial differential equation (PDE) and will be solved to find the steady-state solution as the optimized structure. No sensitivity analysis is needed, whereas complex topological changes can be automatically captured without ad-hoc treatment. We adopt minimum compliance as the model problem. 
Nitsche's method is used to impose boundary and interface conditions, on top of which minimal stabilization is adapted to the time-dependent high-order problem to deal with severely cut elements. 
Moreover, non-design domains, a major source of complex design domains, are studied in detail through the Immersed Boundary-Conformal Method (IBCM) \cite{Wei_IBCM}, which embeds a layer of conformal mesh (representing the non-design domain) to each boundary/interface into a background mesh (representing the design domain). On their coupling interfaces, a unified formulation based on one-sided fluxes is proposed to handle different cases of material distributions across interfaces. 

It is worth mentioning that the cut finite element method (CutFEM), which combines Nitsche's method with ghost-penalty stabilization \cite{BurmanEtAl_CutFEMGeometryPDE}, has also been applied to TO with non-design domains \cite{Burman_NondesignCutTO}. There, the non-design domain is also represented by a conformal mesh, whereas the design domain uses a cut background mesh. But our method differs in (1) the treatment of the interface with one-sided fluxes (thus no stabilization needed), (2) tight integration with CAD and exact representation of given geometric features due to the adoption of IGA, and (3) the phase-field-based TO procedure for improved performance such as the mitigation of the mesh dependency issue and the mass conservation property for readily imposing the material constraint.

The remainder of the paper is organized as follows. 
Section~2 introduces the geometric representation via trimming. 
Section~3 reviews the generalized
Cahn--Hilliard model used for topology optimization.
Section~4 presents the phase-field-based ITO method, together with the imposition of boundary/interface conditions and the stabilized formulation. 
Section~5 develops the IBCM-based treatment of non-design domains, with the emphasis on the one-sided interface coupling. 
Section~6 presents a variety of numerical examples to demonstrate the effectiveness and flexibility of the proposed method. 
Finally, Section~7 concludes the paper and outlines possible directions for the future work. \ref{app:dimensionless} and \ref{app:time_integration} give the dimensionless formulation and the time integration with adaptivity, respectively.

\section{Geometric representation via trimming}
\label{sec:geometric_representation}

In this section, we briefly review the basic ideas of B-spline/NURBS surfaces and trimming, which is the foundation to model complex geometries in CAD systems.

\subsection{B-spline/NURBS surfaces}

B-splines and their rational extensions, Non-Uniform Rational B-Splines (NURBS), are the basis for modeling free-form geometries in Computer-Aided Design (CAD). 
B-splines are defined on a non-decreasing sequence of real numbers called a knot vector. They can be constructed recursively through the Cox--de Boor formula \cite{Piegl_NURBSBook}. Given a knot vector $\Xi=\{\xi_1,\xi_2,\ldots,\xi_{n+p+1}\}$, the corresponding univariate B-spline basis functions of degree $p$ are defined as
\begin{equation}
N_{i,0}(\xi)=
\begin{cases}
1, & \xi_i \le \xi < \xi_{i+1},\\
0, & \text{otherwise},
\end{cases}
\end{equation}
and
\begin{equation}
N_{i,p}(\xi)=
\frac{\xi-\xi_i}{\xi_{i+p}-\xi_i}N_{i,p-1}(\xi)
+
\frac{\xi_{i+p+1}-\xi}{\xi_{i+p+1}-\xi_{i+1}}N_{i+1,p-1}(\xi), \quad i=1,\ldots,n,
\end{equation}
where $n$ is the number of B-spline basis functions. B-splines are piecewise smooth polynomials and possess many desired properties for geometric modeling, such as high-order smoothness, linear independence, local support, non-negative partition of unity, and continuity control \cite{Piegl_NURBSBook}. 

The bivariate basis functions are obtained by tensor products of two univariate ones. Given a control mesh $\{\mathbf{P}_{ij}\}_{i=1,j=1}^{n,m}$, its corresponding B-spline surface is written as
\begin{equation}
\mathbf{S}(\xi,\eta)=\sum_{i=1}^{n}\sum_{j=1}^{m}N_{i,p}(\xi)M_{j,q}(\eta)\mathbf{P}_{ij}, \qquad (\xi,\eta)\in\hat{\Omega}^{0},
\label{eq:bspline_surface}
\end{equation}
where univariate B-splines $\{M_{j,q}\}_{j=1}^{m}$ are defined on another given knot vector $\mathcal{H}=\{\eta_1,\ldots,\eta_{m+q+1}\}$, $p$ and $q$ denote the polynomial degrees in the two parametric directions, and $\hat{\Omega}^{0}=(\xi_1,\xi_{n+p+1})\times(\eta_1,\eta_{m+q+1})$. Equation~\eqref{eq:bspline_surface} defines a smooth geometric mapping from the parameter domain $\hat{\Omega}^{0}$ to the physical domain $\Omega^{0}$.

NURBS extend B-splines by assigning a non-negative weight $w_{ij}$ to each control point. The rational tensor-product basis function is defined as
\begin{equation}
R_{ij}^{p,q}(\xi,\eta)=
\frac{N_{i,p}(\xi)M_{j,q}(\eta)w_{ij}}
{\sum_{a=1}^{n}\sum_{b=1}^{m}N_{a,p}(\xi)M_{b,q}(\eta)w_{ab}},
\end{equation}
and the corresponding NURBS surface is expressed as
\begin{equation}
\mathbf{S}(\xi,\eta)=\sum_{i=1}^{n}\sum_{j=1}^{m}R_{ij}^{p,q}(\xi,\eta)\mathbf{P}_{ij}.
\end{equation}
NURBS can exactly represent conic sections. When all weights are identical, the rational form reduces to a B-spline surface.

\subsection{Trimming}

Trimming provides a simple yet effective way to model complex geometries~\cite{Marussig_TrimmingReview,Kim_TrimmingComplexTopology}. 
Let $\mathbf{F}:\hat{\Omega}^{0}\rightarrow\Omega^{0}$ be a B-spline/NURBS surface patch; see Fig. 1(a) as an example. A trimming curve is a B-spline/NURBS curve lying in $\hat{\Omega}^{0}$, denoted as  $\mathbf{C}: [0,1] \to \hat{\Gamma}^{t} \subset \hat{\Omega}^{0}$. Such curves, together with the boundaries, form oriented loops to specify the visible/active part of the surface patch; see Fig. 1(b), where by convention, the region left to the loops, denoted as $\hat{\Omega}$, will serve as the (active) definition domain of $\mathbf{F}$ after trimming. As a result, the image $\mathbf{F}(\hat{\Omega})$ yields the desired geometry $ \Omega$; see Fig. 1(c). Note that the parameterization (indicated by the mesh lines) is not conformal to the geometric boundary. Indeed, a trimming curve can intersect with the background mesh in an arbitrary way.

In a trimming operation, everything related to the parameterization $\mathbf{F}$, such as knot vectors, B-splines, and control points, remains the same except the definition domain, which changes from $\hat{\Omega}^0$ to $\hat{\Omega}$. In practice, a geometric model is composed of multiple trimmed patches, where trimming curves are often determined through the surface-surface intersection algorithm~\cite{Marussig_TrimmingReview,Piegl_NURBSBook}. As trimming curves are introduced independently of the surface parameterization, it is highly flexible to model holes, cut-outs, and many other  geometric features, which is why it is a universal operation supported in almost every CAD system for modeling complex geometries.

\begin{figure}[htbp]
    \centering
    \newlength{\trimfigheight}
    \setlength{\trimfigheight}{0.27\textwidth}

    \begin{subfigure}[t]{0.31\textwidth}
        \vspace{0pt}
        \centering
        \includegraphics[height=\trimfigheight,keepaspectratio]{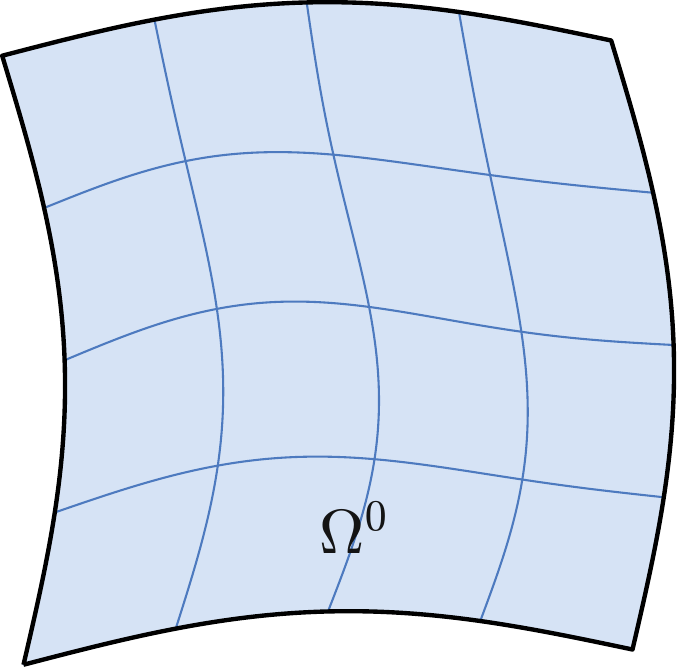}
        \caption{Non-trimmed surface patch in the physical domain}
        \label{fig:trimming_a}
    \end{subfigure}
    \hfill
    \begin{subfigure}[t]{0.31\textwidth}
        \vspace{0pt}
        \centering
        \includegraphics[height=\trimfigheight,keepaspectratio]{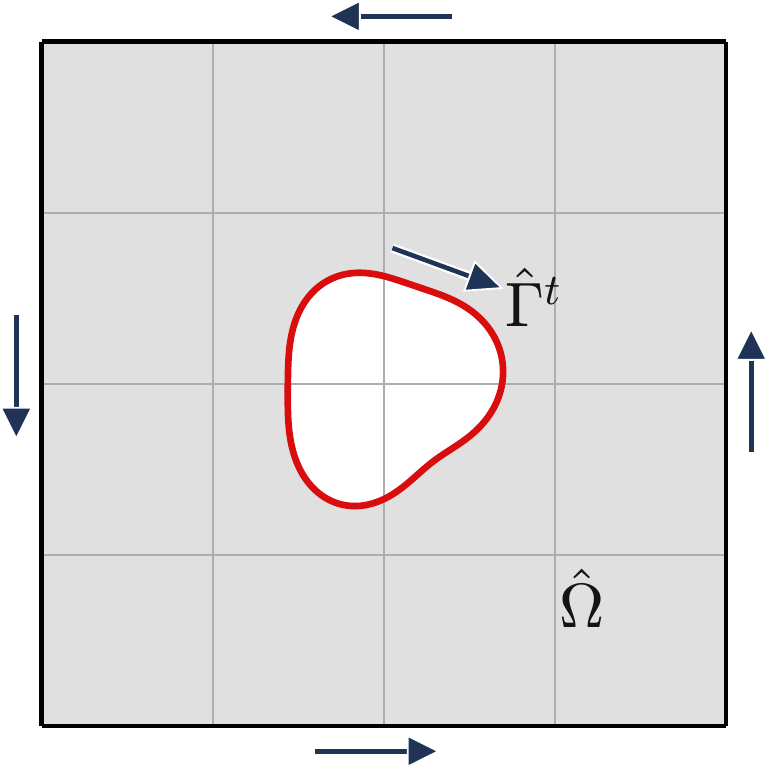}
        \caption{Active parameter domain defined by the oriented loops}
        \label{fig:trimming_b}
    \end{subfigure}
    \hfill
    \begin{subfigure}[t]{0.31\textwidth}
        \vspace{0pt}
        \centering
        \includegraphics[height=\trimfigheight,keepaspectratio]{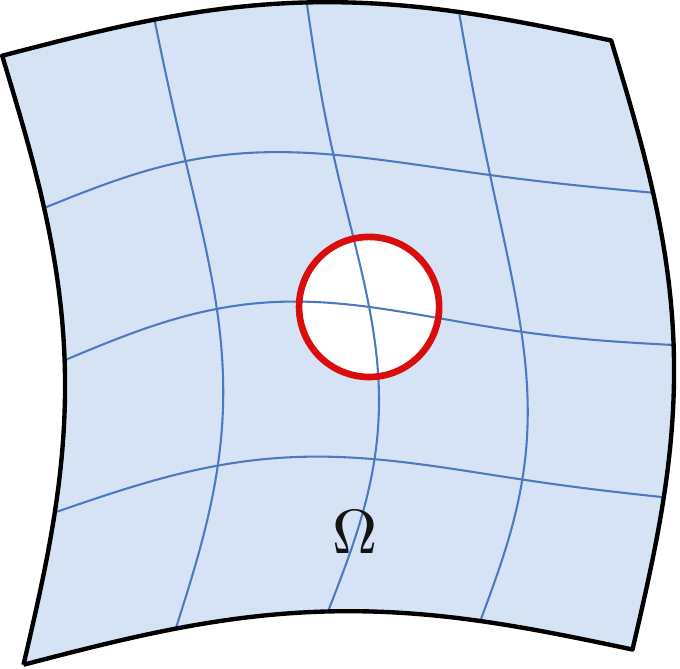}
        \caption{Trimmed surface patch in the physical domain}
        \label{fig:trimming_c}
    \end{subfigure}

    \caption{Illustration of trimming for a single surface patch, where isoparametric lines are used to show the mesh of the underlying parameterization.}
    \label{fig:trimming_illustration}
\end{figure}

The parameter domain $\hat{\Omega}^0$ is partitioned by the knot vectors into a parametric mesh,
\begin{equation}
\begin{aligned}
\hat{\mathcal{M}}^0
={}& \left\{ \hat{K}^0=(\xi_i,\xi_{i+1})\times(\eta_j,\eta_{j+1}) :
\xi_i<\xi_{i+1},\ \eta_j<\eta_{j+1}, \right.\\
&\left. i=1,\ldots,n+p, \text{and} \,j=1,\ldots,m+q \right\},
\end{aligned}
\end{equation}
where $\hat{K}^0$ is a parametric element. The corresponding (physical) mesh is obtained by the geometric mapping $\mathbf{F}$,
\begin{equation}
\mathcal{M}^0 = \{ K^0=\mathbf{F}(\hat{K}^0): \hat{K}^0 \in \hat{\mathcal{M}}^0 \} ,
\end{equation}
where $K^0$ is a (physical) element. 

With trimming, the elements are a result of $K^0 \cap \Omega$ (or equivalently, $\hat{K}^0 \cap \hat{\Omega}$). The mesh of a trimmed geometry is obtained by
\begin{equation}
\mathcal{M} = \{ K:=K^0 \cap \Omega: K^0 \in \mathcal{M}^0\, \text{and}\, K \neq \varnothing \},
\end{equation}
where $K^0$ is called the parent of $K$. Note that when $K^0 \cap \Omega = \varnothing$, $K^0$ is totally outside of $\Omega$ and it is an \emph{inactive element}. On the other hand, $K$ is called a \emph{complete element} if $K=K^0$; it is a \emph{cut element} if $K \neq K^0$ and $K \neq \varnothing$. Clearly, a complete element is fully contained in $\Omega$, whereas a cut element only partially intersects $\Omega$. Cut elements are inevitable due to the arbitrary intersection between trimming curves and the B-spline mesh, and thus, challenging issues for boundary-unfitted methods, such as imposition of boundary conditions and stabilization, remain to be addressed, which are the focus of the present work.

\section{Generalized Cahn--Hilliard model}
\label{sec:generalized_cahn_hilliard}

We review the generalized Cahn--Hilliard model \cite{Dede_PhaseFieldITO} in this section, which is adopted to perform topology optimization (TO) in this work. Minimum compliance is taken as the model problem, which aims to find the stiffest structure under a given amount of material. Equivalently speaking, it minimizes the strain energy while satisfying the equilibrium equations and the material constraint. Let $\Omega\subset\mathbb{R}^{d}$, $d\in\{2,3\}$, be a domain in which the material layout is to be optimized. $\Omega$ is called a \emph{design domain}. Here we first consider the case where $\Omega$ can be represented by tensor-product B-splines. Let $\rho:\Omega\rightarrow[0,1]$ be a density field representing the material layout, where $\rho=1$ denotes pure material, $\rho=0$ means void, and intermediate density values are allowed to facilitate the numerical procedure. Given a prescribed amount of material $V_0$, the minimum compliance problem is stated as
\begin{equation}
\begin{aligned}
{\operatorname{minimize}}\quad
&J_E\bigl(\rho,\mathbf{u}(\rho)\bigr)
=
\int_{\Omega}
\boldsymbol{\sigma}\bigl(\rho,\mathbf{u}(\rho)\bigr)
:\boldsymbol{\varepsilon}\bigl(\mathbf{u}(\rho)\bigr)
\,\mathrm{d}\Omega,
\\[2pt]
\text{subject to}\quad
&\left\{
\begin{array}{@{}l@{\qquad}l@{}}
-\nabla\cdot\boldsymbol{\sigma}\bigl(\rho,\mathbf{u}(\rho)\bigr)
=\mathbf{f}
&\text{in }\Omega,
\\
\mathbf{u}
=\mathbf{g}
&\text{on }\Gamma_D,
\\
\boldsymbol{\sigma}\bigl(\rho,\mathbf{u}(\rho)\bigr)\mathbf{n}
=\mathbf{h}
&\text{on }\Gamma_N,
\end{array}
\right.
\\[5pt]
&\int_{\Omega}\rho\,\mathrm{d}\Omega=V_0.
\end{aligned}
\label{eq:min_compliance_problem_review}
\end{equation}
Here $\mathbf{u}(\rho)$ is the equilibrium displacement associated with the current density field $\rho$, the body force $\mathbf{f}$, the prescribed displacement $\mathbf{g}$ on the Dirichlet boundary $\Gamma_D$, and the traction $\mathbf{h}$ on the Neumann boundary $\Gamma_N$. $\boldsymbol{\sigma}$ and $\boldsymbol{\varepsilon}$ are the Cauchy stress and strain, respectively, and $\bold{n}$ is the unit outward normal to the boundary. Note that $\partial\Omega=\overline{\Gamma_D\cup\Gamma_N}$, $\Gamma_D\cap\Gamma_N=\emptyset$, and the factor $1/2$ in the definition of the strain energy $J_E$ has been omitted. In the equilibrium equation, we consider linear elasticity under the small displacement assumption with an isotropic material model.

Various methods exist to solve Problem~\eqref{eq:min_compliance_problem_review}. In this work, we follow a phase-field-based TO method \cite{Dede_PhaseFieldITO,Sheng_ThinWalledComplex} mainly because (1) it solves the optimization problem by solving a time-dependent PDE and thus does not require a dedicated optimizer, and (2) it can handle arbitrary topological changes without ad-hoc treatment while satisfying the material constraint in a straightforward manner. The method features a generalized Cahn--Hilliard model that equivalently converts the optimization problem to the following time-dependent, fourth-order PDE,
\begin{equation}
\left\{
\begin{aligned}
\frac{\partial\rho}{\partial t}
&=\nabla\cdot\bigl(M(\rho)\nabla z\bigr)
&&\text{in }\Omega\times(0,T],
\\
M(\rho)\nabla z\cdot\mathbf{n}
&=0
&&\text{on }\partial\Omega\times(0,T],
\\
\nabla\rho\cdot\mathbf{n}
&=0
&&\text{on }\partial\Omega\times(0,T],
\\
\rho
&=\rho_0:=\frac{V_0}{|\Omega|}
&&\text{in }\Omega,\quad t=0,
\end{aligned}
\right.
\label{eq:gch_review}
\end{equation}
with
\begin{equation}
\left\{
\begin{aligned}
z(\rho,\mathbf{u})
&=-\frac{\gamma}{E_0}
\frac{\partial\psi_E}{\partial\rho}(\rho,\mathbf{u})
+\frac{\mathrm{d}\psi_B}{\mathrm{d}\rho}(\rho)
+\lambda\frac{\mathrm{d}\psi_I}{\mathrm{d}\rho}(\rho),
\\
\psi_E(\rho,\mathbf{u})
&=\boldsymbol{\sigma}(\rho,\mathbf{u})
:\boldsymbol{\varepsilon}(\mathbf{u}),
\\
\psi_B(\rho)
&=\rho^2(1-\rho)^2
+\beta_1\left(10^{-\beta_2\rho}+10^{\beta_2(\rho-1)}\right),
\\
\psi_I(\rho)
&=\frac{1}{2}\nabla\rho\cdot\nabla\rho.
\end{aligned}
\right.
\label{eq:gch_quantities_review}
\end{equation}
Here $|\Omega|$ denotes the measure of $\Omega$. $\rho$ is the phase field used to represent the density field in Problem~\eqref{eq:min_compliance_problem_review}. In other words, pure material ($\rho=1$) and void ($\rho=0$) are represented as the two phases in this model, whereas the intermediate material state $(0<\rho<1)$ is described by a smeared interface between the two phases. $M(\rho)$ is a non-negative function called mobility and we set it to 1. $\psi_E$ is the strain energy density, whereas $\psi_B$ and $\psi_I$ are the bulk and interface energy density of the classical Cahn--Hilliard model, respectively. Parameters $\gamma$ and $\lambda$ are introduced to balance different terms. $\gamma$ controls the portion of the strain energy relative to the bulk and interface energies. If $\gamma$ is too small, the evolution is dominated by the bulk and interface energies, and it will not sufficiently respect the strain energy to yield a desired structure. On the other hand, an excessively large $\gamma$ may hinder the ability of the Cahn--Hilliard model to drive the material evolution. $\lambda$ weights the interface energy and controls the width of the smeared interface. Smaller $\lambda$ yields a sharper interface but requires a finer mesh to resolve it. The specific choices of $\gamma$ and $\lambda$ are problem dependent and will be discussed later. Parameters $\beta_1$ and $\beta_2$ are used to control, respectively, the
magnitude and steepness of the exponential terms in the bulk
energy density. Following \cite{Dede_PhaseFieldITO}, we set $\beta_1=0.5$ and $\beta_2=50$. Note that the material constraint is imposed as an initial condition $\rho_0$. This is possible because the generalized Cahn--Hilliard model has the mass conservation property.

By finding the steady-state solution to Eq.~\eqref{eq:gch_review}, we can obtain the optimized structure corresponding to Problem~\eqref{eq:min_compliance_problem_review}. Isogeometric analysis (IGA) is well-suited to solve this fourth-order problem. In the boundary-fitted setting, the trial and test function spaces are defined as
\begin{equation}
\begin{aligned}
\mathcal{H}
&=\left\{\phi\in H^2(\Omega):\nabla\phi\cdot\mathbf{n}=0
\text{ on }\partial\Omega\right\},
\\
\mathcal{U}
&=\left\{\mathbf{u}\in[H^1(\Omega)]^d:
\mathbf{u}=\mathbf{g}\text{ on }\Gamma_D\right\},
\\
\mathcal{V}
&=\left\{\mathbf{v}\in[H^1(\Omega)]^d:
\mathbf{v}=\mathbf{0}\text{ on }\Gamma_D\right\},
\end{aligned}
\label{eq:boundary_fitted_spaces}
\end{equation}
where $d\in\{2,3\}$ is the spatial dimension.
The variational formulation of
Eq.~\eqref{eq:gch_review} is to find $\rho\in\mathcal{H}$ and $\mathbf{u}\in\mathcal{U}$ such that
\begin{equation}
\left\{
\begin{aligned}
R_{\rho}(\rho,\mathbf{u})(\phi)
:={}&
\int_{\Omega}
\frac{\partial\rho}{\partial t}\,\phi\,\mathrm{d}\Omega
-\int_{\Omega}
z\,\nabla\cdot\bigl(M(\rho)\nabla\phi\bigr)\,\mathrm{d}\Omega = 0
,
&&\forall\phi\in\mathcal{H},
\\
R_{\mathbf{u}}(\rho,\mathbf{u})(\mathbf{v})
:={}&
\int_{\Omega}
\boldsymbol{\sigma}(\rho,\mathbf{u})
:\boldsymbol{\varepsilon}(\mathbf{v})\,\mathrm{d}\Omega
-\int_{\Omega}\mathbf{f}\cdot\mathbf{v}\,\mathrm{d}\Omega
-\int_{\Gamma_N}\mathbf{h}\cdot\mathbf{v}\,\mathrm{d}\Gamma = 0
,
&&\forall\mathbf{v}\in\mathcal{V}.
\end{aligned}
\right.
\label{eq:boundary_fitted_residuals}
\end{equation}
This coupled nonlinear system is discretized in time by the generalized-$\alpha$ method. At each time step, the system of nonlinear equations is solved by Newton's method. The entire solution procedure does not require sensitivity analysis or user intervention. 

The resulting phase field exhibits a smeared interface. Within the smeared interface where intermediate phase/density values ($0<\rho<1$) are involved, the corresponding material property is scaled according to the Solid Isotropic Material with Penalization (SIMP) \cite{BendsoeSigmund_MaterialInterpolation}. More specifically, we have
\begin{equation}
\boldsymbol{\sigma}(\rho,\mathbf{u})
=g(\rho)\,\mathbb{C}_0:\boldsymbol{\varepsilon}(\mathbf{u}), \quad \boldsymbol{\varepsilon}(\mathbf{u}) = \frac{1}{2} (\nabla \mathbf{u} + \nabla \mathbf{u}^T),
\end{equation}
where $\mathbb{C}_0$ is the elasticity tensor of the pure material. 
$g(\rho)$ is given as
\begin{equation}
\begin{gathered}
g(\rho)=f(\rho)^P,\qquad P=3,
\\[3pt]
f(\rho)=
\begin{cases}
\rho_{\min},
& \rho<\rho_1,
\\
\rho_{\min}+(1-\rho_{\min})\,b(\rho),
& \rho_1\leq\rho<\rho_6,
\\
1,
& \rho\geq\rho_6,
\end{cases}
\\[3pt]
\end{gathered}
\label{eq:material_interpolation_function}
\end{equation}
where \(b(\rho)\) is a \(C^2\)-continuous cubic B-spline curve with knot vector
$\bigl\{
\rho_1,\allowbreak
\rho_1,\allowbreak \rho_1, \allowbreak\rho_1,\allowbreak \rho_2, \allowbreak \rho_3, \rho_4, \allowbreak \rho_5, \rho_6, \allowbreak \rho_6,\allowbreak \rho_6, \allowbreak \rho_6\allowbreak
\bigr\}$
and control points
\(\{0,0,0,0,1,1,1,1\}\). We set
$
\rho_{\min}=0.1,\allowbreak
\rho_1=0.055,\allowbreak
\rho_2=0.065,\allowbreak
\rho_3=0.075,\allowbreak
\rho_4=0.925,\allowbreak
\rho_5=0.935,\allowbreak
\rho_6=0.945.
$

However, constructing a test function space according to
Eq.~\eqref{eq:boundary_fitted_spaces} is challenging for complex
surfaces even in the boundary-fitted setting \cite{Sheng_ThinWalledComplex}, and it becomes generally impossible to strongly impose the essential conditions in the boundary-unfitted setting. Therefore, the essential boundary conditions,
including both the displacement field and the phase field, will need to be
enforced in the weak form through, for instance, Nitsche's method. Moreover,
severely cut elements need to be dealt with to retain stability.

\section{Topology optimization with trimmed geometries}
\label{sec:trimmed_topology_optimization}

In this section, we introduce the proposed isogeometric topology optimization (ITO) method with design domains represented by trimmed geometries. In what follows, we start by imposing essential boundary conditions using Nitsche's method. We then proceed to extend the minimal stabilization method \cite{Buffa_MinimalStabilization} to deal with the involved high-order PDE.

\subsection{Imposition of essential conditions by Nitsche's method}
\label{subsec:nitsche_essential_conditions}

When design domains are represented by trimmed geometries, essential boundary conditions, such as $\nabla\rho\cdot\mathbf{n}=0$ and $\mathbf{u}=\mathbf{g}$, need to be imposed weakly. We adopt Nitsche's formulation for this purpose. 

Recall that in the generalized Cahn--Hilliard model, two fields are coupled together and they need to be solved simultaneously. Now we introduce $\mathcal{H}=H^2(\Omega)$ and $\mathcal{V}=[H^1(\Omega)]^d$ as the test function spaces of the phase field and the displacement, respectively, where no essential conditions are imposed in this definition. Two residuals are defined accordingly. First, given $\boldsymbol{\mu}=[\gamma,\lambda]^{\mathsf T}$ and test function $\phi\in\mathcal{H}$, the phase-field residual is defined as
\begin{equation}
\begin{aligned}
R_{\rho}(\rho,\mathbf{u};\boldsymbol{\mu})(\phi)
:=\;&
\int_{\Omega}
\frac{\partial\rho}{\partial t}\,\phi\,\mathrm{d}\Omega
-
\int_{\Omega}
z(\rho,\mathbf{u};\boldsymbol{\mu})
\nabla\!\cdot\!\bigl(M(\rho)\nabla\phi\bigr)\,\mathrm{d}\Omega
\\[3pt]
&+
\underbrace{
\int_{\partial\Omega}
z(\rho,\mathbf{u};\boldsymbol{\mu})M(\rho)
(\nabla\phi\cdot\mathbf{n})\,\mathrm{d}\Gamma
}_{\text{consistency term}}
\\[3pt]
&+
\underbrace{
\int_{\partial\Omega}
\frac{C_{\mathrm{pen}}^{\rho}}{h}
(\nabla\rho\cdot\mathbf{n})
(\nabla\phi\cdot\mathbf{n})\,\mathrm{d}\Gamma
}_{\text{penalty term}},
\end{aligned}
\label{eq:Rrho_nitsche}
\end{equation}
where $h$ is a mesh-size parameter independent of trimming, and $C_{\mathrm{pen}}^{\rho}>0$ is a penalty parameter. There are two additional terms compared to Eq.~\eqref{eq:boundary_fitted_residuals} due to the fact that $\nabla\phi\cdot\mathbf{n}\neq0$. The consistency term naturally arises from integration by parts, whereas the penalty term imposes the boundary condition $\nabla\rho\cdot\mathbf{n}=0$. Also note that the symmetry term, which usually appears in Nitsche's formulation, is not introduced here because the system is nonlinear and unsymmetric.

Second, given test function $\mathbf{v}\in\mathcal{V}$, the displacement residual is defined in a classical way as
\begin{equation}
\begin{aligned}
R_{\mathbf{u}}(\rho,\mathbf{u})(\mathbf{v})
:=\;&
\int_{\Omega}
\boldsymbol{\sigma}(\rho,\mathbf{u})
:\boldsymbol{\varepsilon}(\mathbf{v})\,\mathrm{d}\Omega
-
\int_{\Omega}\mathbf{f}\cdot\mathbf{v}\,\mathrm{d}\Omega
-
\int_{\Gamma_N}\mathbf{h}\cdot\mathbf{v}\,\mathrm{d}\Gamma
\\[3pt]
&-
\underbrace{
\int_{\Gamma_D}
\boldsymbol{\sigma}(\rho,\mathbf{u})\mathbf{n}
\cdot\mathbf{v}\,\mathrm{d}\Gamma
}_{\text{consistency term}}
-
\underbrace{
\int_{\Gamma_D}
\boldsymbol{\sigma}(\rho,\mathbf{v})\mathbf{n}
\cdot(\mathbf{u}-\mathbf{g})\,\mathrm{d}\Gamma
}_{\text{symmetry term}}
\\[3pt]
&+
\underbrace{
\int_{\Gamma_D}
\frac{C_{\mathrm{pen}}^u}{h}
(\mathbf{u}-\mathbf{g})\cdot\mathbf{v}\,\mathrm{d}\Gamma
}_{\text{penalty term}},
\end{aligned}
\label{eq:Ru_nitsche}
\end{equation}
where $C_{\mathrm{pen}}^u>0$ is another penalty parameter and we have a non-zero test function $\mathbf{v}\neq\mathbf{0}$ on $\Gamma_D$ in general.

The goal is to find $\mathbf{u}\in\mathcal{V}$ and $\rho\in\mathcal{H}$ such that
\begin{equation}
\left\{
\begin{aligned}
R_{\rho}(\rho,\mathbf{u};\boldsymbol{\mu})(\phi)&=0,
&&\forall\phi\in\mathcal{H},
\\
R_{\mathbf{u}}(\rho,\mathbf{u})(\mathbf{v})&=0,
&&\forall\mathbf{v}\in\mathcal{V}.
\end{aligned}
\right.
\label{eq:trimmed_coupled_problem}
\end{equation}
While it is possible to introduce penalty parameters $C_{\mathrm{pen}}^{\rho}$ and $C_{\mathrm{pen}}^u$ local to each element to guarantee the well-posedness of the discretized problem \cite{Buffa_MinimalStabilization}, they can be arbitrarily large when a cut element becomes arbitrarily small. Therefore, a robust formulation regardless of trimming is needed.

\subsection{Minimal stabilization for the generalized Cahn--Hilliard problem}
\label{subsec:minimal_stabilization_high_order}

We extend minimal stabilization \cite{Buffa_MinimalStabilization} to obtain a stabilized formulation of Problem~\eqref{eq:trimmed_coupled_problem}. Minimal stabilization was first introduced to deal with elliptic problems on trimmed geometries, and later it was  generalized to overlapping multi-patch geometries \cite{Antolin_OverlapMinimal}, the Stokes problem \cite{Wei_StokesOverlap}, and shell problems \cite{Guarino_ShellStabilization}. It guarantees well-posedness regardless of trimming conditions and provides a priori estimates with optimal convergence rates.

\begin{figure}[htbp]
    \centering
    \includegraphics[width=0.42\textwidth]{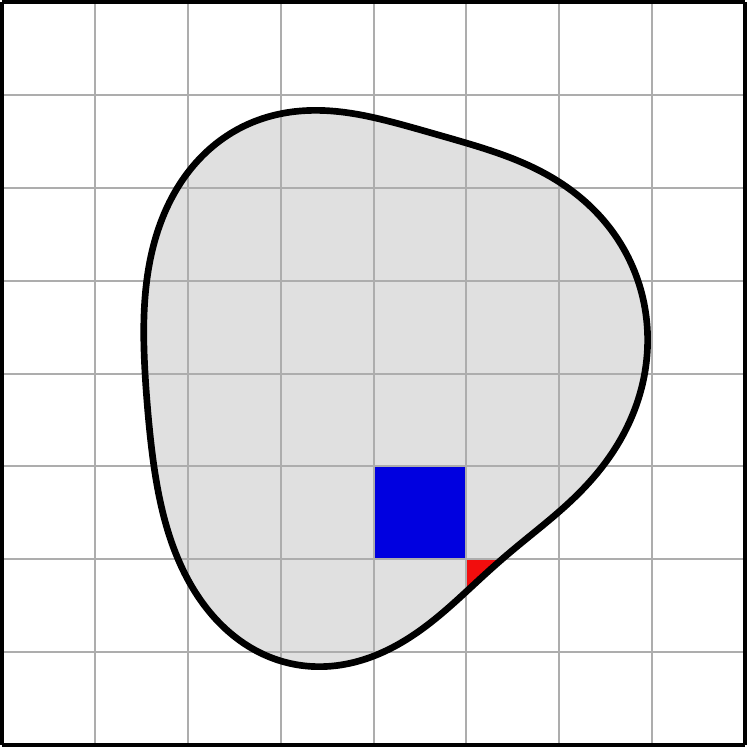}
    \caption{Schematic of a bad cut element (red) and its associated good neighbor (blue).}
    \label{fig:good_bad_ms}
\end{figure}

Recall that in Section~\ref{sec:geometric_representation}, we have classified three types of elements: cut elements, complete elements, and inactive elements. Among them, inactive elements are not used in either design or analysis, so we will not consider them in what follows. For every element $K$ and its parent $K^0$ before trimming, we define its effective area/volume ratio as
\begin{equation}
\eta_K=\frac{|K|}{|K^0|}.
\label{eq:area_ratio}
\end{equation}
Given a threshold $\theta\in(0,1)$, $K$ is called a \emph{bad element} if $\eta_K<\theta$; otherwise it is a good element, which can be either a complete element or a cut element with a big effect area/volume. In practice, we set $\theta \leq 0.1$. We assume that for each bad element $K_b$, there exists a good element in its direct neighborhood, so it is called a good neighbor ($K_g$) to $K_b$; see Fig.~\ref{fig:good_bad_ms}. This assumption generally holds unless the background mesh is too coarse with respect to the trimming features. In this case, mesh refinement will be performed.

It has been shown that it is the boundary/interface flux terms (i.e., normal derivatives) from the bad elements that cause the instability issue \cite{Buffa_MinimalStabilization}. Minimal stabilization replaces such terms with the counterparts from their good neighbors via polynomial extension. In the case of the generalized Cahn--Hilliard problem, we need to stabilize the consistency term and the penalty term of Eq.~\eqref{eq:Rrho_nitsche}. The treatment of Eq.~\eqref{eq:Ru_nitsche} follows the same as in \cite{Buffa_MinimalStabilization} and will be briefly commented later.

Here we focus on the treatment of Eq.~\eqref{eq:Rrho_nitsche}. Three terms need to be stabilized, including $\nabla\phi\cdot\mathbf{n}$, $\nabla\rho\cdot\mathbf{n}$, and $z(\rho,\mathbf{u};\boldsymbol{\mu})$. Among them, the treatment of the first two terms is the same, whereas $z(\rho,\mathbf{u};\boldsymbol{\mu})$ is nonlinear and needs additional care. We start with $\nabla\phi\cdot\mathbf{n}$. Given the test function $\phi$ and a bad element $K_b$, we first move to its good neighbor $K_g$. Let $K_g^{\mathrm{bb}}$ be the bounding box of $K_g$ in the physical domain and $\{b_l(\bold{x})\}$ be a set of Bernstein polynomials defined on $K_g^{\text{bb}}$. It is sufficient to choose the degrees of $b_l(\bold{x})$ to be the same as those of the spline basis functions. For simplicity, we assume the same degree $p$ in every spatial direction. $\{b_l(\bold{x})\}$ spans the polynomial space $\mathbb{Q}^p(K_g^{\text{bb}})$. Next, $\phi$ is projected onto $\mathbb{Q}^p(K_g^{\text{bb}})$. As a result, we have
\begin{equation}
\mathcal{P}(\phi)
=
\sum_{l=1}^{(p+1)^d}c_{\phi,l}b_l(\mathbf{x}),
\qquad \mathbf{x}\in K_g^{\text{bb}},
\label{eq:projection_good}
\end{equation}
where $\mathcal{P}$ denotes an $L^2$-projection. Finally, we move back to the bad element by polynomial extension of $\{b_l(\bold{x})\}$. Equivalently speaking, the definition domain of Eq.~\eqref{eq:projection_good} is enlarged by incorporating the bad element $K_b$. As a result, we have
\begin{equation}
\tilde{\phi}
:=
\mathcal{E} ( \mathcal{P} (\phi) )
=
\sum_{l=1}^{(p+1)^d}c_{\phi,l}b_l(\mathbf{x}),
\qquad \mathbf{x}\in K_g^{\text{bb}}\cup\overline{K_b},
\label{eq:extension_good}
\end{equation}
where $\mathcal{E}$ represents the extension operator and $\overline{K_b}$ denotes the closure of $K_b$. The flux terms in Eq.~\eqref{eq:Rrho_nitsche}, such as $\nabla\phi\cdot \mathbf{n}$, will be evaluated using $\tilde{\phi}$. As a result, the stabilized flux term becomes $\nabla\tilde{\phi}\cdot \mathbf{n}$. $\nabla\rho\cdot \mathbf{n}$ is treated in the same manner with its stabilized counterpart $\nabla \tilde{\rho} \cdot \mathbf{n}$.

Next, recall that in Eq.~\eqref{eq:gch_quantities_review}, the nonlinear term $z(\rho,\mathbf{u};\boldsymbol{\mu})$ depends on three components, $\psi_E(\rho,\mathbf{u})$, $\psi_B(\rho)$,
and $\psi_I(\rho)$. Following the same spirit as before, we stabilize all their spatial derivatives (with respect to $x$, $y$ and $z$) on the boundary/interface. Note that $\psi_E(\rho,\mathbf{u})$ and $\psi_I(\rho)$ are defined with $\nabla\mathbf{u}$ and $\nabla\rho$, respectively, whereas $\psi_B(\rho)$ does not involve $\nabla \rho$. Therefore, $z(\rho,\mathbf{u};\boldsymbol{\mu})$ is stabilized as
\begin{equation}
\widetilde{z}(\rho,\mathbf{u};\boldsymbol{\mu})
:=
-\frac{\gamma}{E_0}
\frac{\partial\psi_E}{\partial\rho}
\bigl(\rho,\widetilde{\mathbf{u}}\bigr)
+
\frac{\mathrm{d}\psi_B}{\mathrm{d}\rho}(\rho)
+
\lambda
\frac{\mathrm{d}\psi_I}{\mathrm{d}\rho}
\bigl(\widetilde{\rho}\bigr).
\label{eq:z_stabilized_minimal}
\end{equation}
Note that the density dependence in the strain-energy contribution is not
reconstructed. Thus, the density interpolation $g(\rho)$ is evaluated using
the original phase field on $K_b$. The strain
$\boldsymbol{\varepsilon}(\widetilde{\mathbf{u}})$ is evaluated using the
reconstructed displacement field. The bulk contribution is also left
unchanged because $\psi_B$ depends only on the pointwise value of $\rho$
and does not involve spatial derivatives. For the interface contribution,
$\mathrm{d}\psi_I/\mathrm{d}\rho(\widetilde{\rho})
=-\Delta\widetilde{\rho}$, it is evaluated using the reconstructed phase field. 

\begin{remark}
\label{rem:unified_stabilization}
Using $\tilde{\rho}$ everywhere on the right-hand side of
Eq.~\eqref{eq:z_stabilized_minimal} provides another feasible
stabilization, leading to
\begin{equation}
\widetilde{z}(\rho,\mathbf{u};\boldsymbol{\mu})
:=
-\frac{\gamma}{E_0}
\frac{\partial\psi_E}{\partial\rho}
\bigl(\widetilde{\rho},\widetilde{\mathbf{u}}\bigr)
+
\frac{\mathrm{d}\psi_B}{\mathrm{d}\rho}
\bigl(\widetilde{\rho}\bigr)
+
\lambda
\frac{\mathrm{d}\psi_I}{\mathrm{d}\rho}
\bigl(\widetilde{\rho}\bigr).
\label{eq:z_stabilized_full}
\end{equation}
Based on our numerical tests in Section~\ref{sec:numerical_examples},
no noticeable difference between 
Eqs.~(\ref{eq:z_stabilized_minimal}, \ref{eq:z_stabilized_full}) can be observed in terms of accuracy and conditioning. This is because polynomial extension is of high accuracy. In practice, Eq.~\eqref{eq:z_stabilized_full} treats different terms in a unified way, so it simplifies implementation and thus is adopted.
\end{remark}

In the case of Eq.~\eqref{eq:Ru_nitsche}, only the consistency term and the symmetry term need stabilization, whereas the penalty term does not need it because it does not involve a normal derivative. The procedure follows the same. As a result, the stabilized formulations are obtained as follows,
\begin{equation}
\begin{aligned}
\tilde{R}_{\rho}(\rho,\mathbf{u};\boldsymbol{\mu})(\phi)
={}&
\int_{\Omega}
\frac{\partial\rho}{\partial t}\,\phi\,\mathrm{d}\Omega
-
\int_{\Omega}
z(\rho,\mathbf{u};\boldsymbol{\mu})
\nabla\!\cdot\!\bigl(M(\rho)\nabla\phi\bigr)\,\mathrm{d}\Omega
\\
&+
\int_{\partial\Omega}
\tilde z(\rho,\mathbf{u};\boldsymbol{\mu})M(\rho)
(\nabla\tilde{\phi}\cdot\mathbf{n})\,\mathrm{d}\Gamma
\\
&+
\int_{\partial\Omega}
\frac{C_{\mathrm{pen}}^{\rho}}{h}
(\nabla\tilde{\rho}\cdot\mathbf{n})
(\nabla\tilde{\phi}\cdot\mathbf{n})\,\mathrm{d}\Gamma,
\\[6pt]
\tilde{R}_{\mathbf{u}}(\rho,\mathbf{u})(\mathbf{v})
={}&
\int_{\Omega}
\boldsymbol{\sigma}(\rho,\mathbf{u})
:\boldsymbol{\varepsilon}(\mathbf{v})\,\mathrm{d}\Omega
-
\int_{\Omega}\mathbf{f}\cdot\mathbf{v}\,\mathrm{d}\Omega
-
\int_{\Gamma_N}\mathbf{h}\cdot\mathbf{v}\,\mathrm{d}\Gamma
\\
&-
\int_{\Gamma_D}
\boldsymbol{\sigma}(\rho,\tilde{\mathbf{u}})\mathbf{n}
\cdot\mathbf{v}\,\mathrm{d}\Gamma
-
\int_{\Gamma_D}
\boldsymbol{\sigma}(\rho,\tilde{\mathbf{v}})\mathbf{n}
\cdot(\mathbf{u}-\mathbf{g})\,\mathrm{d}\Gamma
\\
&+
\int_{\Gamma_D}
\frac{C_{\mathrm{pen}}^u}{h}
(\mathbf{u}-\mathbf{g})\cdot\mathbf{v}\,\mathrm{d}\Gamma,
\end{aligned}
\end{equation}
where the stabilized quantities reduce to their original forms when stabilization is not needed. 

In the implementation, the test function $\phi$ (or the trial function $\rho$) appears as a basis function defined elementwise. A subtle aspect regarding the indices of involved basis functions needs additional care. Let $\mathcal{I}_b$ and $\mathcal{I}_g$ denote the index sets of basis functions defined on $K_b$ and its $K_g$, respectively. Obviously, we have $\mathcal{I}_b \neq \mathcal{I}_g$. In evaluating the stabilized flux terms ($\nabla\tilde{\rho}\cdot \mathbf{n}$), basis functions from $\mathcal{I}_b$ will not be used. Note that they will still present in other bulk terms, so they are active degrees of freedom (DOF). On the other hand, basis functions from $\mathcal{I}_g$ will come into play in the stabilized formulation. Therefore, the pattern of the stiffness matrix will be changed slightly around the boundaries/interfaces.

In practice, the coupled problem is solved in dimensionless
form. The corresponding derivation is presented in
\ref{app:dimensionless}.

\section{Treatment of non-design domains}
\label{sec:ibcm_nondesign}

Certain important geometric features, such as bolt holes and heat-transfer interfaces, are not allowed to be modified during TO. Such regions are called \emph{non-design domains} and they are always fully occupied with material \cite{Burman_NondesignCutTO}. As a result, coupling of design and non-design domains is necessary when evaluating structural responses.

We find that the immersed boundary-conformal method (IBCM) \cite{Wei_IBCM} serves as a perfect candidate in this regard. IBCM creates a layer of mesh conformal to an important geometric boundary/interface, embeds it to a background Cartesian grid, and then weakly couples the two independent regions (the conformal layer and the cut background mesh) as a whole through Nitsche's formulation. 

In the context of TO that involves non-design domains, the conformal layer will be responsible for representing the non-design domain $\Omega^{\text{nd}}$, whereas the cut background mesh will be used to handle the design domain $\Omega^{\text{d}}$. However, their interface, $\Gamma_I:=\Omega^{\text{nd}}\cap\Omega^{\text{d}}$, needs a unified treatment to account for different material cases on the side of $\Omega^{\text{d}}$, including solid (i.e., pure material), intermediate material states, and void. Note that the other side of $\Gamma_I$ always has solid material because it is $\Omega^{\text{nd}}$.

\begin{figure}[!t]
    \centering
    \includegraphics[width=0.5\textwidth]{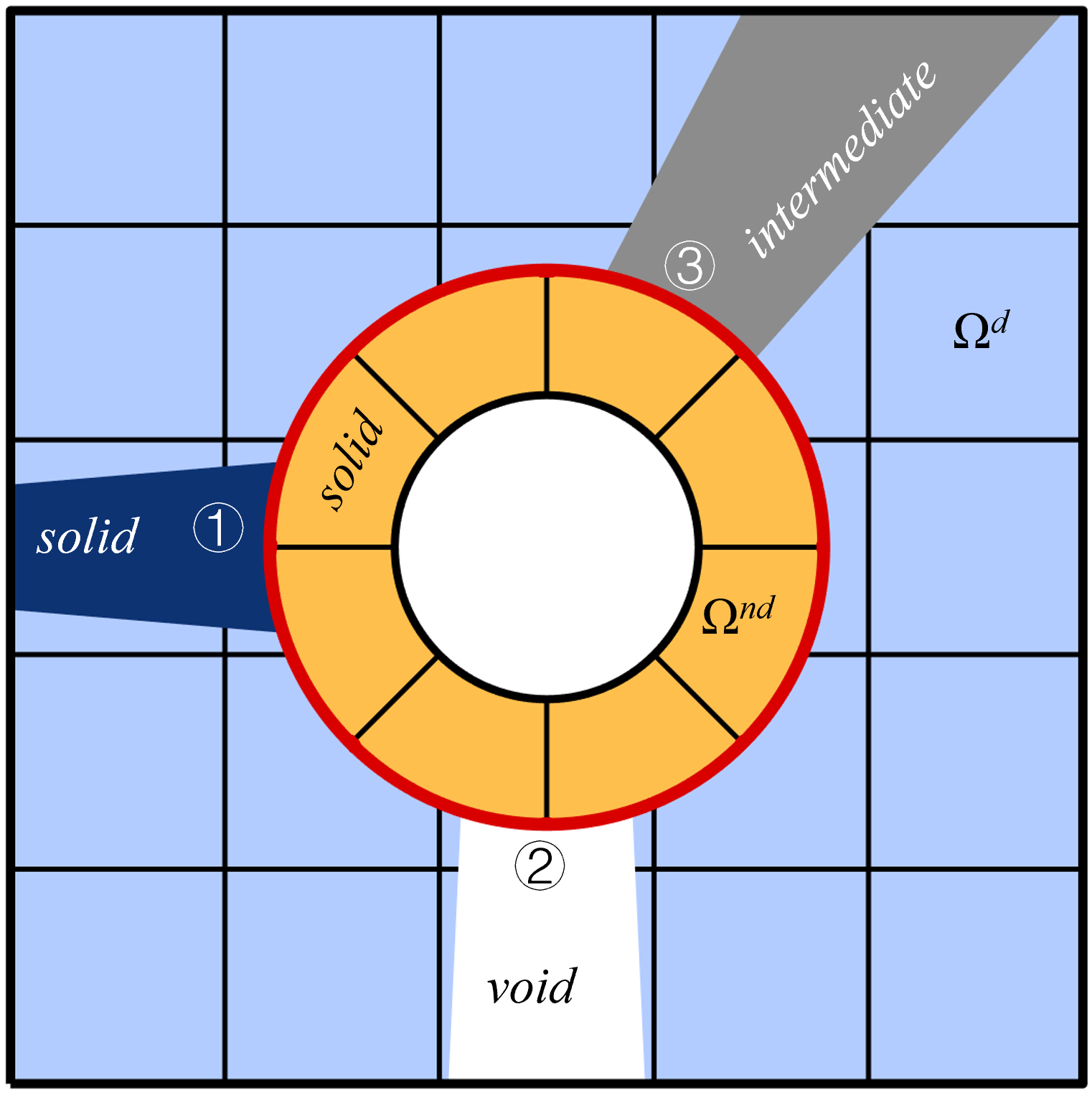}
    \caption{Schematic of the coupling between the design domain (blue) and the
    non-design domain (orange) via IBCM. For the interface (the red circle), while the non-design-domain side always has solid material, the design-domain side may have different material states: solid, void, and intermediate.}
    \label{fig:ibcm_nondesign_regimes}
\end{figure}

For example, in Fig.~\ref{fig:ibcm_nondesign_regimes}, the important feature that needs to be preserved is the hole in the center. IBCM first creates a conformal layer (the orange region) to the hole to serve as the non-design domain. Next, this layer is embedded into a regular Cartesian grid. As a result, part of the grid is trimmed away and the remaining part (the blue region) will be coupled with the conformal layer across their interface (the red circle). Note that for the interface, the design domain side may have different cases in terms of the material distribution: solid (pure material), intermediate, and void. We discuss their unified treatment next.

We still need to simultaneously solve the generalized Cahn--Hilliard problem and the equilibrium equation, but the computational domains are different for the two problems. The former is solved only in $\Omega^{\text{d}}$ with the boundary condition $\nabla\rho\cdot \bold{n}=0$ imposed on $\partial\Omega^{\text{d}}$. Note that $\Gamma_I \subseteq \partial\Omega^{\text{d}}$. The solution procedure follows exactly the same as in Section~\ref{sec:trimmed_topology_optimization}. On the other hand, the latter is solved in $\Omega^{\text{d}} \cup \Omega^{\text{nd}}$ with transmission conditions imposed on $\Gamma_I$. More specifically, the equilibrium equations are given as
\begin{equation}
\left\{
\begin{aligned}
-\nabla\cdot
\boldsymbol{\sigma}(\rho,\mathbf{u}^{\mathrm d})
&=\mathbf{f}
&&\text{in }\Omega^{\mathrm d},
\\
-\nabla\cdot
\boldsymbol{\sigma}(\mathbf{u}^{\mathrm{nd}})
&=\mathbf{f}
&&\text{in }\Omega^{\mathrm{nd}},
\\
\mathbf{u}^{\mathrm d}-\mathbf{u}^{\mathrm{nd}}
&=\mathbf{0}
&&\text{on }\Gamma_I,
\\
\boldsymbol{\sigma}(\mathbf{u}^{\mathrm{nd}})
\mathbf{n}^{\mathrm{nd}}
-
\boldsymbol{\sigma}(\rho,\mathbf{u}^{\mathrm d})
\mathbf{n}^{\mathrm {nd}}
&=\mathbf{0}
&&\text{on }\Gamma_I,
\\
\mathbf{u}
&=\mathbf{g}
&&\text{on }\Gamma_D,
\\
\boldsymbol{\sigma}(\rho,\mathbf{u})\mathbf{n}
&=\mathbf{h}
&&\text{on }\Gamma_N,
\end{aligned}
\right.
\label{eq:coupled_system_nondesign}
\end{equation}
where $\mathbf{n}^{\mathrm{nd}}$ is the outwards unit normal to $\partial\Omega^{\text{nd}}$, and $\mathbf{u}$ can be $\mathbf{u}^{\mathrm d}$ or $\mathbf{u}^{\mathrm {nd}}$. The Dirichlet and Neumann boundaries $\Gamma_D$ and $\Gamma_N$ may lie on $\partial\Omega^{\mathrm d}$ or $\partial\Omega^{\mathrm{nd}}$. Note that we always have $\rho=1$ in $\Omega^{\text{nd}}$, so the dependency of $\rho$ in $\boldsymbol{\sigma}$ is removed in this region.

We focus on the treatment of the two transmission conditions on $\Gamma_I$, whereas the other terms are handled in the same manner as in Section~\ref{sec:trimmed_topology_optimization}. Coupling the two sides of $\Gamma_I$ is expected to behave as follows. When it is solid material ($\rho = 1$) on the side of $\Omega^{\text{d}}$, the two transmission conditions should be treated in the same classical way of IBCM. On the other hand, when it is void ($\rho=0$) on the side of $\Omega^{\text{d}}$, the displacement condition should be discarded and a traction-free boundary condition should be recovered, i.e.,
$\boldsymbol{\sigma}(\mathbf{u}^{\mathrm{nd}})
\mathbf{n}^{\mathrm{nd}}=\mathbf{0}$. When it is intermediate ($0 < \rho < 1$), the imposition of the two transmission conditions should be penalized towards either of the two extremes ($\rho = 0$ or 1).

To treat these different cases in a unified manner, we propose the weak form of Problem~\eqref{eq:coupled_system_nondesign} as follows,

\begin{align}
R_{\mathbf{u}}(\rho,\mathbf{u})(\mathbf{v})
:={}&
\int_{\Omega^{\mathrm d}}
\boldsymbol{\sigma}(\rho,\mathbf{u}^{\mathrm d})
:
\boldsymbol{\varepsilon}(\mathbf{v}^{\mathrm d})\,\mathrm{d}\Omega
+
\int_{\Omega^{\mathrm{nd}}}
\boldsymbol{\sigma}(\mathbf{u}^{\mathrm{nd}})
:
\boldsymbol{\varepsilon}(\mathbf{v}^{\mathrm{nd}})\,\mathrm{d}\Omega
\nonumber\\
&-
\int_{\Omega^{\mathrm d}}
\mathbf{f}\cdot\mathbf{v}^{\mathrm d}\,\mathrm{d}\Omega
-
\int_{\Omega^{\mathrm{nd}}}
\mathbf{f}\cdot\mathbf{v}^{\mathrm{nd}}\,\mathrm{d}\Omega
-
\int_{\Gamma_N}
\mathbf{h}\cdot\mathbf{v}\,\mathrm{d}\Gamma
\nonumber\\
&-
\int_{\Gamma_D}
\boldsymbol{\sigma}(\rho,\mathbf{u})\mathbf{n}
\cdot\mathbf{v}\,\mathrm{d}\Gamma
-
\int_{\Gamma_D}
\boldsymbol{\sigma}(\rho,\mathbf{v})\mathbf{n}
\cdot(\mathbf{u}-\mathbf{g})\,\mathrm{d}\Gamma
\nonumber\\
&+
\int_{\Gamma_D}
\frac{C_{\mathrm{pen}}^u}{h}
(\mathbf{u}-\mathbf{g})\cdot\mathbf{v}\,\mathrm{d}\Gamma
\nonumber\\
&-
\int_{\Gamma_I}
\left\langle
\boldsymbol{\sigma}(\rho,\mathbf{u})\mathbf{n}^{\mathrm{nd}}
\right\rangle
\cdot[\mathbf{v}]\,\mathrm{d}\Gamma
\nonumber\\
&-
\int_{\Gamma_I}
\left\langle
\boldsymbol{\sigma}(\rho,\mathbf{v})\mathbf{n}^{\mathrm{nd}}
\right\rangle
\cdot[\mathbf{u}]\,\mathrm{d}\Gamma
\nonumber\\
&+
\beta
\int_{\Gamma_I}
g(\rho)
\left(
\frac{1}{h^\mathrm{d}}+\frac{1}{h^{\mathrm{nd}}}
\right)
[\mathbf{u}]\cdot[\mathbf{v}]\,\mathrm{d}\Gamma,
\label{eq:ibcm_nondesign_residual}
\end{align}
where we have adopted Nitsche's method to impose both the boundary conditions and the transmission conditions. The boundary conditions follow the same as in Section~\ref{sec:trimmed_topology_optimization}, whereas we introduce a scaling factor, $g(\rho)$ from Eq.~\eqref{eq:material_interpolation_function}, for the transmission conditions. The average and jump terms are defined as
\begin{equation}
\begin{aligned}
\left\langle
\boldsymbol{\sigma}(\rho,\mathbf{u})
\mathbf{n}^{\mathrm{nd}}
\right\rangle
&:=
g(\rho)\,
\boldsymbol{\sigma}(\mathbf{u}^{\mathrm{nd}})
\mathbf{n}^{\mathrm{nd}},
\\
[\mathbf{u}]
&:=
\mathbf{u}^{\mathrm{nd}}-\mathbf{u}^{\mathrm d}.
\end{aligned}
\label{eq:ibcm_average_jump}
\end{equation}
\noindent
Note that we have adopted one-sided flux in the definition of the average term. This choice will be commented on later. The trimming-independent element size parameters $h^{\mathrm d}$ and $h^{\mathrm{nd}}$ correspond to the meshes of $\Omega^{\mathrm d}$ and $\Omega^{\mathrm{nd}}$, respectively. Following \cite{Wei_IBCM}, the penalty parameter $\beta$ is taken as
\begin{equation}
\beta
=
6\left(p^{\max}\right)^2
\times
8\left(3\lambda_{\mathrm{L}}+2\mu_{\mathrm{L}}\right),
\label{eq:ibcm_penalty_beta}
\end{equation}
where $p^{\max}$ is the maximum degree of involved splines, and
$\lambda_{\mathrm{L}}$ and $\mu_{\mathrm{L}}$ are the Lam\'e constants of the solid material.

It is the scaling factor $g(\rho)$ that enables a unified treatment of the transmission conditions in the different material cases of $\Omega^{\mathrm d}$. We can observe that when $\rho=1$, we have $g(\rho)=1$ and Eq.~\eqref{eq:ibcm_nondesign_residual} coincides with the classical formulation of IBCM. When $\rho=0$, we have $g(\rho)=0$. Thus, all the interface terms vanish, so the displacement condition is discarded and the traction-free condition is achieved. Moreover, when $0<\rho<1$, the degree of penalization is controlled by $g(\rho)$ in a smooth manner. Note that in theory, we can make $\rho_{\min}=0$  in the definition of $g(\rho)$, but in practice, we choose a small positive value for $\rho_{\min}$.

\begin{remark}
The CutFEM-based TO \cite{Burman_NondesignCutTO} treats non-design domains in a similar setting. The non-design domain is also represented by a conformal mesh, and it is coupled with the cut background mesh of the design domain through Nitsche's method. The main difference lies in how the material state enters the transmission conditions. In \cite{Burman_NondesignCutTO}, the interface terms are weighted by the densities on both sides, so that either side may become void. Such generality is not needed here, as $\Omega^{\mathrm{nd}}$ is always solid. It thus suffices to scale the interface terms by $g(\rho)$ from the side of $\Omega^{\mathrm{d}}$. Moreover, the flux can be taken only from the untrimmed conformal layer, so that no stabilization is needed for the transmission conditions; see the next remark. Another difference stems from the phase-field model. Since $\Gamma_I$ is also a trimmed boundary of the phase-field problem, $\nabla\rho\cdot\mathbf{n}=0$ needs to be imposed there. Such a condition does not arise in the density-based TO of \cite{Burman_NondesignCutTO}.
\end{remark}

\begin{remark}
In general, the average term is defined as
\begin{equation}
\begin{aligned}
\left\langle
\boldsymbol{\sigma}(\rho,\mathbf{u})
\mathbf{n}^{\text{nd}}
\right\rangle
:={}&
\omega\,
\boldsymbol{\sigma}(\rho,\mathbf{u}^{\text{d}})
\mathbf{n}^{\text{d}}
+
(1-\omega)\,g(\rho)\,
\boldsymbol{\sigma}(\mathbf{u}^{\text{nd}})
\mathbf{n}^{\text{nd}},
\qquad \omega\in[0,1].
\end{aligned}
\end{equation}
Our choice here sets $\omega=0$, which yields a one-sided flux only from the side of $\Omega^{\mathrm{nd}}$. Recall that $\Omega^{\mathrm{nd}}$ is represented as a conformal layer and it does not involve trimming. Therefore, the flux from this side does not need any stabilization. In contrast, when $\omega\neq 0$, we end up with a two-sided flux, which is the case of \cite{Burman_NondesignCutTO}. The flux from $\Omega^{\mathrm{d}}$, which is represented by a cut background mesh, will enter the formulation and thus needs stabilization. To have a simpler formulation, we adopt the one-sided flux. We will show that the choice does not affect the TO results.
\end{remark}

\begin{remark}
\label{rem:dimensionless}
In practice, the coupled problem is solved in a dimensionless form. The derivation follows a standard way \cite{Dede_PhaseFieldITO,Sheng_ThinWalledComplex} and it is given in~\ref{app:dimensionless}. In particular, the parameters $\lambda$ and $\gamma$ in Eq.~\eqref{eq:gch_quantities_review} are specified through the dimensionless parameters $\bar{\lambda}$ and $\bar{\gamma}$ as
\begin{equation*}
\lambda=\bar{\lambda}h^2,
\qquad
\gamma=\bar{\gamma}\gamma_E,
\end{equation*}
where $h$ is the size of the background mesh and $\gamma_E$ balances the strain energy against the bulk and interface energies of the initial design; see its definition in~\ref{app:dimensionless}. $J^{\star}$, $J_E^{\star}$, $J_B^{\star}$, and $J_I^{\star}$ are the dimensionless forms of $J$, $J_E$, $J_B$, and $J_I$, respectively.
\end{remark}

\section{Numerical examples}
\label{sec:numerical_examples}

This section presents several numerical tests to demonstrate the
stability, accuracy, and flexibility of the proposed method. The implementation is developed in C++.
Sparse matrix assembly and linear algebra operations are carried out with Eigen
\cite{Guennebaud_EigenLibrary}; shared-memory parallelism is implemented using
OpenMP \cite{OpenMPARB_OpenMPSpecification}; and sparse direct-solver kernels are
accelerated with Intel oneAPI Math Kernel Library
\cite{Intel_oneMKLDeveloperReference}.

\subsection{Biharmonic problem in 3D}
\label{subsec:biharmonic_3d}

We first study minimal stabilization in the biharmonic problem, which is a linear fourth-order problem and serves as the foundation of the nonlinear generalized Cahn--Hilliard problem. In particular, we consider a unit cube domain with homogeneous boundary conditions,
\begin{equation}
\left\{
\begin{aligned}
\Delta^2 u &= q
&& \text{in }\Omega = (0,1)^3,\\
u &= 0
&& \text{on }\partial\Omega,\\
\Delta u &= 0
&& \text{on }\partial\Omega.
\end{aligned}
\right.
\label{eq:biharmonic_strong_form}
\end{equation}
The source term and the corresponding exact solution are
\begin{equation}
\begin{aligned}
q(x,y,z)
&=\sin(\pi x)\sin(\pi y)\sin(\pi z),
\\
u_{\mathrm{ex}}(x,y,z)
&=\frac{1}{9\pi^4}
\sin(\pi x)\sin(\pi y)\sin(\pi z).
\end{aligned}
\label{eq:biharmonic_source_exact}
\end{equation}

\begin{figure}[htbp]
    \centering
    \includegraphics[width=0.46\textwidth]{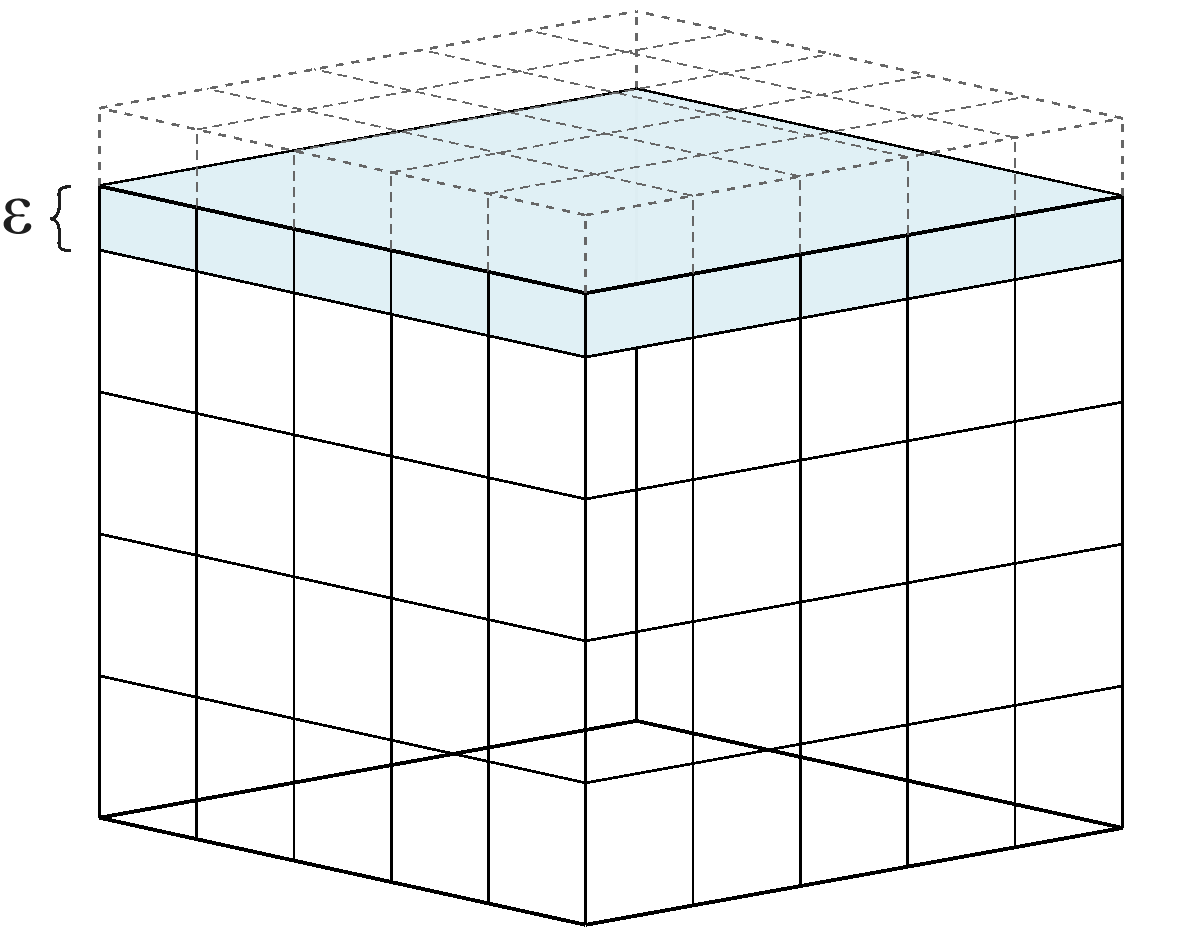}
        \caption{Biharmonic problem in a unit cube. The blue plane cuts the
    background mesh, where the dashed lines indicate the part of the mesh
    that is trimmed away, and parameter $\varepsilon$ is used to control
    the size of cut elements.}
    \label{fig:biharmonic_3d_mesh}
\end{figure}

The unit cube is represented by a Cartesian grid that is cut by a plane on the top; see
Fig.~\ref{fig:biharmonic_3d_mesh}. The parameter $\varepsilon$ is used to control the size of cut elements. 

Let $V_h\subset H^2(\Omega)$ denote the spline space. Nitsche's method is used to impose the boundary conditions.
The weak form of Eq.~\eqref{eq:biharmonic_strong_form}
seeks $u_h\in V_h$ such that
\begin{equation}
a_h(u_h,v_h)=\ell_h(v_h),\qquad \forall v_h\in V_h,
\label{eq:biharmonic_weak_problem}
\end{equation}
where
\begin{equation}
\begin{aligned}
a_h(u_h,v_h)
={}&\int_\Omega\Delta u_h\,\Delta v_h\,\mathrm{d}\Omega
\\
&+\int_{\partial\Omega}
u_h\bigl(\nabla(\Delta v_h)\cdot\mathbf n\bigr)
\,\mathrm{d}\Gamma
+\int_{\partial\Omega}
v_h\bigl(\nabla(\Delta u_h)\cdot\mathbf n\bigr)
\,\mathrm{d}\Gamma
\\
&+\frac{C_{\mathrm{pen}}}{h^3}\int_{\partial\Omega}u_h\,v_h\,\mathrm{d}\Gamma,
\\
\ell_h(v_h)
={}&\int_\Omega v_h\,q\,\mathrm{d}\Omega.
\end{aligned}
\label{eq:biharmonic_weak_forms}
\end{equation}
Note that the formulation involves high-order fluxes on the boundary, $\nabla(\Delta u_h)\cdot\mathbf n$, which is stabilized according to Section~\ref{subsec:minimal_stabilization_high_order}. The penalty term scales as $h^{-3}$, and we set $C_{\mathrm{pen}}=1$.

\begin{figure}[!t]
    \centering

    \begin{subfigure}[t]{0.49\textwidth}
        \centering
        \includegraphics[width=\linewidth]{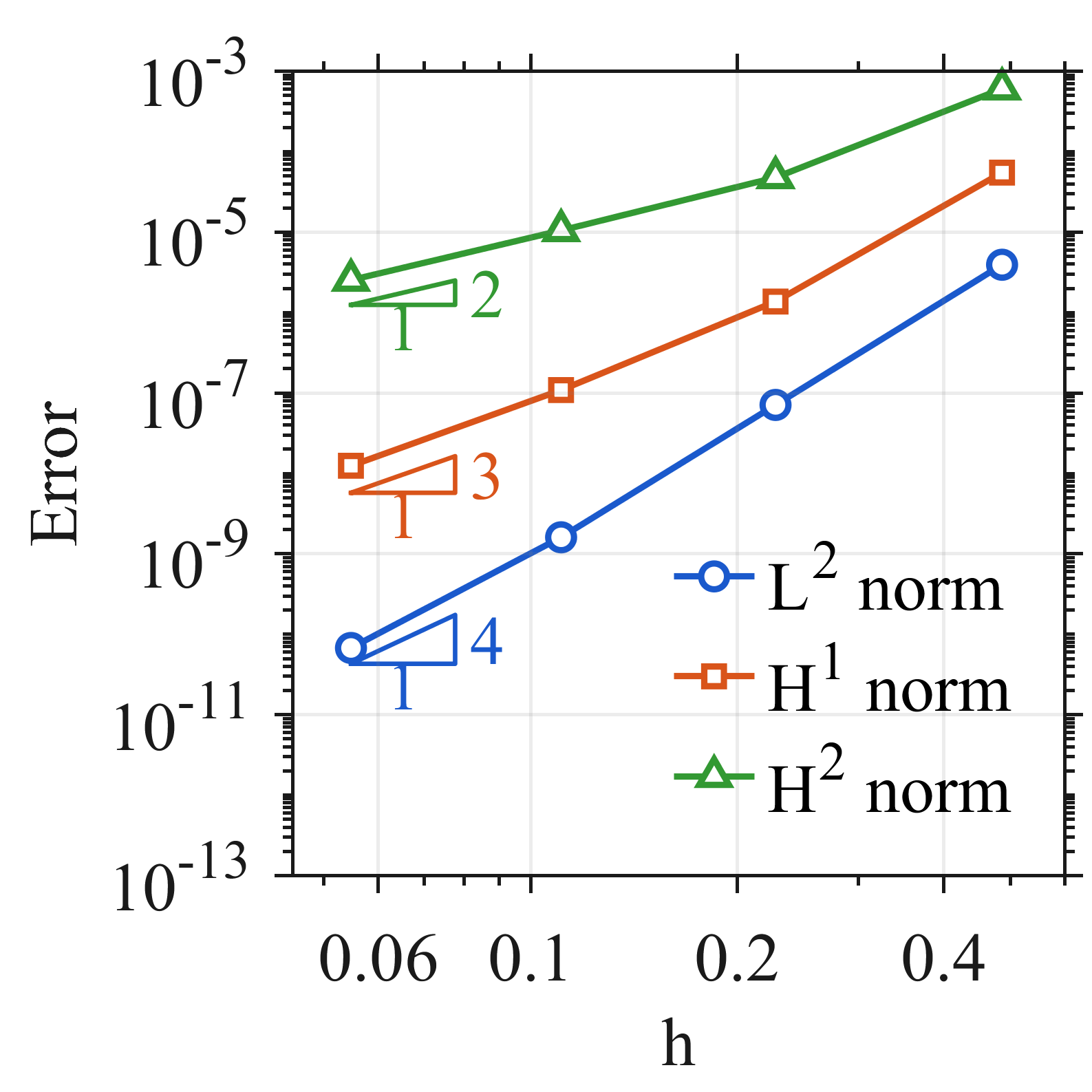}
        \caption{Tri-cubic B-splines}
        \label{fig:biharmonic3d_norm}
    \end{subfigure}\hfill%
    \begin{subfigure}[t]{0.49\textwidth}
        \centering
        \includegraphics[width=\linewidth]{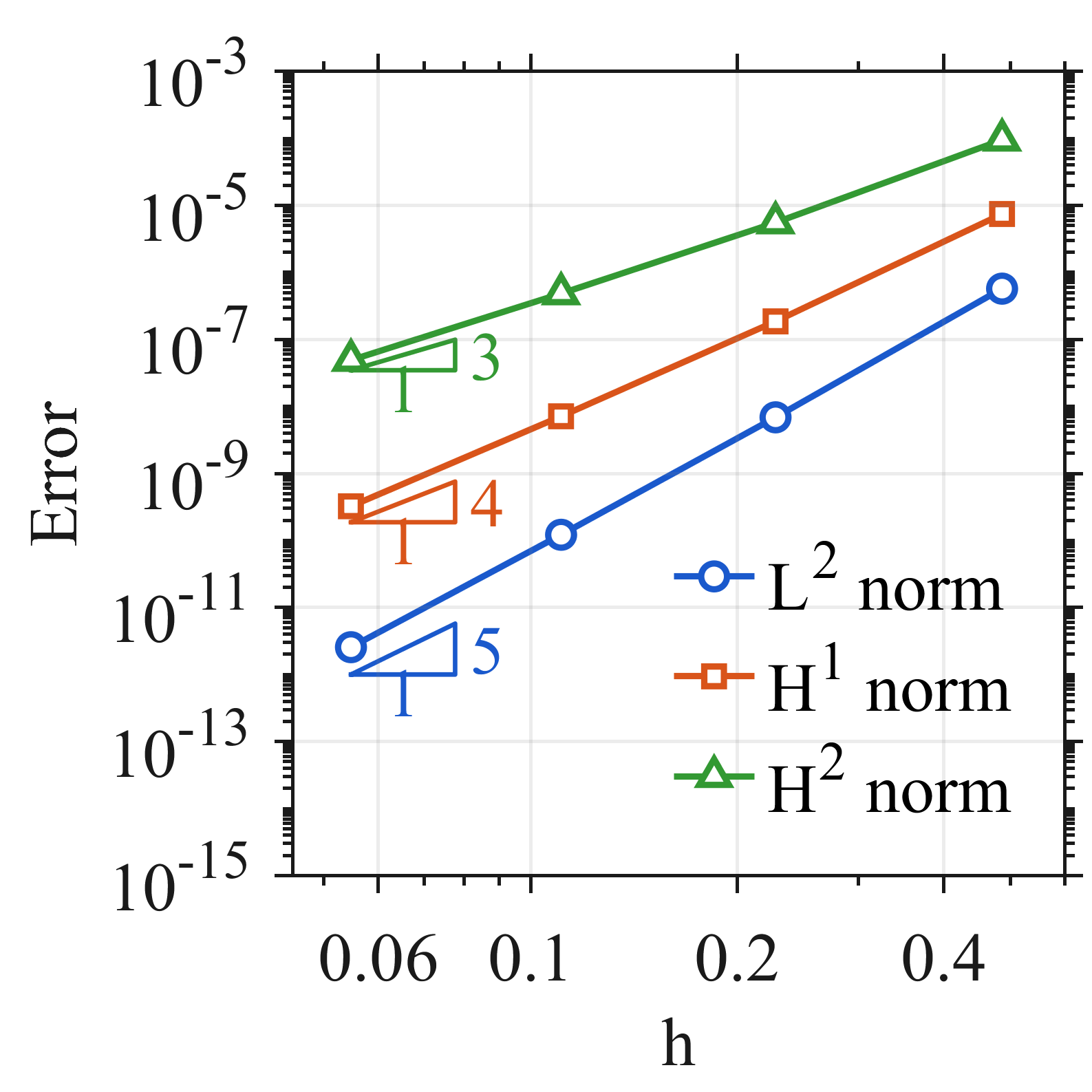}
        \caption{Tri-quartic B-splines}
        \label{fig:biharmonic3d_norm_p4}
    \end{subfigure}

    \caption{Convergence plots of errors in different norms with respect to
    the element size of the background mesh.}
    \label{fig:biharmonic_3d_results}
\end{figure}

We first check the convergence of the stabilized method. Tri-cubic and
tri-quartic B-splines are used because the third order normal derivatives are involved. We set $\varepsilon=10^{-9}$ in both cases. The $L^2$, $H^1$ and $H^2$ errors are
shown with respect to the element size ($h$) of the background mesh; see
Fig.~\ref{fig:biharmonic_3d_results}. Expected optimal
convergence rates are achieved in all cases, indicating that minimal
stabilization preserves the expected accuracy in the high-order problem
with high-order boundary fluxes.

\begin{figure}[!t]
    \centering
    \includegraphics[width=0.6\textwidth]{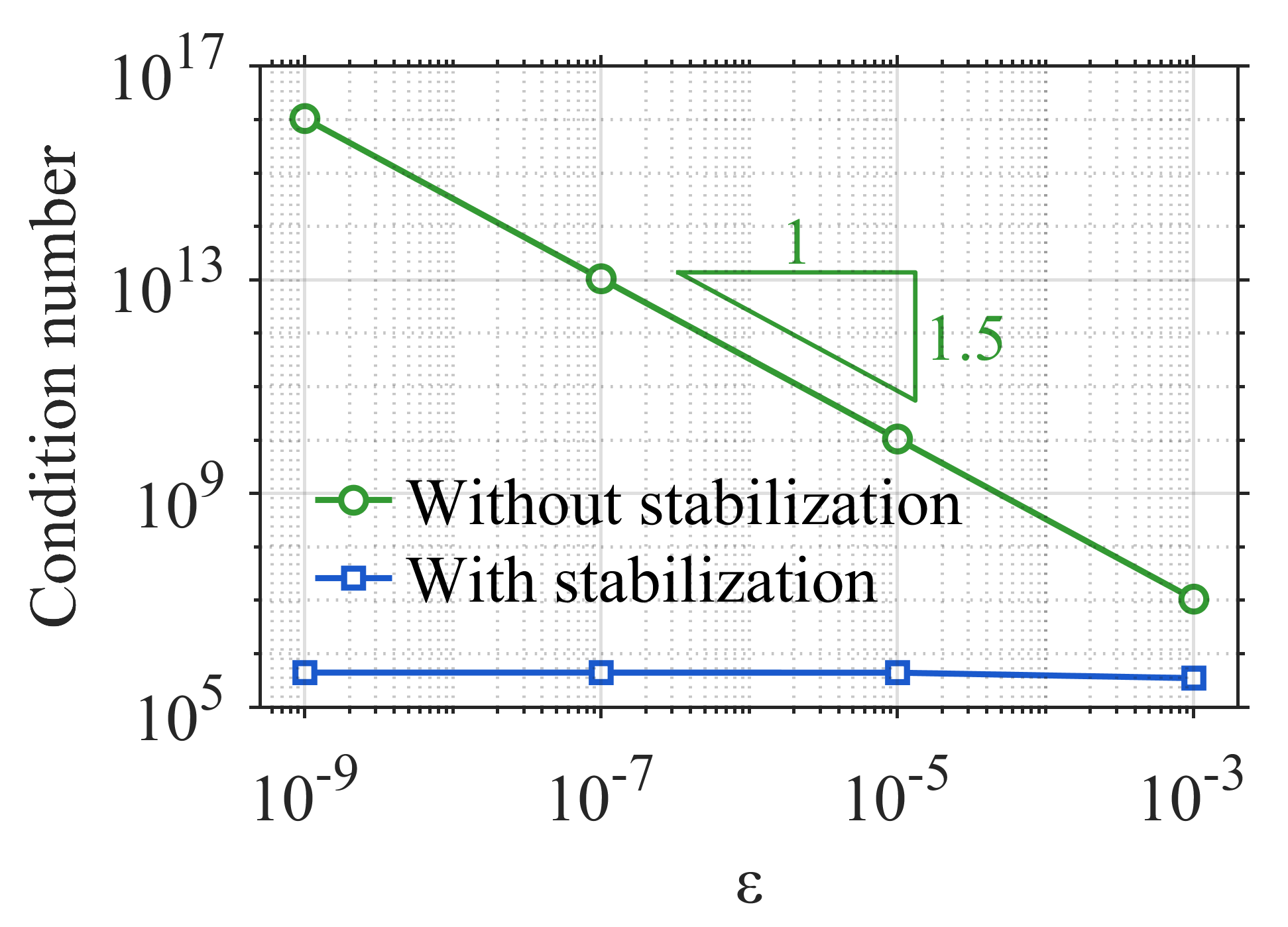}
        \caption{Conditioning of stiffness matrices with and without
    stabilization (tri-cubic B-splines).}
    \label{fig:biharmonic_3d_condition}
\end{figure}

We next study the conditioning of the stiffness matrix as cut elements
become arbitrarily small. The condition numbers are computed on the
diagonally rescaled matrices. Given stiffness matrix $\mathbf A$, its
rescaled matrix is defined as
\begin{equation}
\begin{aligned}
\widetilde{\mathbf{A}}
&=
\mathbf{D}^{-1/2}
\mathbf{A}
\mathbf{D}^{-1/2},
\qquad
\mathbf{D}
=
\operatorname{diag}(\mathbf{A}),
\end{aligned}
\label{eq:diagonally_rescaled_matrix}
\end{equation}
where $\mathbf D$ is the diagonal matrix of $\mathbf A$. In every test, $h$
is fixed and $\varepsilon$ varies from $10^{-3}$ to $10^{-9}$. The results
are shown in Fig.~\ref{fig:biharmonic_3d_condition}. We observe that the
condition number explodes without stabilization, whereas it remains nearly
constant with minimal stabilization. This result demonstrates the
effectiveness of minimal stabilization in the high-order problem.

\begin{figure}[!t]
    \centering
    
    \begin{subfigure}[t]{0.495\textwidth}
        \centering
        \vspace{0pt}
        \parbox[t][0.534\linewidth][c]{\linewidth}{\centering
        \includegraphics[width=\linewidth]{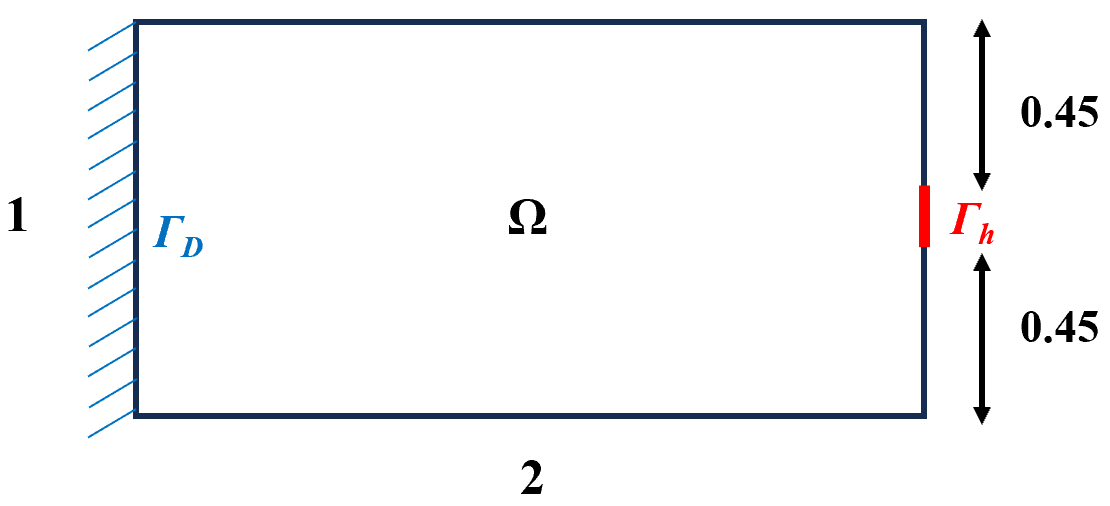}}
        \caption{Problem setup}
        \label{fig:cantilever_problem}
    \end{subfigure}\hfill%
    \begin{subfigure}[t]{0.495\textwidth}
        \centering
        \vspace{0pt}
        \includegraphics[width=\linewidth]{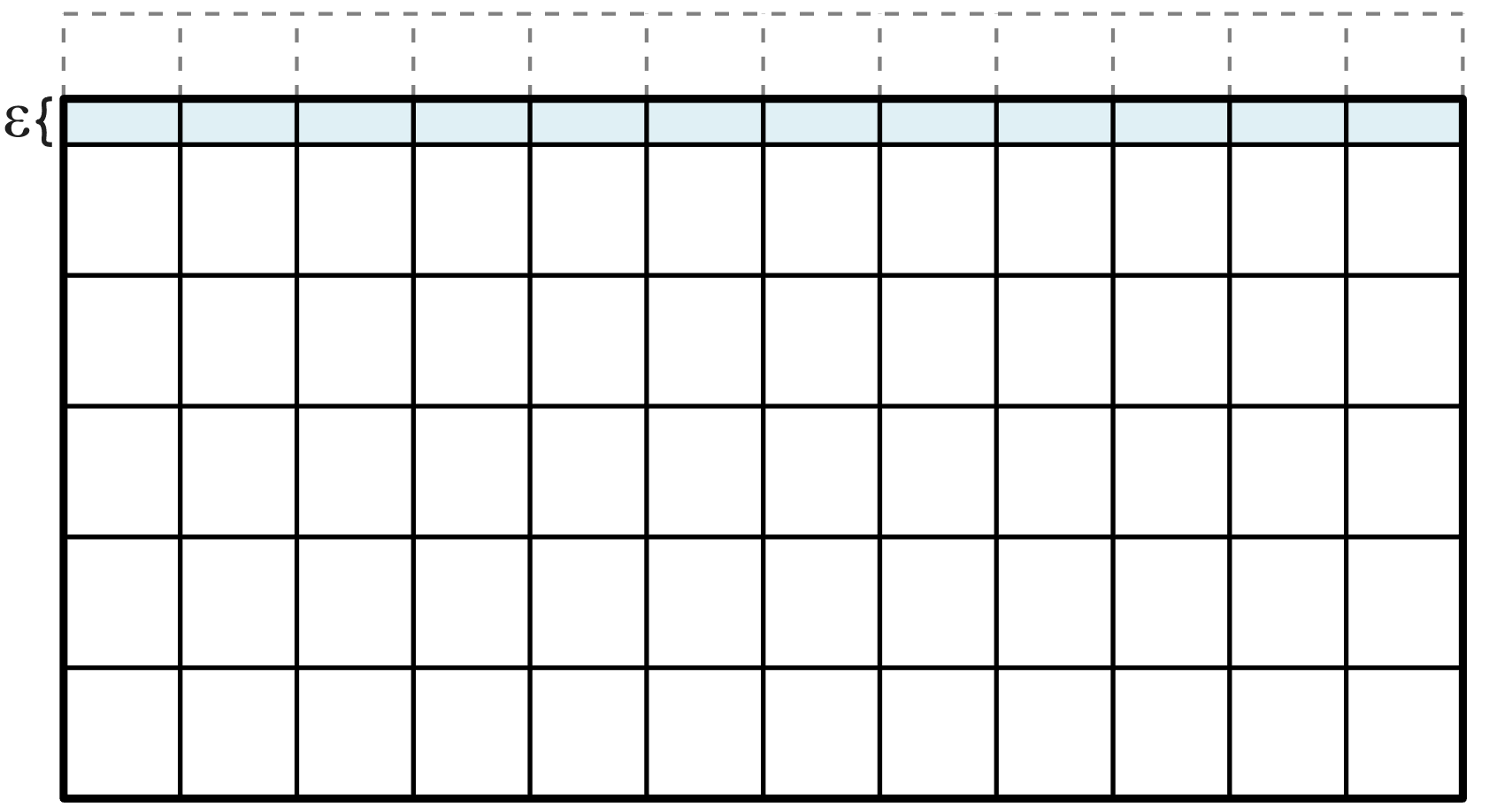}
        \caption{Background mesh and cut elements (shaded blue)}
        \label{fig:cantilever_epsilon_mesh_global}
    \end{subfigure}

    \par\medskip

    \begin{subfigure}[t]{0.495\textwidth}
        \centering
        \includegraphics[width=\linewidth]{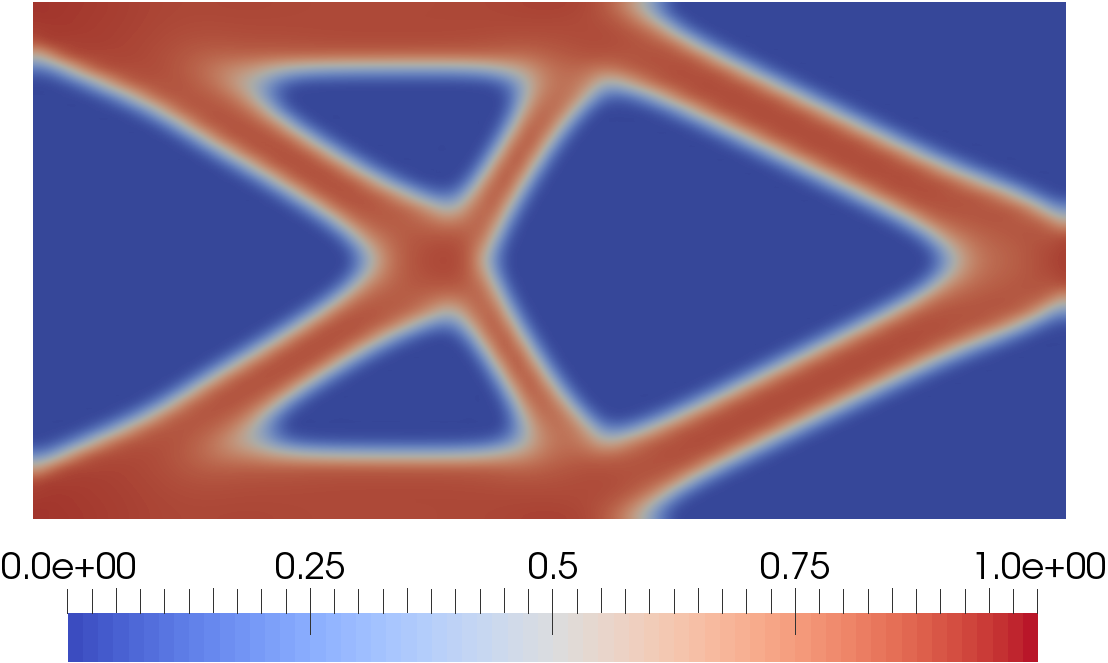}
        \caption{Optimized structure at minimum $J_E^{\star}$ with minimal stabilization}
        \label{fig:cantilever_epsilon_result}
    \end{subfigure}\hfill%
    \begin{subfigure}[t]{0.495\textwidth}
        \centering
        \includegraphics[width=\linewidth]{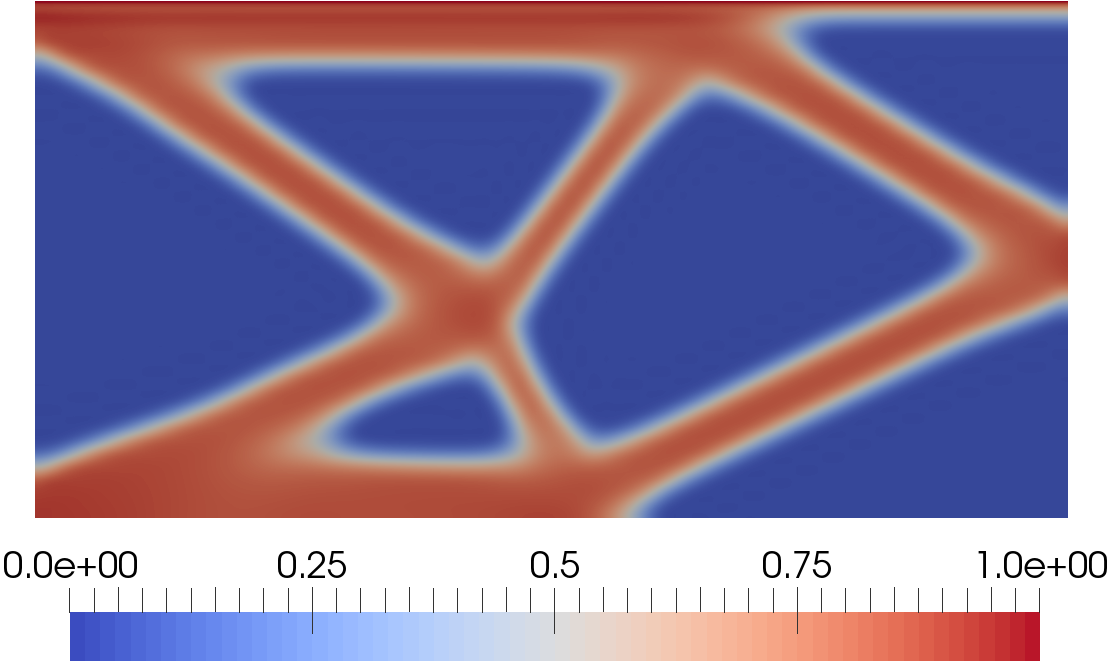}
        \caption{Phase field without stabilization}
        \label{fig:cantilever_epsilon_result_nostable}
    \end{subfigure}

    \caption{TO on a cantilever beam, whose representation involves extremely small cut
    elements ($\varepsilon=10^{-9}$). In (a), $\Omega$ is the design domain, the left
    boundary ($\Gamma_D$) is fixed, and a downward traction is applied on $\Gamma_h$.}
    \label{fig:cantilever_epsilon_results}
\end{figure}

\subsection{Cantilever beam}
\label{subsec:cantilever_stabilization}

We next study the phase-field-based TO on a cantilever beam when
cut elements become arbitrarily small. The nonlinear, generalized Cahn--Hilliard
problem is solved. The problem setup is shown in
Fig.~\ref{fig:cantilever_epsilon_results}(\subref{fig:cantilever_problem}). The design domain is a $2\,\mathrm{m}\times1\,\mathrm{m}$ rectangle. The left edge $\Gamma_D$ is fixed. A uniform downward traction of magnitude
$2.0\times10^8$ N/m is applied to the red segment $\Gamma_h$ on the right edge.
The plane-strain assumption is adopted with Young's modulus $E=2.0\times10^{11}$ Pa
and Poisson's ratio $\nu=1/3$. Bicubic B-splines are used. The background
mesh resolution is fixed
($h=0.035$ m) as the focus here is the effect of small cut elements on TO, so the
variable is mainly $\varepsilon$. The prescribed material fraction is $0.35$.
The phase-field parameters are set to $\bar{\lambda}=2$ and $\bar{\gamma}=50$, and the penalty parameter is $C_{\mathrm{pen}}^{\rho}=1$.

\begin{figure}[!t]
    \centering
    \begin{subfigure}[t]{0.48\textwidth}
        \centering
        \includegraphics[width=\linewidth,height=0.7\linewidth,keepaspectratio]{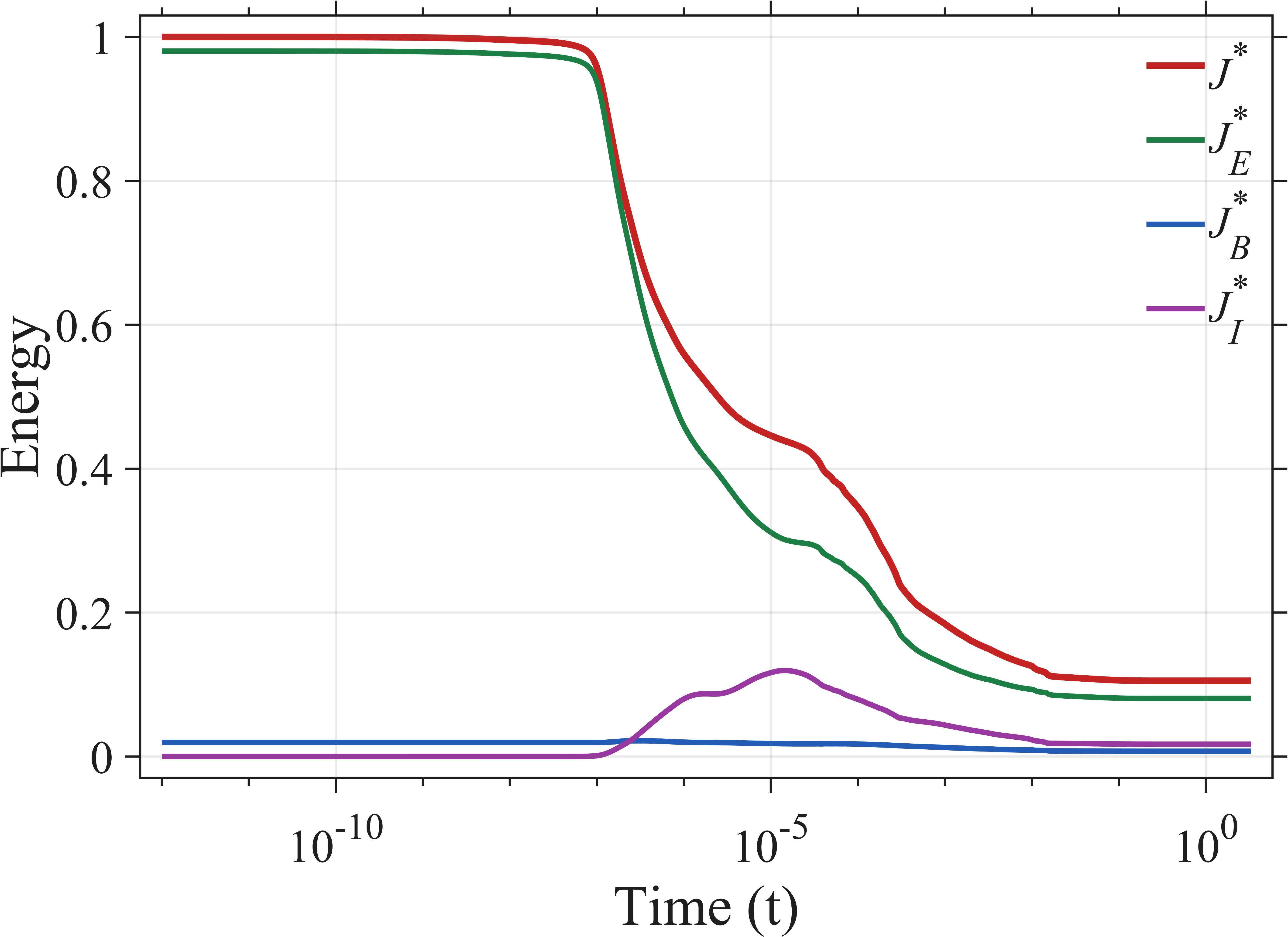}
        \caption{Normalized energy with minimal stabilization}
        \label{fig:cantilever_energy_history}
    \end{subfigure}\hfill%
    \begin{subfigure}[t]{0.48\textwidth}
        \centering
        \includegraphics[width=\linewidth,height=0.7\linewidth,keepaspectratio]{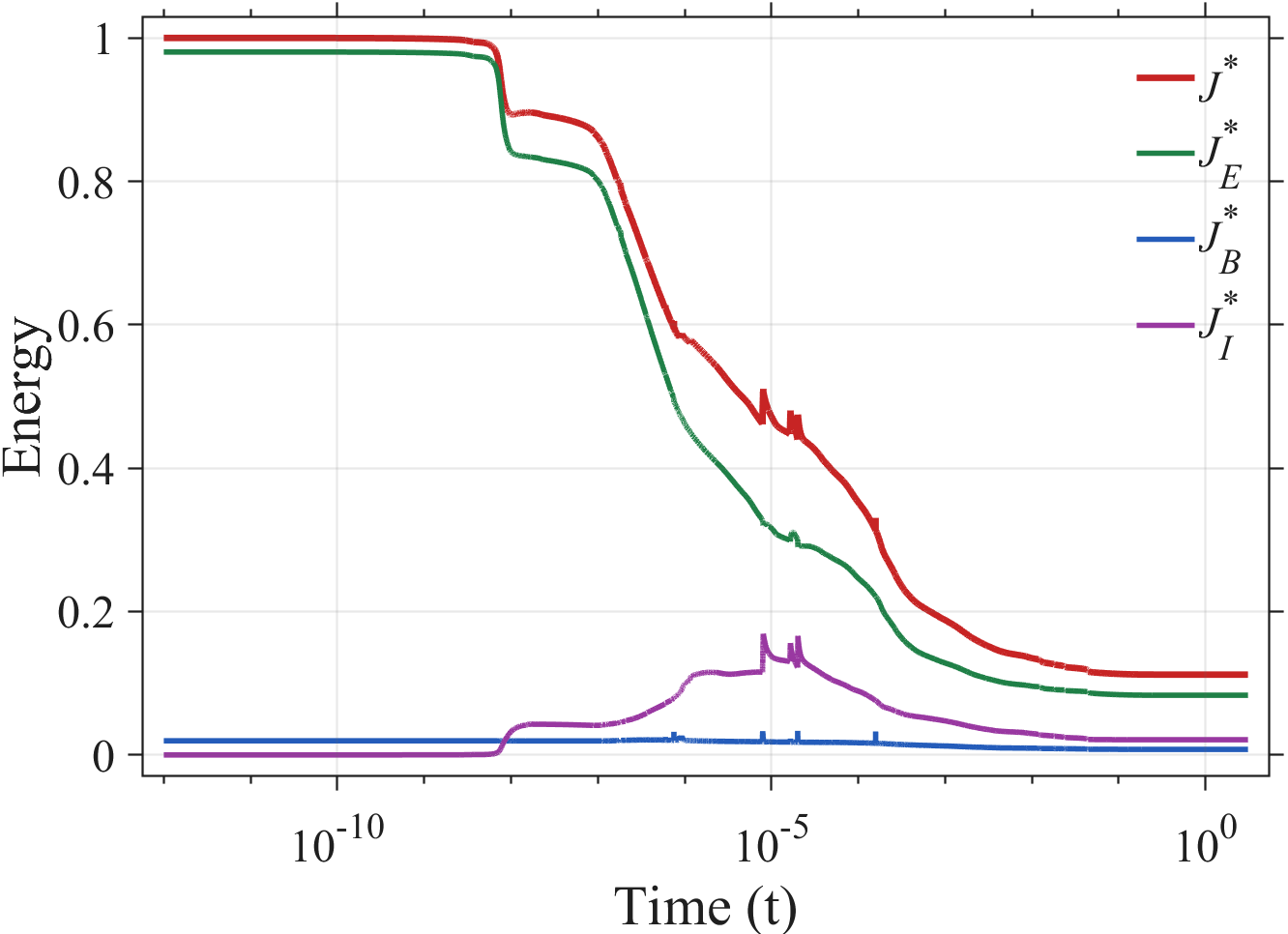}
        \caption{Normalized energy without stabilization}
        \label{fig:cantilever_energy_history_nostable}
    \end{subfigure}

    \par\medskip

    \begin{subfigure}[t]{0.48\textwidth}
        \centering
        \includegraphics[width=\linewidth,height=0.7\linewidth,keepaspectratio]{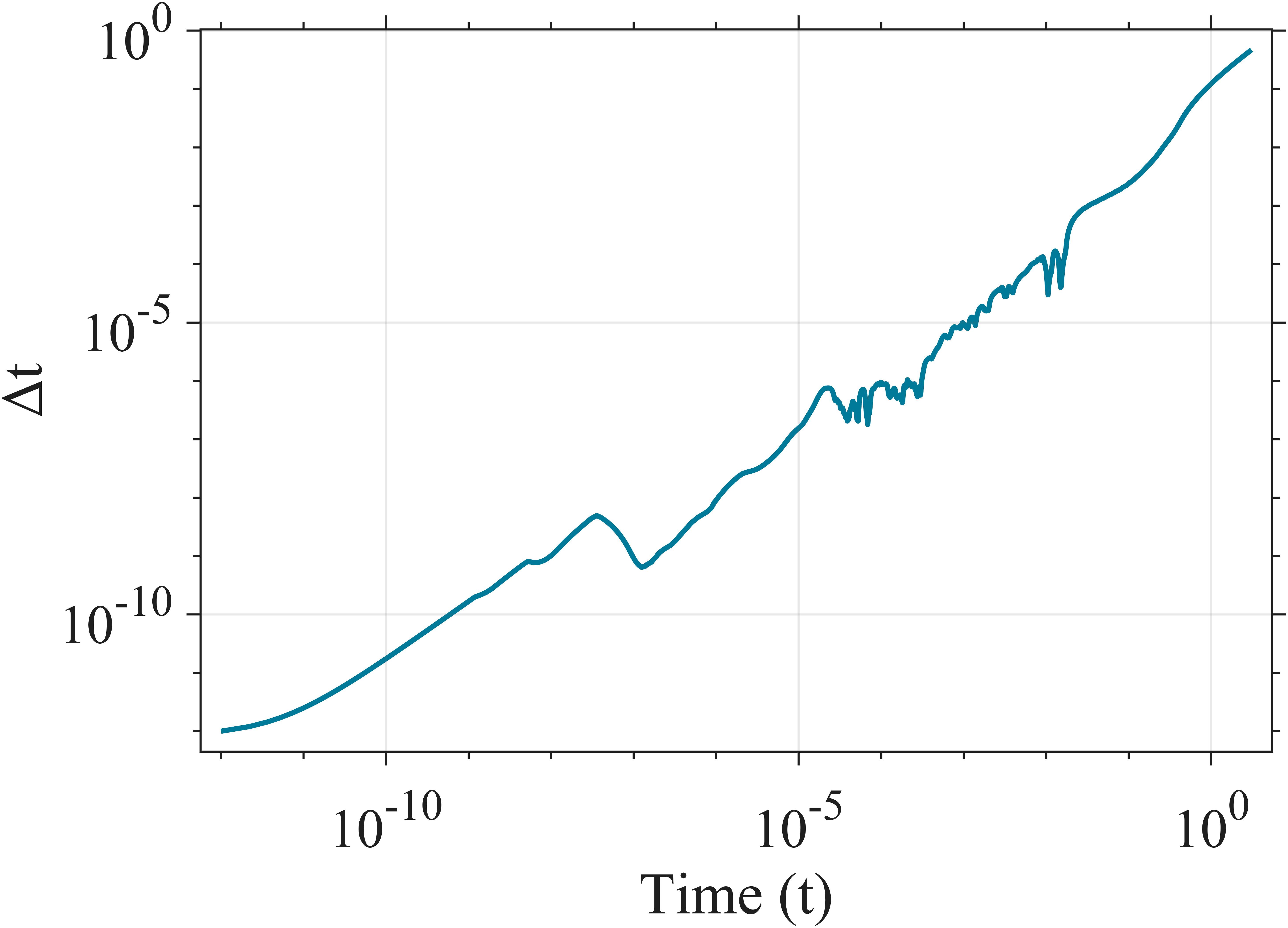}
        \caption{$\Delta t$ with minimal stabilization}
        \label{fig:cantilever_dt_stable}
    \end{subfigure}\hfill%
    \begin{subfigure}[t]{0.48\textwidth}
        \centering
        \includegraphics[width=\linewidth,height=0.7\linewidth,keepaspectratio]{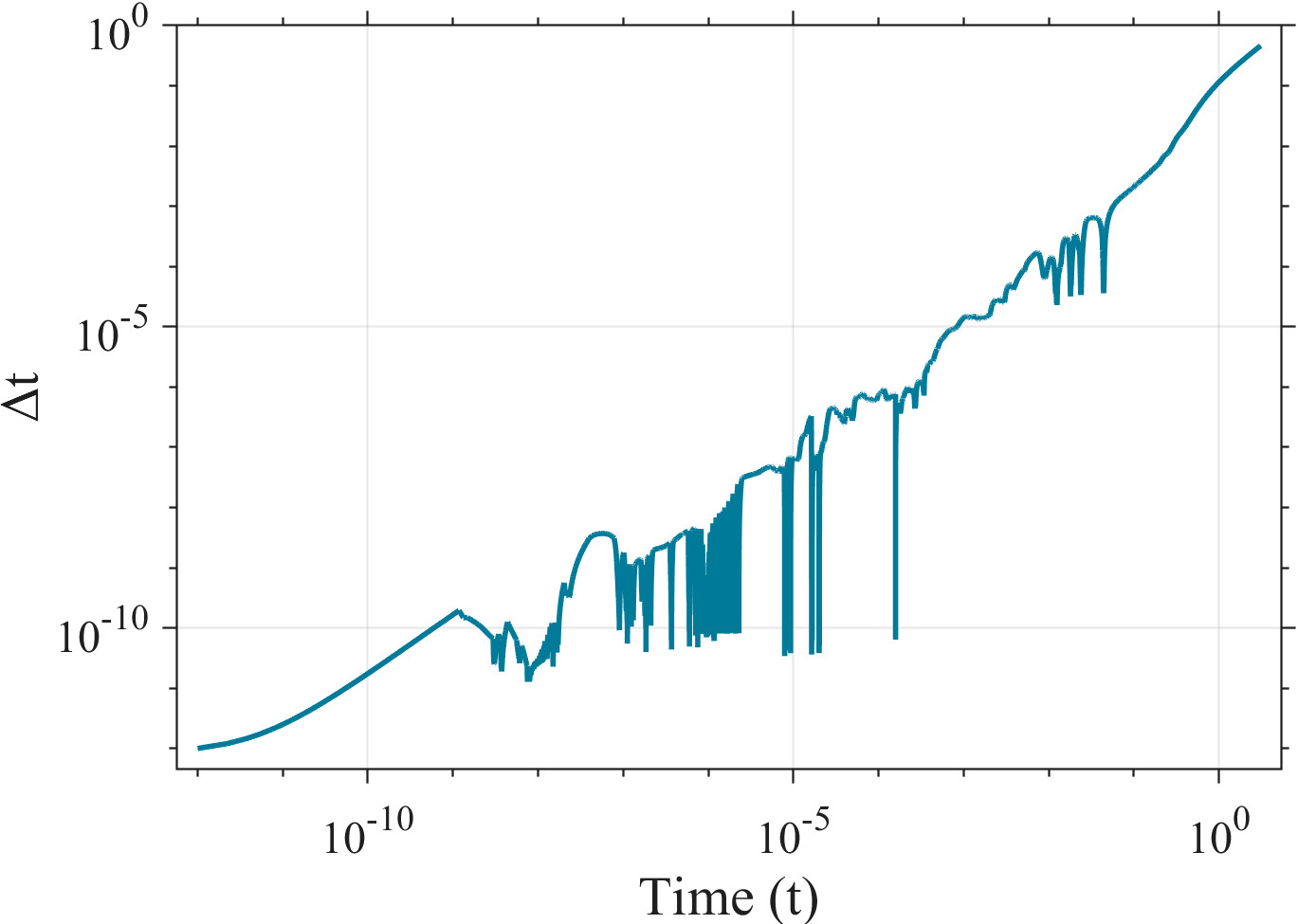}
        \caption{$\Delta t$ without stabilization}
        \label{fig:cantilever_dt_nostable}
    \end{subfigure}
    \caption{Normalized energy evolution (a, b) and adaptive time steps (c, d) of
    the phase-field-based TO on the cantilever beam with extremely small cut
    elements ($\varepsilon=10^{-9}$).}
    \label{fig:cantilever_energy_histories}
\end{figure}

We first consider a cut mesh shown in
Fig.~\ref{fig:cantilever_epsilon_results}(\subref{fig:cantilever_epsilon_mesh_global}),
where $\varepsilon=10^{-9}$ is chosen to introduce extremely small cut elements. The results with stabilization and without stabilization are shown in
Fig.~\ref{fig:cantilever_epsilon_results}(\subref{fig:cantilever_epsilon_result},
\subref{fig:cantilever_epsilon_result_nostable}), respectively. We observe that the
expected structure (when compared to \cite{Dede_PhaseFieldITO})
is obtained with stabilization, whereas the structural layout is problematic
without stabilization. The reason is the insufficient imposition of the
boundary condition $\nabla\rho\cdot\mathbf{n}=0$, which guarantees monotonic
decrease in the total energy \cite{Wu_CahnHilliardReview}. This can be verified by
the energy evolution of the generalized Cahn--Hilliard model; see
Fig.~\ref{fig:cantilever_energy_histories}(\subref{fig:cantilever_energy_history},
\subref{fig:cantilever_energy_history_nostable}). Without stabilization, we
observe that around $t=10^{-5}$, the total energy $J^{\star}$ fluctuates up and
down, diverging the optimization path to a different local minimum.
In contrast, $J^{\star}$ monotonically decreases in the stabilized case.
Moreover, the tangent matrix of Newton's method (already diagonally scaled)
becomes nearly singular in the non-stabilized case, introducing relatively larger errors also in
the linear solver.

The difference is also reflected in the adaptive time stepping, for which we
follow \cite{Dede_PhaseFieldITO} to enhance efficiency; see
\ref{app:time_integration} for details. $\Delta t$ is adjusted
such that the error, estimated from the difference between the generalized-$\alpha$ and
backward Euler solutions, stays below a prescribed tolerance. With
stabilization, $\Delta t$ increases overall as the evolution proceeds, with
only mild reductions; see
Fig.~\ref{fig:cantilever_energy_histories}(\subref{fig:cantilever_dt_stable}).
Without stabilization, $\Delta t$ repeatedly drops by several orders of
magnitude; see
Fig.~\ref{fig:cantilever_energy_histories}(\subref{fig:cantilever_dt_nostable}).
This is because the estimated error is dominated by the DOFs associated with
the cut elements. Small cut elements lead to small stiffness terms, and thus
large variations in their corresponding control variables, which in turn
introduce large errors.

\begin{figure}[!t]
    \centering
    \includegraphics[width=0.6\textwidth]{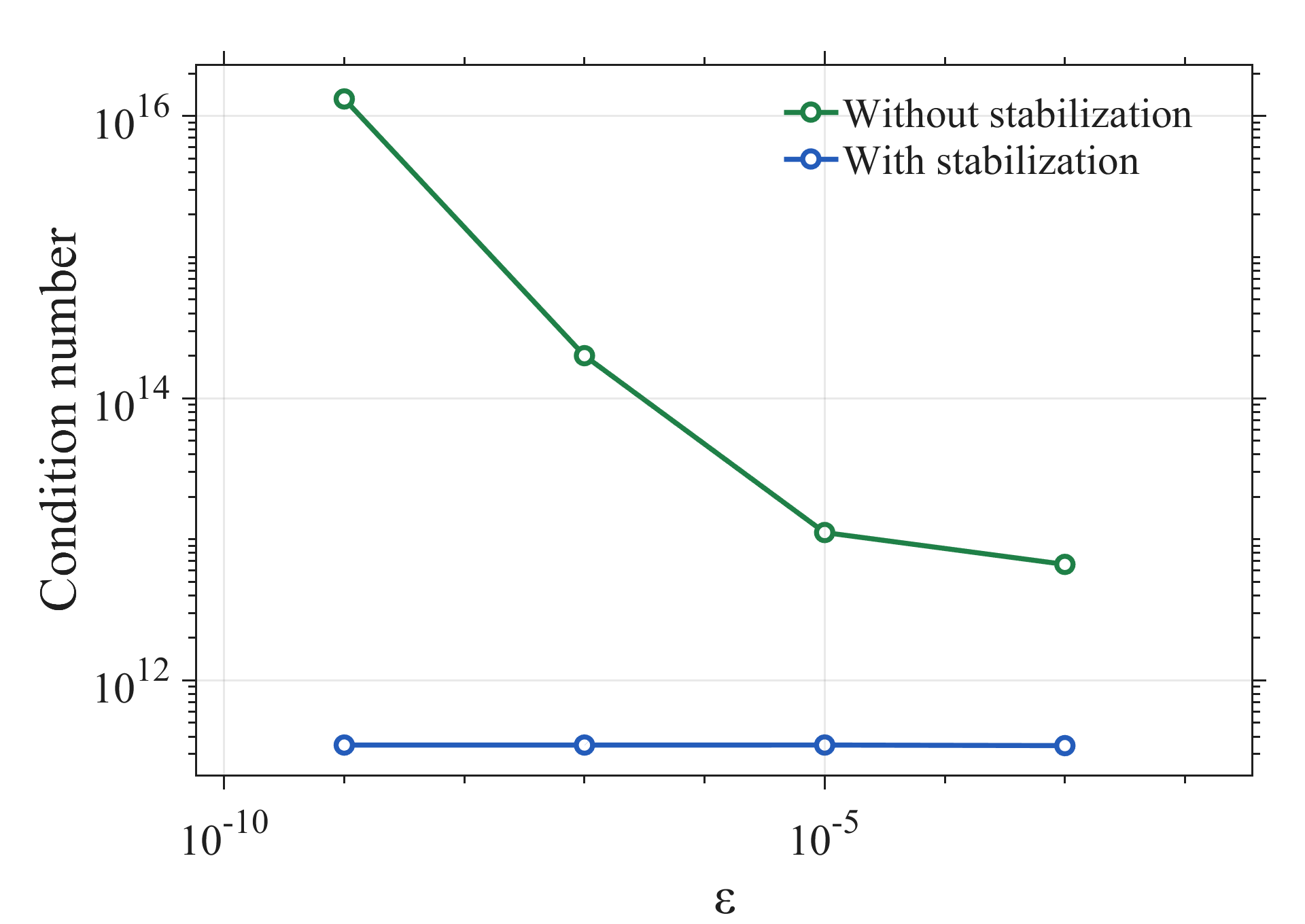}
    \caption{Conditioning of the tangent matrix of Newton's method at the first
    time step in the cantilever problem. Condition numbers are computed on the
    diagonally rescaled matrices.}
    \label{fig:cantilever_condition_number}
\end{figure}

Next, we study the conditioning of the tangent matrix (after diagonal
rescaling) of Newton's method by varying $\varepsilon$. In particular, we
choose the matrix at the first time step. The results are shown in
Fig.~\ref{fig:cantilever_condition_number}. Without stabilization, the
condition number grows rapidly as $\varepsilon$ decreases, whereas with minimal
stabilization, the condition number remains nearly constant, demonstrating the
effectiveness of the proposed method.

\begin{figure}[!t]
    \centering

    \begin{subfigure}[t]{0.495\textwidth}
        \centering
        \vspace{0pt}
        \includegraphics[width=\linewidth,height=0.59\linewidth,keepaspectratio]{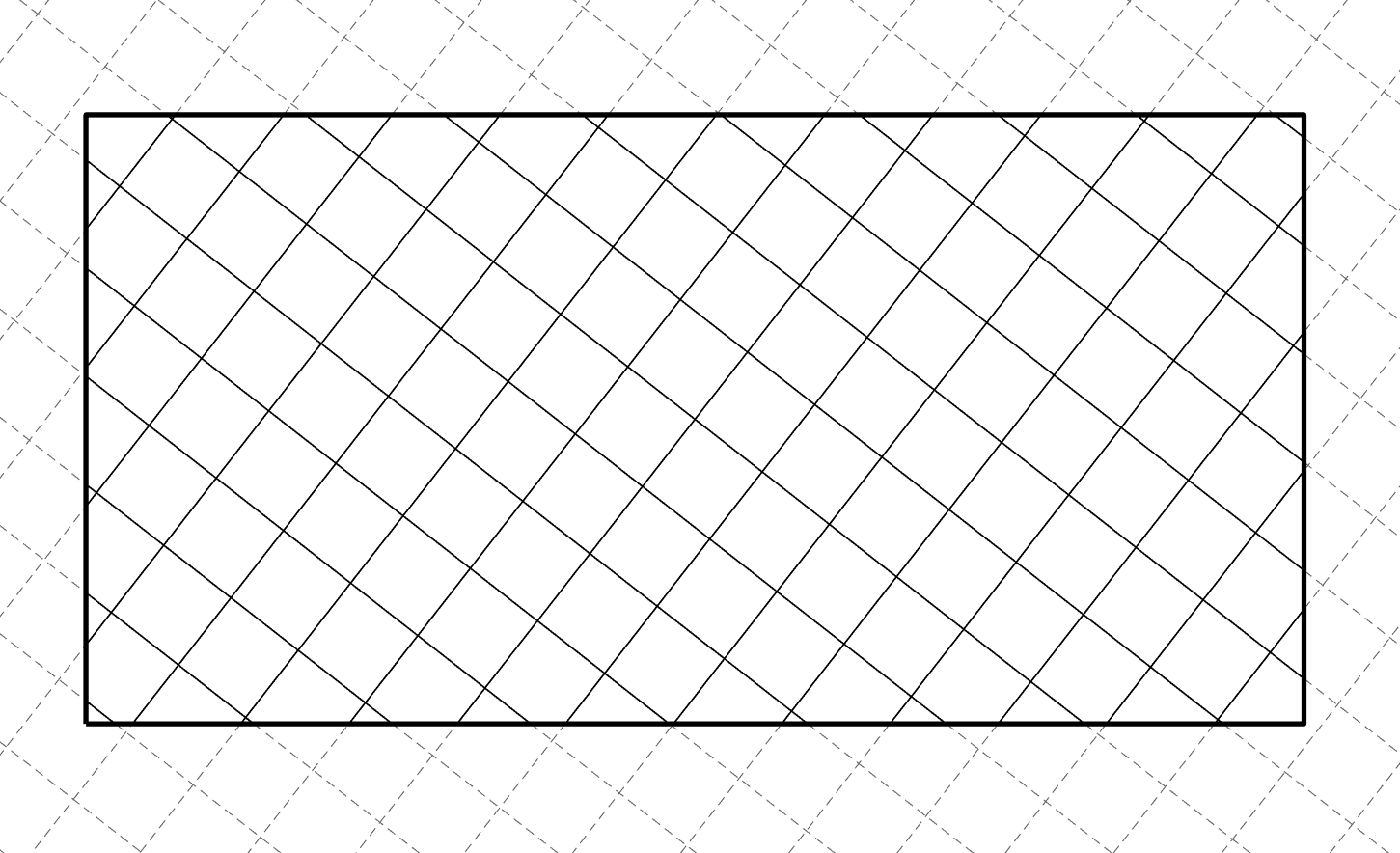}
        \caption{Background mesh rotated by $37.2^\circ$}
        \label{fig:cantilever_rotated_mesh_schematic}
    \end{subfigure}\hfill%
    \begin{subfigure}[t]{0.495\textwidth}
        \centering
        \vspace{0pt}
                \parbox[t][0.59\linewidth][b]{\linewidth}{\centering
        \includegraphics[width=0.909\linewidth]{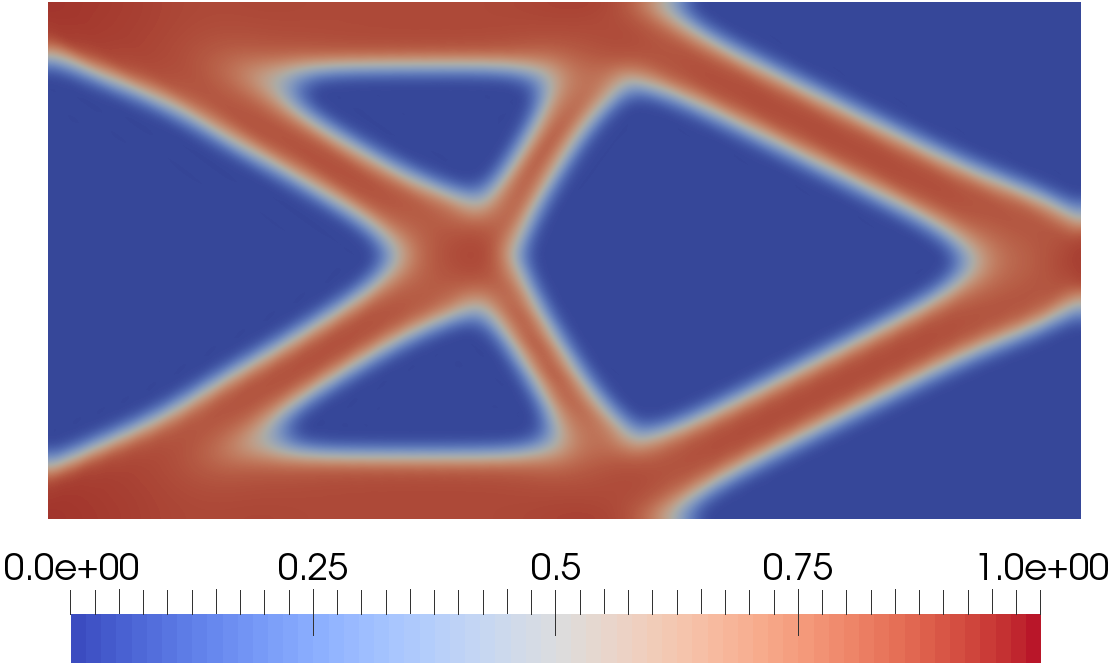}}
        \caption{Optimized topology at minimum $J_E^{\star}$}
        \label{fig:cantilever_rotated_result}
    \end{subfigure}

    \caption{Cantilever beam with a rotated background mesh.}
    \label{fig:cantilever_rotated_results}
\end{figure}

We then study a more general trimming case where the background mesh is
rotated by $37.2^\circ$, as shown in
Fig.~\ref{fig:cantilever_rotated_results}(\subref{fig:cantilever_rotated_mesh_schematic}).
The threshold for classifying bad elements is taken as $\theta=10^{-3}$; see
Eq.~\eqref{eq:area_ratio}. As $\Gamma_D$ is trimmed in this case, $\mathbf{u}=\mathbf{0}$ is imposed by Nitsche's method with $C_{\mathrm{pen}}^{u}=300E$. The optimized topology obtained with minimal
stabilization is shown in
Fig.~\ref{fig:cantilever_rotated_results}(\subref{fig:cantilever_rotated_result}),
where again, the expected structural layout is achieved.

\subsection{Non-design domain as a hole}
\label{subsec:one_hole_cantilever}

\begin{figure}[p]
    \centering

    \begin{subfigure}[t]{0.48\textwidth}
        \centering
        \vspace{0pt}
        \parbox[t][0.70\linewidth][c]{\linewidth}{\centering
        \includegraphics[width=\linewidth]{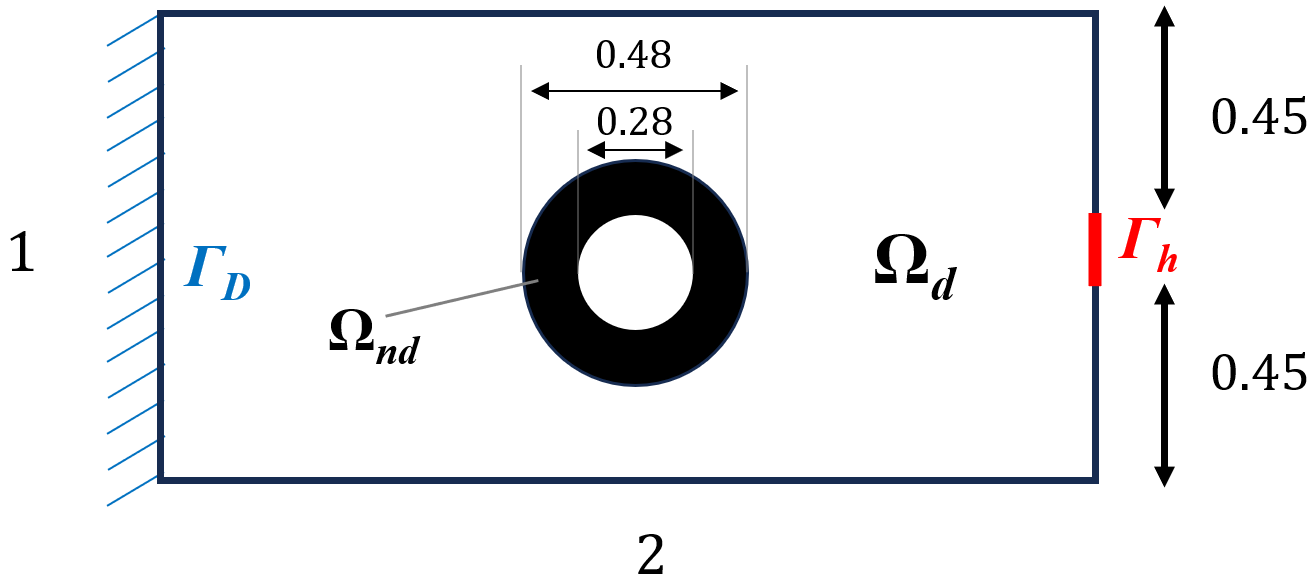}}
        \caption{Problem setup}
        \label{fig:one_hole_cantilever_geometry}
    \end{subfigure}\hfill%
    \begin{subfigure}[t]{0.48\textwidth}
        \centering
        \vspace{0pt}
        \parbox[t][0.70\linewidth][c]{\linewidth}{\centering
        \includegraphics[width=\linewidth,height=0.70\linewidth,keepaspectratio]{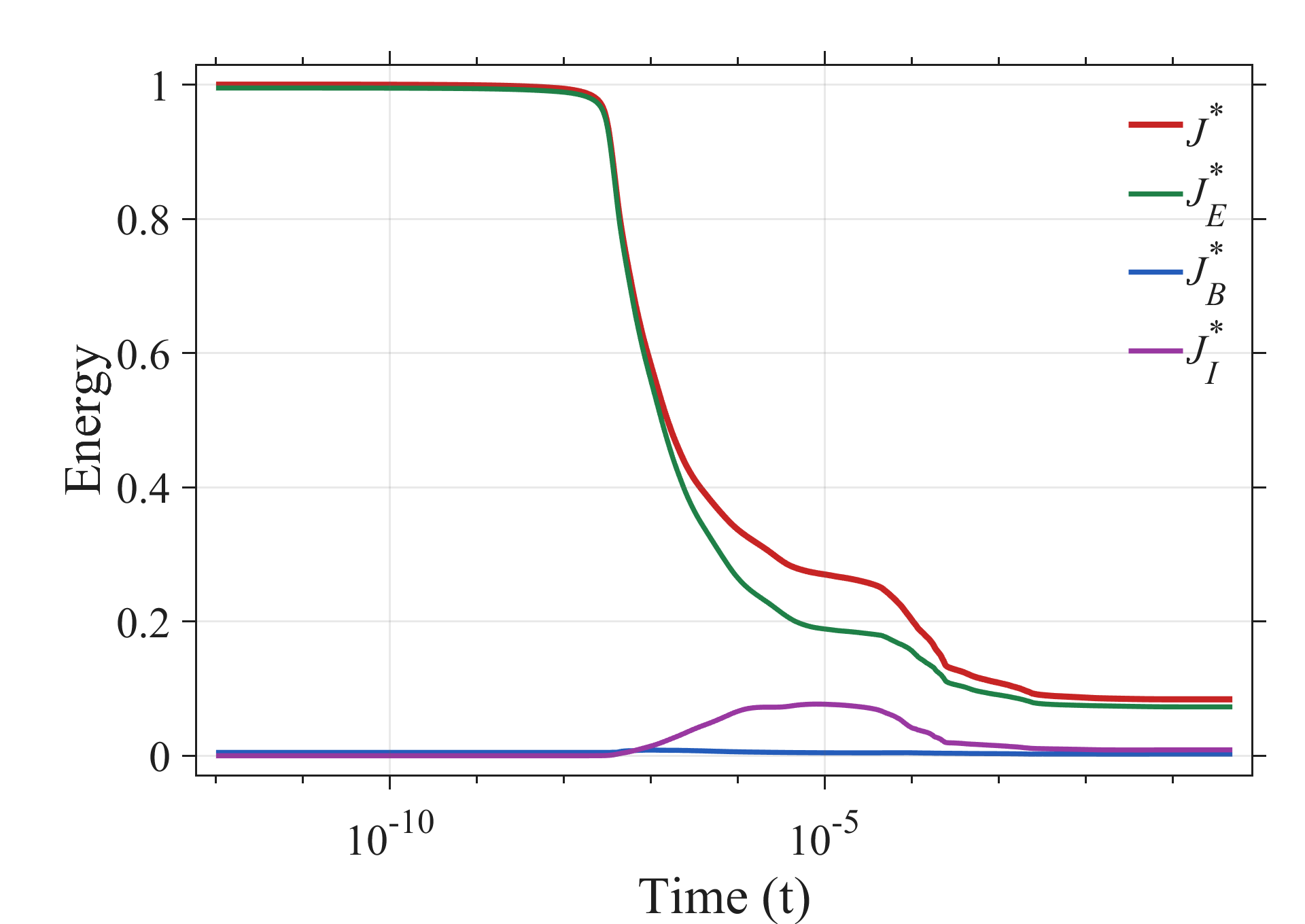}}
        \caption{Normalized energy evolution}
        \label{fig:one_hole_cantilever_energy}
    \end{subfigure}

    \par\medskip

    \begin{subfigure}[t]{0.48\textwidth}
        \centering
        \vspace{0pt}
        \parbox[t][0.662\linewidth][t]{\linewidth}{\centering
        \includegraphics[width=\linewidth]{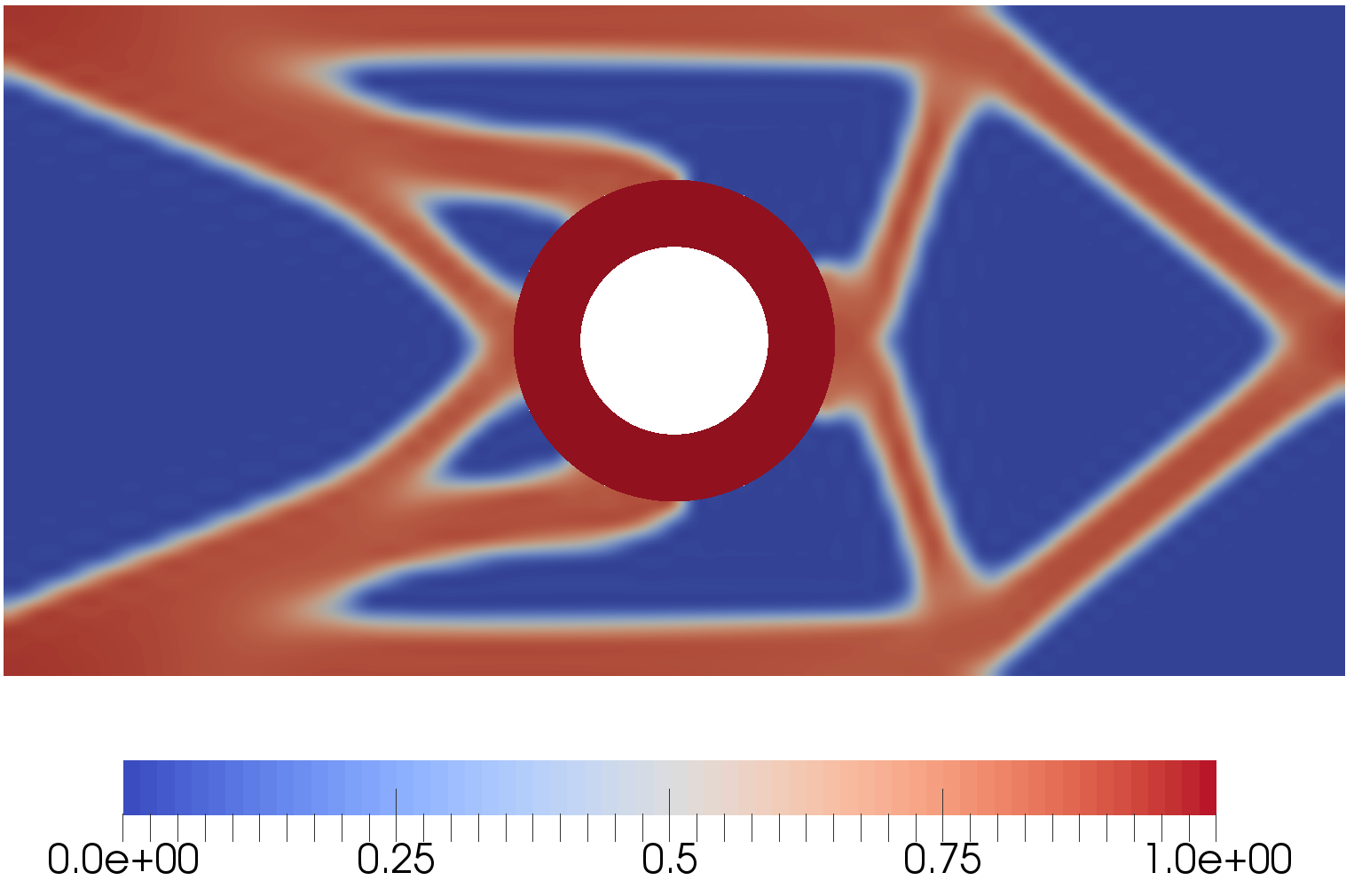}}
        \caption{Our method}
        \label{fig:one_hole_cantilever_phasefield}
    \end{subfigure}\hfill%
    \begin{subfigure}[t]{0.48\textwidth}
        \centering
        \vspace{0pt}
        \parbox[t][0.662\linewidth][t]{\linewidth}{\centering
        \includegraphics[width=\linewidth]{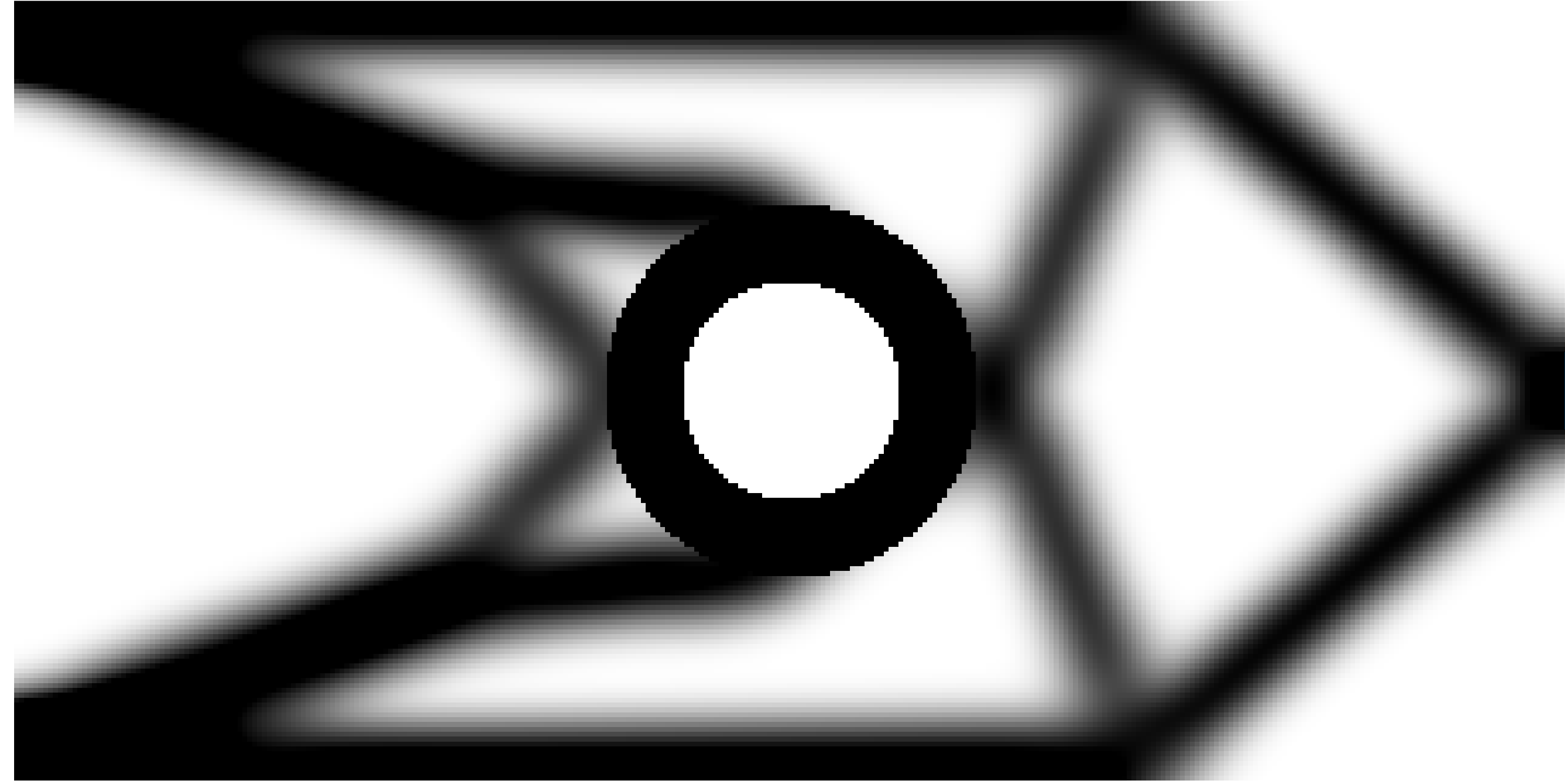}}
        \caption{SIMP}
        \label{fig:one_hole_cantilever_simp}
    \end{subfigure}

    \caption{TO with a non-design domain as a hole.}
    \label{fig:one_hole_cantilever_results}

    \vspace{\floatsep}

    \begin{subfigure}[t]{0.34\textwidth}
        \centering
        \includegraphics[width=\linewidth]{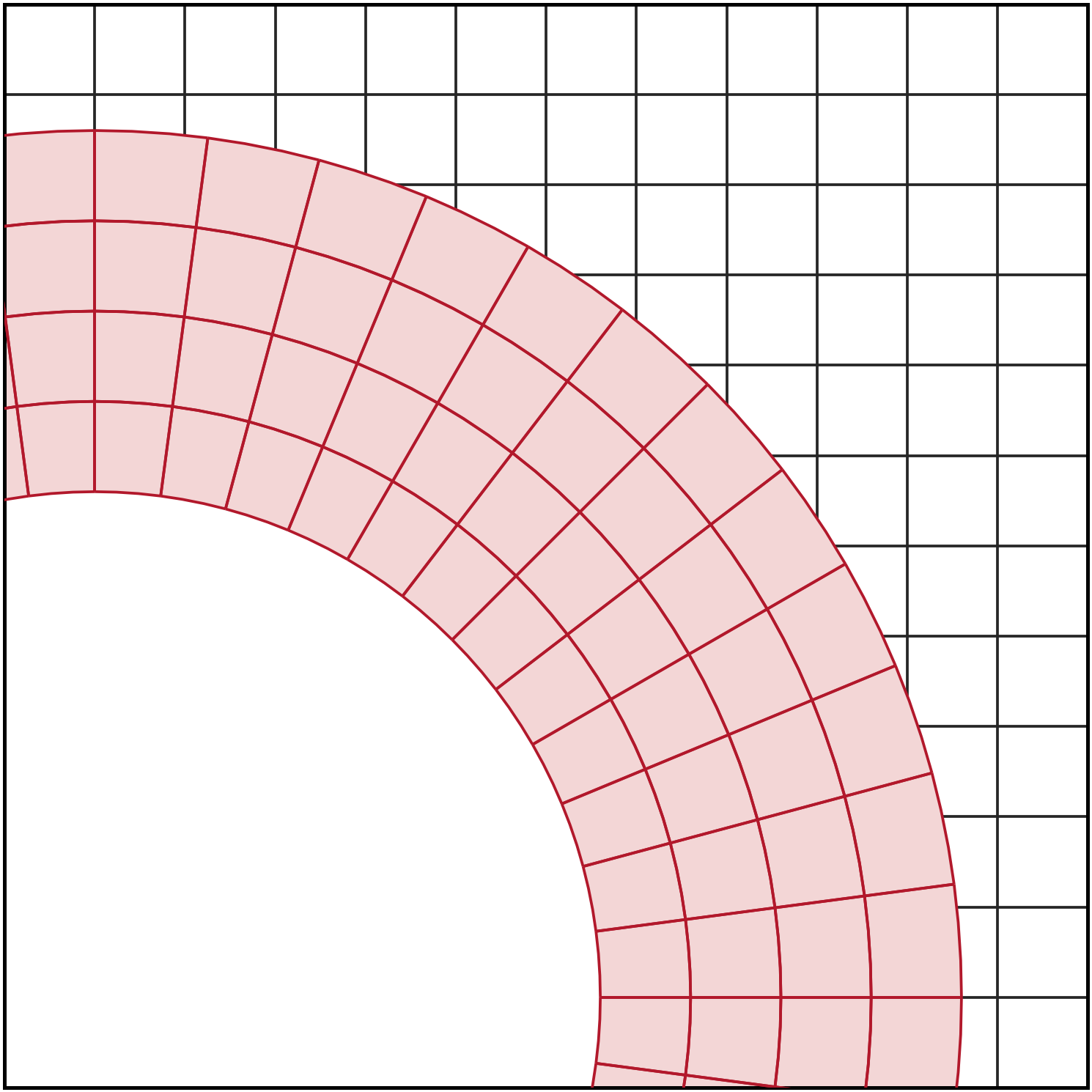}
        \caption{Our method}
        \label{fig:one_hole_cantilever_mesh_pf}
    \end{subfigure}\hspace{0.08\textwidth}%
    \begin{subfigure}[t]{0.34\textwidth}
        \centering
        \includegraphics[width=\linewidth]{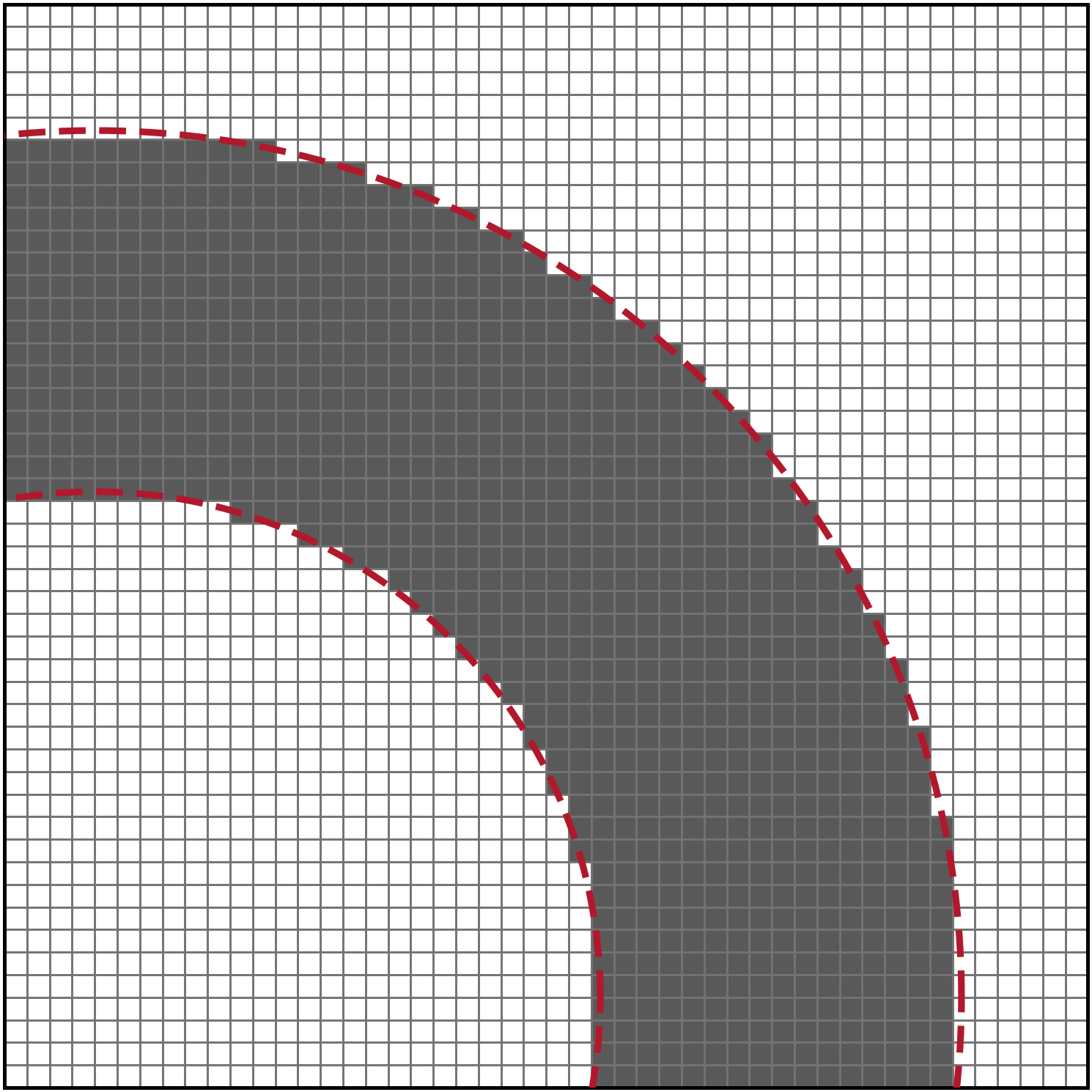}
        \caption{SIMP}
        \label{fig:one_hole_cantilever_mesh_simp}
    \end{subfigure}

    \caption{Meshes around the upper-right quarter of the
    non-design domain: (a) The background mesh (black) and the conformal layer
    (red); (b) Passive elements (gray) representing the annulus, where the dashed lines
    indicate the exact boundaries of the annulus.}
    \label{fig:one_hole_cantilever_meshes}
\end{figure}

We now study the treatment of non-design domains. As shown in
Fig.~\ref{fig:one_hole_cantilever_results}(\subref{fig:one_hole_cantilever_geometry}),
a hole in the cantilever beam is a feature to be preserved: the annular region
around it is a non-design domain $\Omega^{\text{nd}}$ that remains solid,
whereas the rest of the region, i.e., the design domain $\Omega^{\text{d}}$, can
be optimized. The design domain is discretized by an $80\times40$ background
mesh, whereas
$\Omega^{\text{nd}}$ is discretized by a conformal layer of $4\times48$
elements; see
Fig.~\ref{fig:one_hole_cantilever_meshes}(\subref{fig:one_hole_cantilever_mesh_pf}).
Clearly, the background mesh is cut by the outer boundary of
$\Omega^{\text{nd}}$. Bicubic B-splines and biquadratic NURBS are used for
the background mesh and the conformal layer, respectively.
The material parameters are $E=1$ Pa and $\nu=0.3$. The left edge $\Gamma_D$ is
fixed, and a uniform downward traction $F=-0.001$ N/m is applied to the middle
segment $\Gamma_h$ of the right edge. The prescribed material fraction is
$0.34$, and the phase-field parameters are $\bar{\lambda}=2$ and
$\bar{\gamma}=200$. The penalty parameter is $C_{\mathrm{pen}}^{\rho}=10$.

The results are shown in
Fig.~\ref{fig:one_hole_cantilever_results}(\subref{fig:one_hole_cantilever_energy}--\subref{fig:one_hole_cantilever_simp}).
Expected energy evolutions and structural layout are achieved. For
comparison, the reference result is obtained by an open-source code
\cite{Sigmund_99LineTO,Andreassen_88LineTO}, which is based on the Solid
Isotropic Material Penalization (SIMP) method. A $320\times160$ linear
mesh is used in the finite element analysis for SIMP.

Domains are represented differently in our method and SIMP. Our method
features exact representations of conic sections (e.g., circles) and
flexible treatment of design and non-design domains via trimming. In
contrast, the curved boundaries are approximated by staircases in SIMP, and
moreover, a conformal mesh is needed; see
Fig.~\ref{fig:one_hole_cantilever_meshes}(\subref{fig:one_hole_cantilever_mesh_simp}),
where the non-design domain is handled through the so-called passive elements
because their corresponding densities are fixed and do not serve as design
variables.

\subsection{Multiple non-design domains}
\label{subsec:multiple_nondesign}

We now demonstrate that the IBCM-based method can deal with multiple
features as non-design domains in a flexible way. Two examples are considered:
a rocker arm inspired by \cite{Yamada_FictitiousInterface} and a three-hole
cantilever beam inspired by \cite{Li_BSplineContact}. Their problem setups are
shown in Fig.~\ref{fig:multi_nd_setups}.

\newlength{\multiNDsetupheight}
\setlength{\multiNDsetupheight}{0.40\textwidth}
\begin{figure}[!t]
    \centering

    \begin{subfigure}[t]{0.55\textwidth}
        \centering
        \vspace{0pt}
        \includegraphics[height=\multiNDsetupheight,width=\linewidth,keepaspectratio]{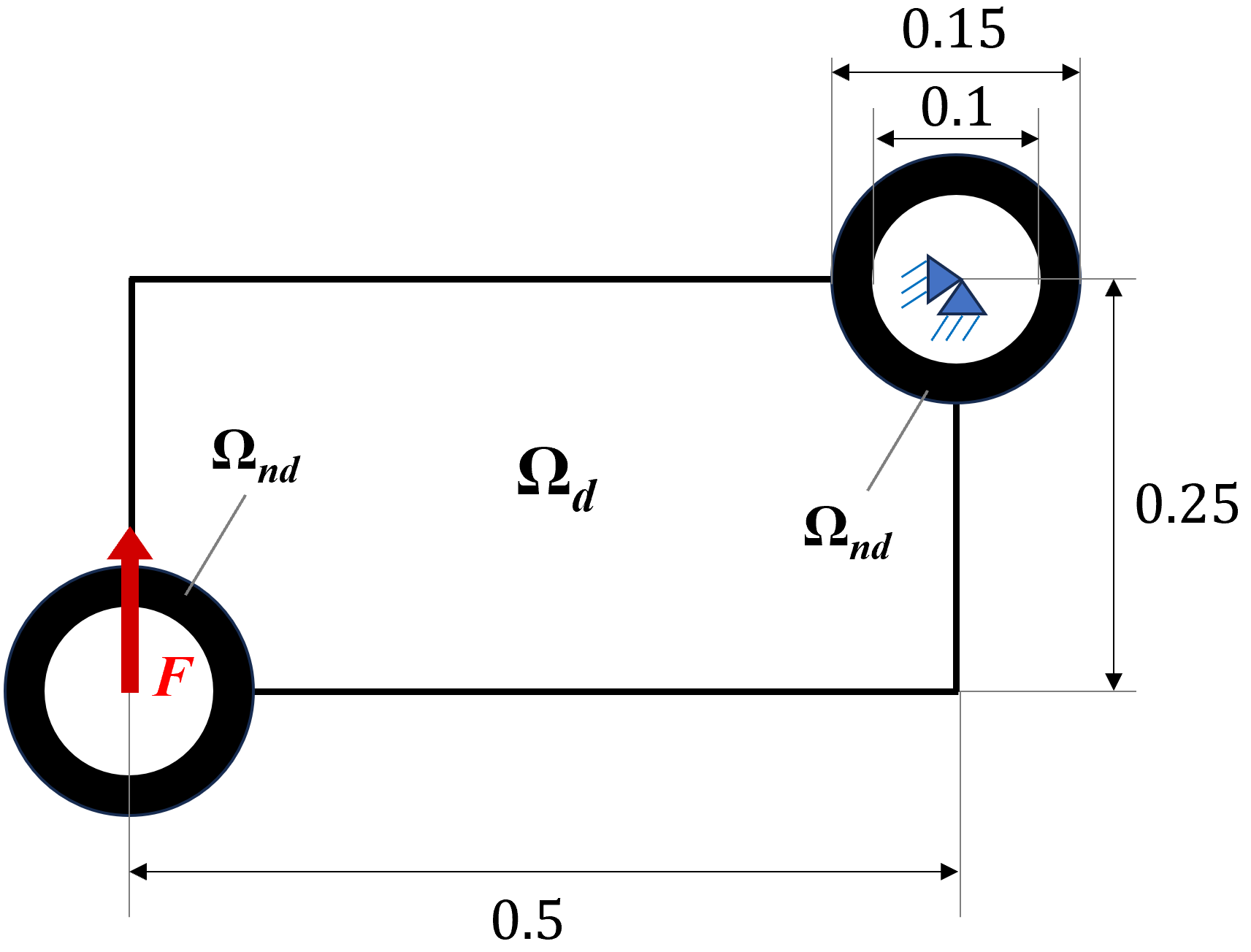}
        \caption{Rocker arm}
        \label{fig:rocker_arm_geometry}
    \end{subfigure}\hfill%
    \begin{subfigure}[t]{0.42\textwidth}
        \centering
        \vspace{0pt}
        \includegraphics[height=\multiNDsetupheight,width=\linewidth,keepaspectratio]{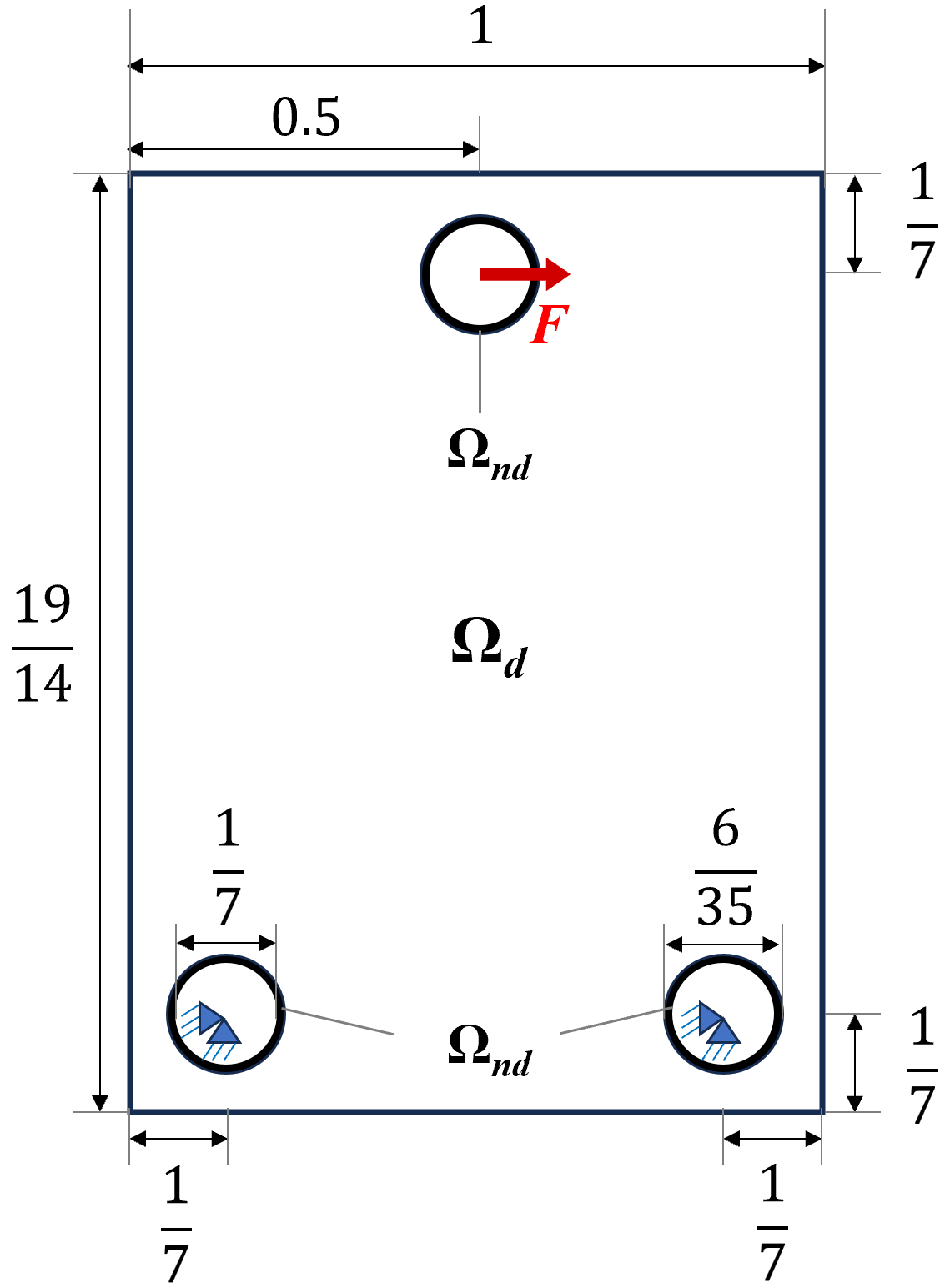}
        \caption{Three-hole cantilever beam}
        \label{fig:three_hole_geometry}
    \end{subfigure}

    \caption{Problem setups with multiple non-design domains.}
    \label{fig:multi_nd_setups}
\end{figure}

First, in the rocker arm problem, the load is applied to the lower-left ring,
whereas the inner boundary of the upper-right ring is fixed. The two rings
have the same dimension. The design domain of size $0.5\,\mathrm{m}\times0.25\,\mathrm{m}$ is discretized by a $186\times93$ background
mesh ($h=0.0038$ m), and each non-design ring is discretized by $8\times80$
elements. The material parameters are $E=2.0\times10^{12}$ Pa and $\nu=0.31$.
The prescribed material fraction is $0.45$, and the phase-field parameters are
$\bar{\lambda}=20$ and $\bar{\gamma}=500$. The penalty parameter is $C_{\mathrm{pen}}^{\rho}=100$.

\begin{figure}[!t]
    \centering

    \begin{subfigure}[t]{0.48\textwidth}
        \centering
        \vspace{0pt}
        \includegraphics[width=\linewidth]{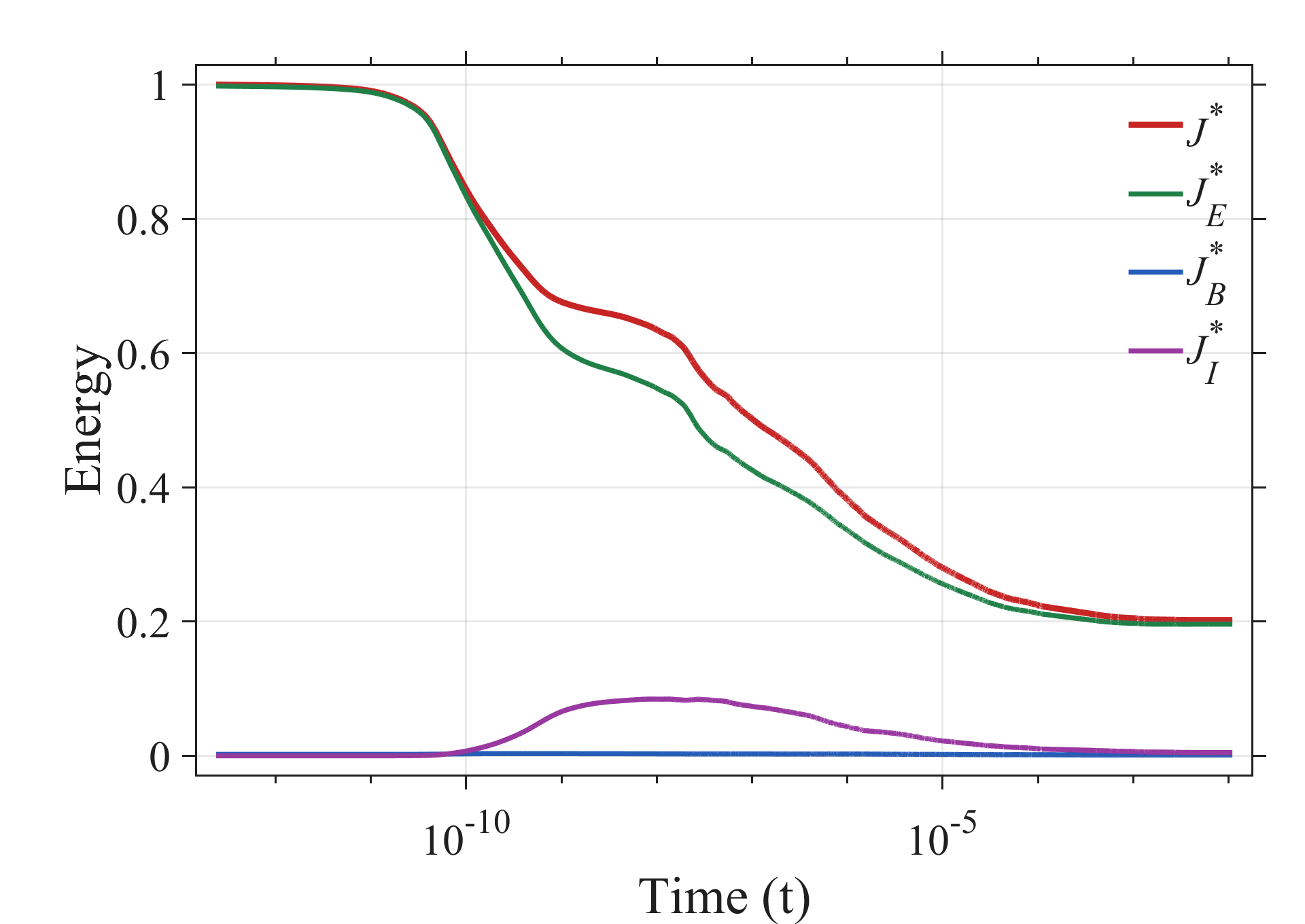}
        \caption{Rocker arm}
        \label{fig:rocker_arm_energy}
    \end{subfigure}\hfill%
    \begin{subfigure}[t]{0.48\textwidth}
        \centering
        \vspace{0pt}
        \includegraphics[width=\linewidth]{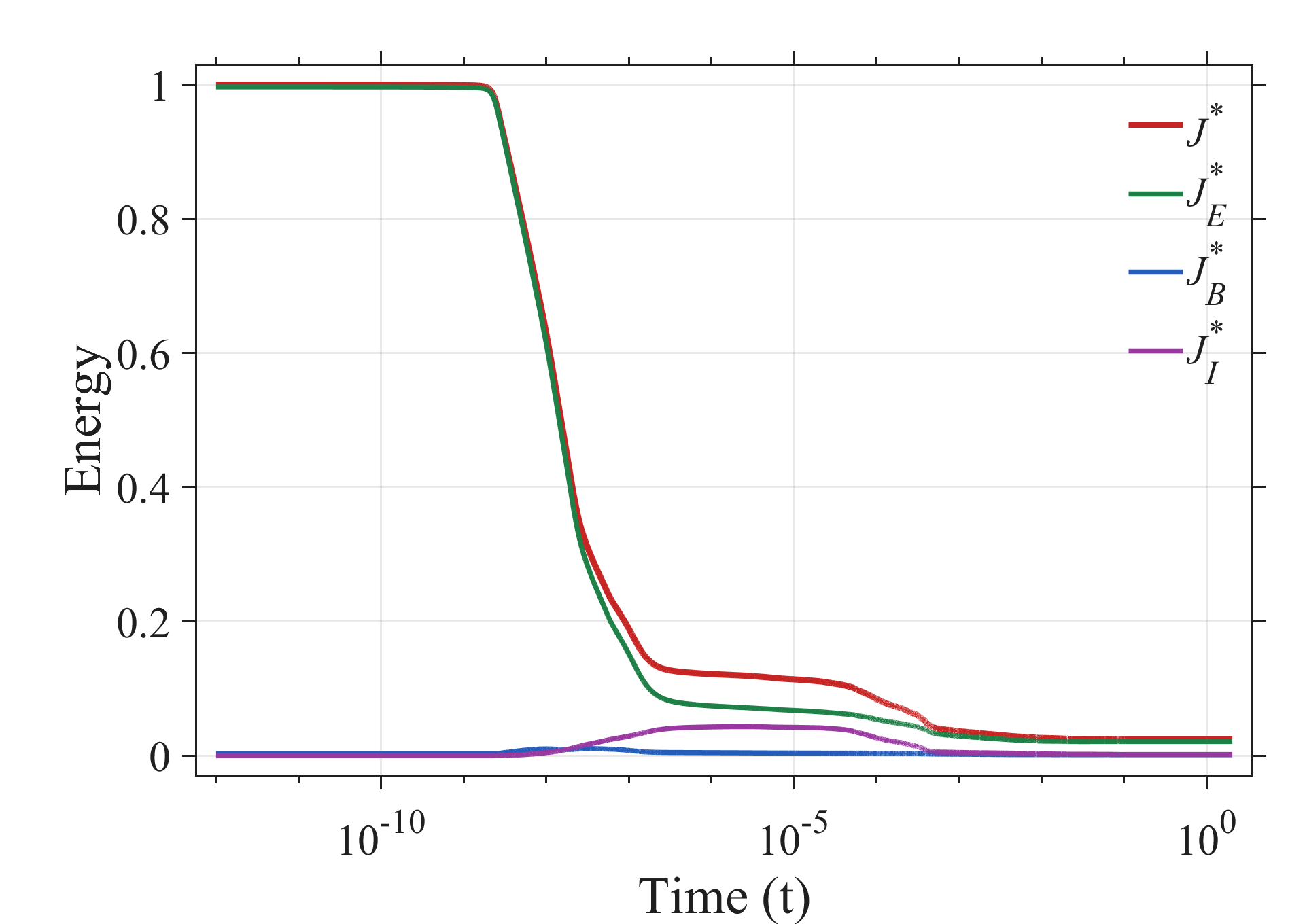}
        \caption{Three-hole cantilever beam}
        \label{fig:three_hole_energy}
    \end{subfigure}

    \caption{Normalized energy evolution with multiple non-design domains.}
    \label{fig:multi_nd_energies}
\end{figure}

\begin{figure}[!t]
    \centering

    \begin{subfigure}[t]{0.495\textwidth}
        \centering
        \vspace{0pt}
        \includegraphics[width=\linewidth]{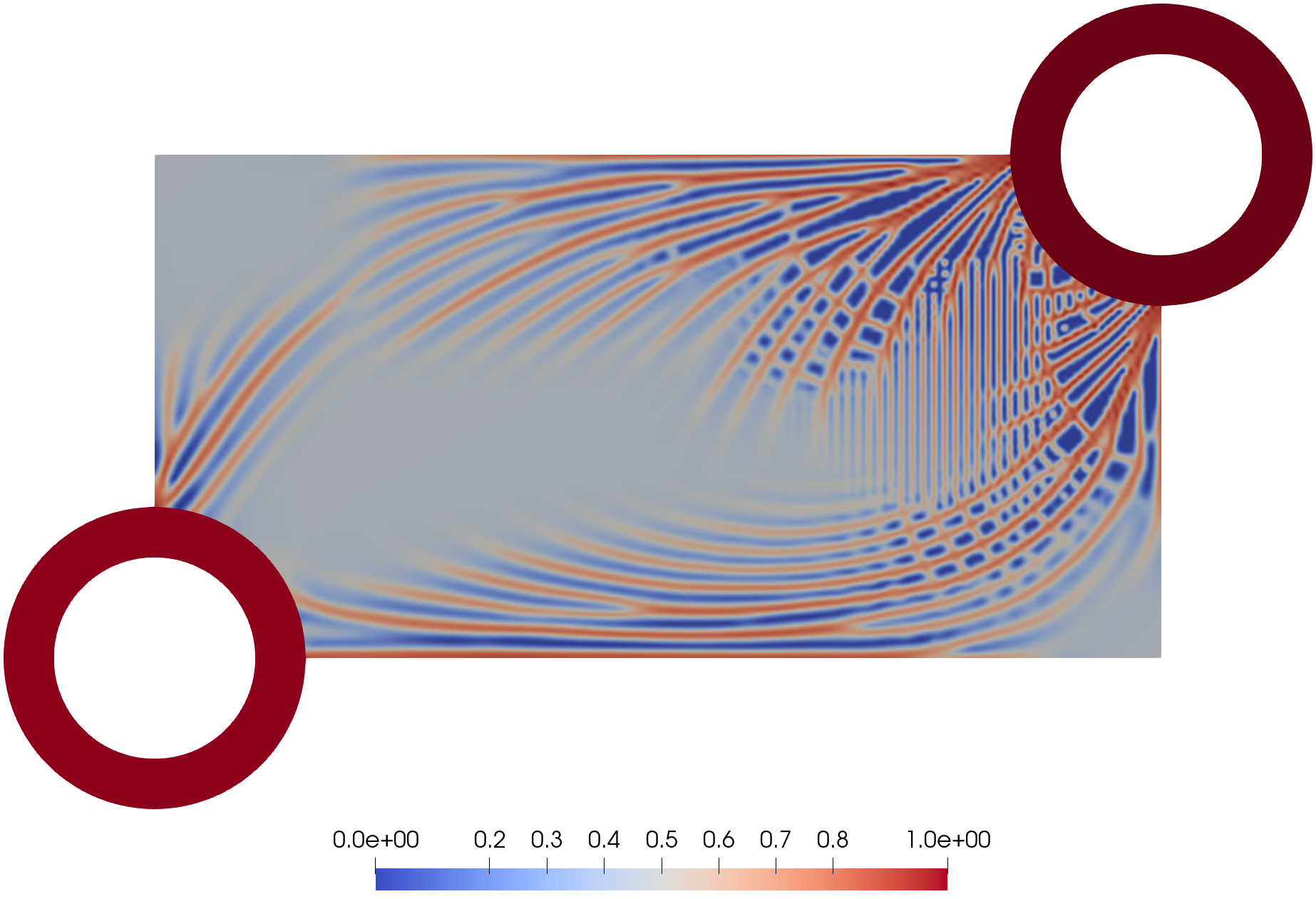}
        \caption{$t=3.12\times10^{-7}$}
        \label{fig:rocker_arm_t1}
    \end{subfigure}\hfill%
    \begin{subfigure}[t]{0.495\textwidth}
        \centering
        \vspace{0pt}
        \includegraphics[width=\linewidth]{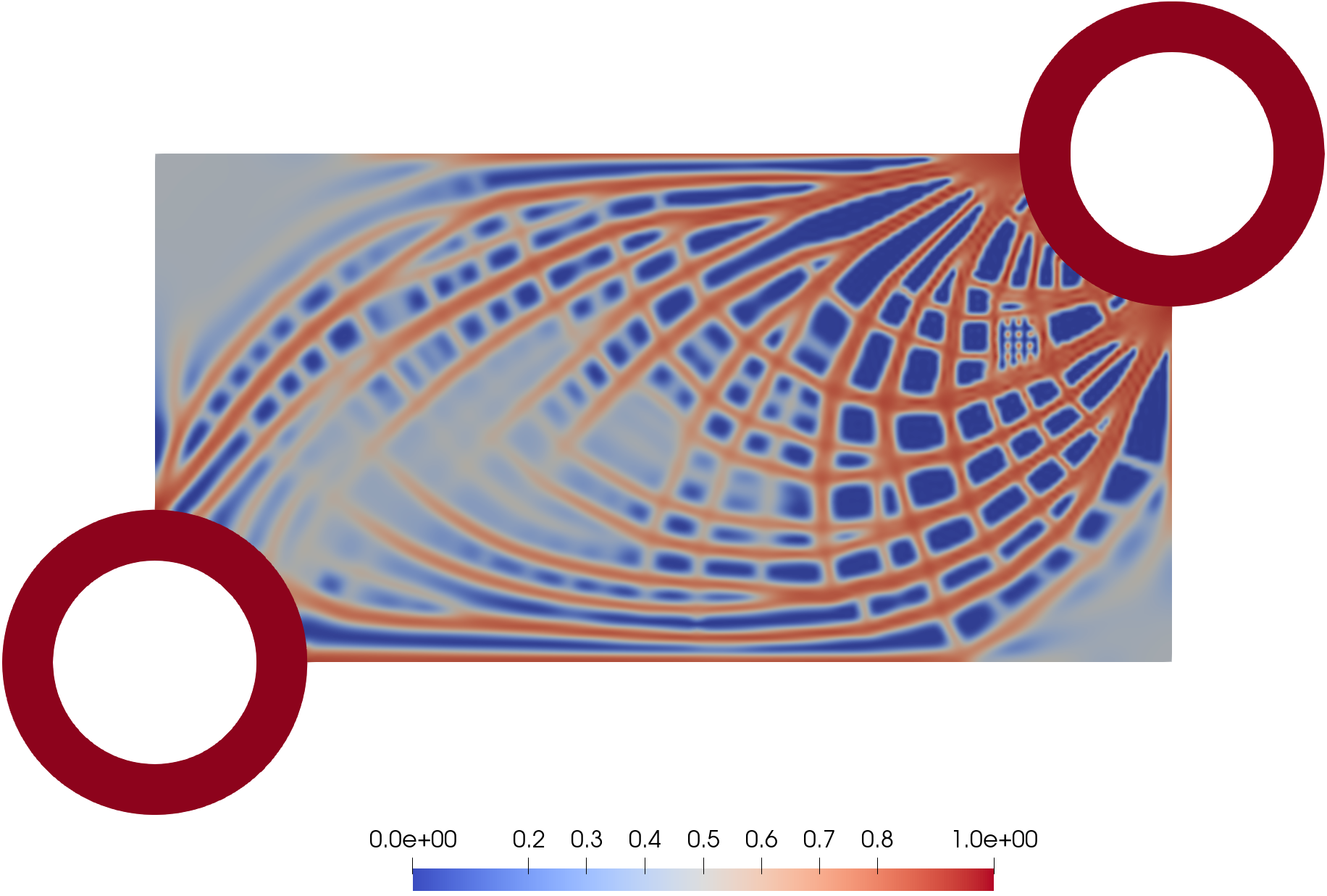}
        \caption{$t=1.38\times10^{-6}$}
        \label{fig:rocker_arm_t2}
    \end{subfigure}

    \par\vspace{6pt}

    \begin{subfigure}[t]{0.495\textwidth}
        \centering
        \vspace{0pt}
        \includegraphics[width=\linewidth]{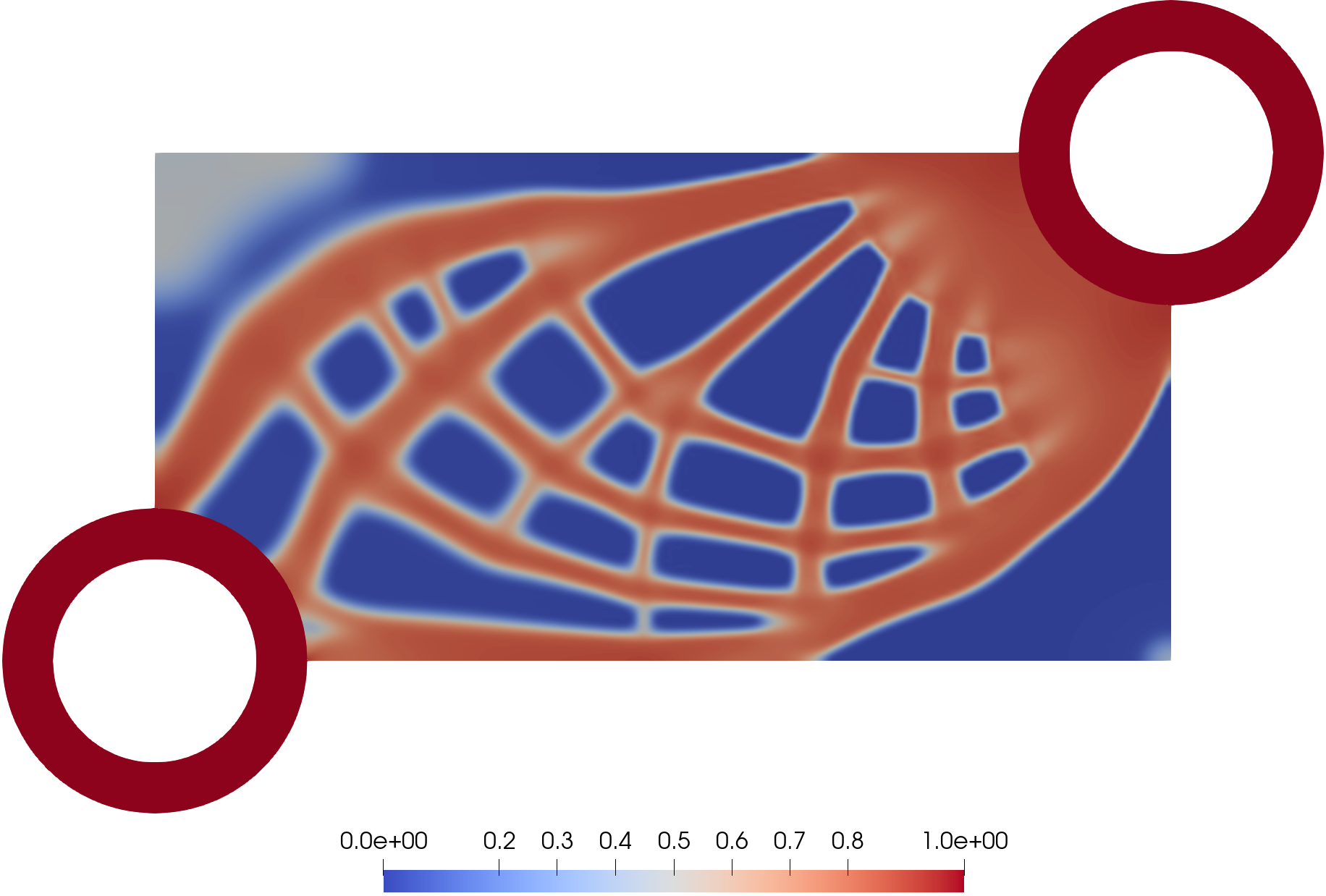}
        \caption{$t=4.84\times10^{-5}$}
        \label{fig:rocker_arm_t3}
    \end{subfigure}\hfill%
    \begin{subfigure}[t]{0.495\textwidth}
        \centering
        \vspace{0pt}
        \includegraphics[width=\linewidth]{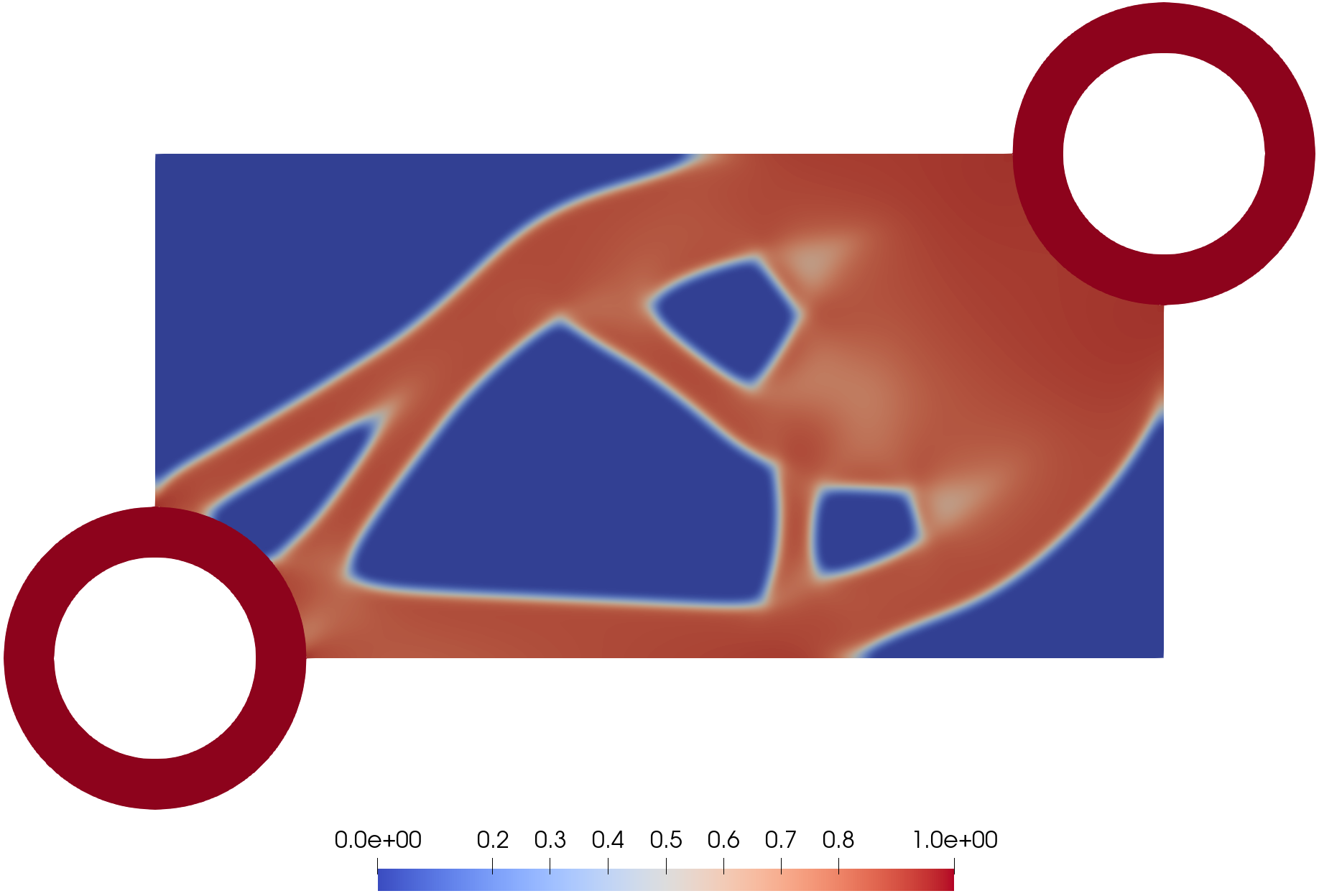}
        \caption{$t=3.32\times10^{-3}$ at minimum $J_E^{\star}$}
        \label{fig:rocker_arm_t4}
    \end{subfigure}

    \caption{Evolution of the topology for the rocker-arm problem.}
    \label{fig:rocker_arm_evolution}
\end{figure}

The normalized energies are shown in
Fig.~\ref{fig:multi_nd_energies}(\subref{fig:rocker_arm_energy}), where
$J^{\star}$ monotonically decreases. The topology evolution is shown in Fig.~\ref{fig:rocker_arm_evolution}. Compared with
\cite{Yamada_FictitiousInterface}, where the problem is solved in 3D on an
unstructured mesh conforming to the rings, we achieve a similar result and do not need a conformal mesh.

The second multi-feature example is a three-hole cantilever beam shown in
Fig.~\ref{fig:multi_nd_setups}(\subref{fig:three_hole_geometry}), where the
upper ring carries the load, and the inner boundaries of the two lower rings
are fixed. The design domain of size $1\,\mathrm{m}\times1.357\,\mathrm{m}$ is discretized by a $112\times152$
background mesh ($h=0.0126$ m), and each ring (i.e., the non-design domain) is discretized by
$2\times48$ elements. The Young's modulus is $E=207.4\times10^9$ Pa, the Poisson's ratio is $\nu=0.3$, and the prescribed material fraction is
$0.2213$. The phase-field parameters are $\bar{\lambda}=1$ and
$\bar{\gamma}=300$, and the penalty parameter is $C_{\mathrm{pen}}^{\rho}=100$.

\begin{figure}[htbp]
    \centering

    \begin{subfigure}[t]{0.32\textwidth}
        \centering
        \vspace{0pt}
        \includegraphics[width=\linewidth]{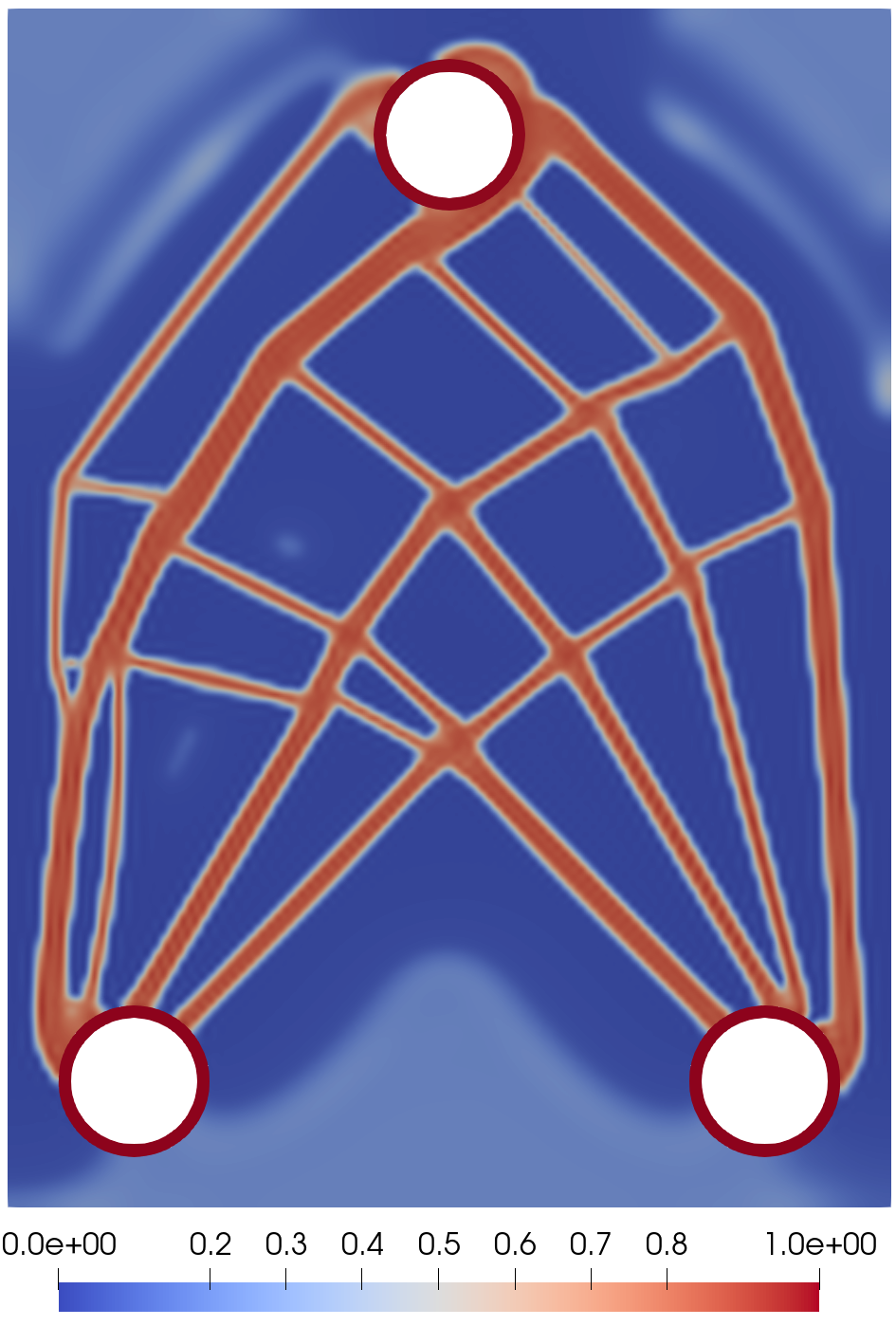}
        \caption{$t=5.43\times10^{-4}$}
        \label{fig:three_hole_t1}
    \end{subfigure}\hfill%
    \begin{subfigure}[t]{0.32\textwidth}
        \centering
        \vspace{0pt}
        \includegraphics[width=\linewidth]{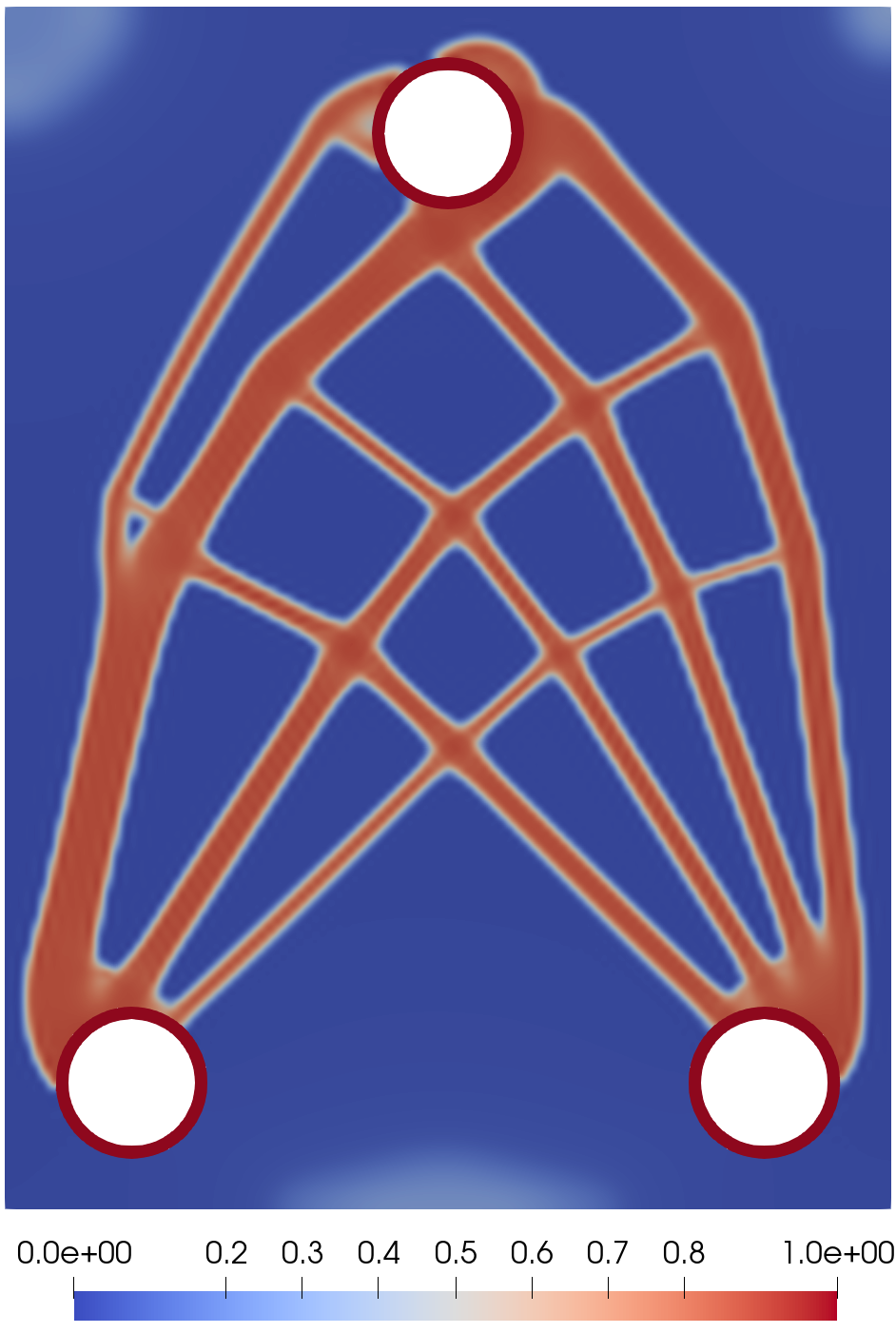}
        \caption{$t=3.84\times10^{-3}$}
        \label{fig:three_hole_t2}
    \end{subfigure}\hfill%
    \begin{subfigure}[t]{0.32\textwidth}
        \centering
        \vspace{0pt}
        \includegraphics[width=\linewidth]{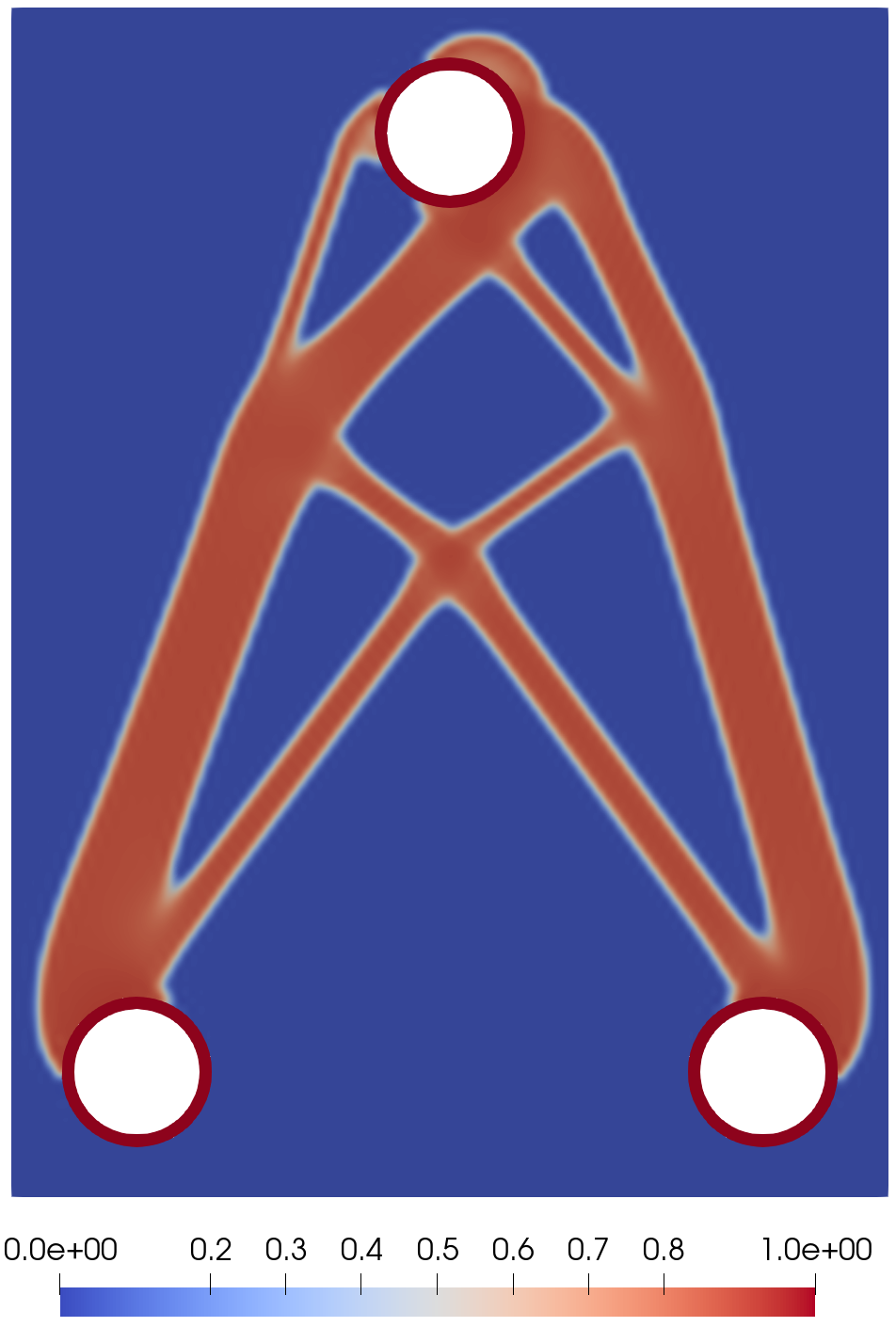}
        \caption{$t=3.29\times10^{-2}$ at minimum $J_E^{\star}$}
        \label{fig:three_hole_t3}
    \end{subfigure}

    \caption{Evolution of the topology for the three-hole cantilever beam.}
    \label{fig:three_hole_results}
\end{figure}

The normalized energies are shown in
Fig.~\ref{fig:multi_nd_energies}(\subref{fig:three_hole_energy}), where
$J^{\star}$ also monotonically decreases. The topology evolution is shown in
Fig.~\ref{fig:three_hole_results}, where redundant branches are gradually
removed, and the three preserved rings are finally connected by a clear load-transfer path. These examples show that the proposed formulation can
treat multiple non-design domains in a unified and flexible way.

\section{Conclusion}

In this work, we present a stabilized isogeometric topology optimization
framework for complex design domains represented by trimmed geometries.
Minimum compliance problem is studied through the generalized Cahn--Hilliard
model coupled with the equilibrium equation. Nitsche's method is adopted to impose boundary and interface conditions, where minimal stabilization is
extended to the nonlinear, time-dependent, and the fourth-order problem.
Only the flux terms from bad elements are stabilized, by replacing them with
their counterparts from good neighbors via polynomial extension. Moreover,
non-design domains are handled by the Immersed Boundary-Conformal Method
through a one-sided flux. This flux is taken from the untrimmed conformal
layer of the non-design domain, so it needs no stabilization. A unified
treatment is developed to impose transmission conditions across interfaces
between design and non-design domains, accounting for material-dependent
conditions on the side of the design domain.

A variety of numerical examples demonstrate the accuracy and robustness of
the method for severely cut meshes and the flexibility in retaining multiple
structural features during topology optimization. In particular, minimal
stabilization proves effective in the three-dimensional biharmonic problem, where optimal convergence rates and trimming-independent conditioning are
achieved. In the cantilever beam with
extremely small cut elements, it effectively imposes the boundary condition
$\nabla\rho\cdot\mathbf{n}=0$, which guarantees that expected structures can
be obtained. In the examples with non-design domains, the features are
represented exactly by conformal layers, and can be treated independently
and flexibly through the cut background mesh.

In the future, we will focus on robust numerical integration for
three-dimensional trimmed geometries and extend the proposed method to
handle complex design and non-design domains represented by B-reps.

\section*{Acknowledgements}

Haoyu Su and Xiaodong Wei are partially supported by National Natural
Science Foundation of China (No. 12494550/12494555 and No. 12571408).

\appendix
\section{Dimensionless formulation and parameter selection}
\label{app:dimensionless}

A dimensionless formulation facilitates the selection of the phase-field parameters
to balance the different energy contributions in the coupled problem.
We adopt the scaling conventions of
\cite{Dede_PhaseFieldITO,Sheng_ThinWalledComplex} and extend them to
the coupling of design and non-design domains.

\subsection{Dimensionless form of the coupled system}
\label{app:dimensionless_system}

With reference length, time, mobility, body-force and traction
scales $L_0$, $T_0$, $M_0$, $f_0$ and $h_0$, respectively, we define
\begin{equation}
\begin{aligned}
\mathbf x^{\star}&=\frac{\mathbf x}{L_0},
& t^{\star}&=\frac{t}{T_0},
& \mathbf u^{\star}&=\frac{\mathbf u}{L_0},
& \rho^{\star}&=\rho,
\\
\mathbf f^{\star}&=\frac{\mathbf f}{f_0},
& \mathbf h^{\star}&=\frac{\mathbf h}{h_0},
& M^{\star}&=\frac{M}{M_0},
& h^{\star}&=\frac{h}{L_0}.
\end{aligned}
\label{eq:app_dimensionless_variables}
\end{equation}
A superscript $\star$ denotes the dimensionless quantity. Accordingly, $\Omega^{\mathrm d,\star}$, $\Omega^{\mathrm{nd},\star}$, $\Gamma_D^{\star}$, $\Gamma_N^{\star}$ and $\Gamma_I^{\star}$ are the corresponding domains and boundaries in the coordinates $\mathbf x^{\star}$, with the measures $\mathrm{d}\Omega^{\star}$ and $\mathrm{d}\Gamma^{\star}$, and $\nabla^{\star}=L_0\nabla$ and $\Delta^{\star}=L_0^2\Delta$ are the gradient and Laplacian with respect to $\mathbf x^{\star}$.
We have
\begin{equation}
\begin{aligned}
\boldsymbol\varepsilon^{\star}(\mathbf u^{\star})
&=\boldsymbol\varepsilon(\mathbf u),
& \boldsymbol\sigma^{\star}&=\frac{\boldsymbol\sigma}{E_0},
\\
\psi_E^{\star}&=\frac{\psi_E}{E_0},
& \psi_B^{\star}&=\psi_B,
&
\psi_I^{\star}&=L_0^2\psi_I.
\end{aligned}
\label{eq:app_energy_density_scaling}
\end{equation}
Eq.~\eqref{eq:app_energy_density_scaling} applies to both
$\Omega^{\mathrm d}$ and $\Omega^{\mathrm{nd}}$.
The dimensionless formulation depends on
$\mathbf D=(D_1,\ldots,D_8)$, with
\begin{equation}
\begin{aligned}
D_1&=\frac{L_0^4}{T_0\lambda M_0},
& D_2&=\frac{L_0^2}{\lambda},
&
D_3&=\gamma D_2,
& D_4&=\frac{E_0}{h_0},
\\
D_5&=\frac{f_0L_0}{h_0},
& D_6&=\frac{C_{\mathrm{pen}}^{\rho}}{\lambda M_0},
&
D_7&=\frac{C_{\mathrm{pen}}^u}{h_0},
& D_8&=\frac{\beta}{h_0}.
\end{aligned}
\label{eq:app_dimensionless_parameters}
\end{equation}
Here, $D_1$ scales the time-derivative term; $D_2$ and $D_3$ weigh
the bulk and strain energies relative to the interface energy,
respectively; $D_4$ and $D_5$ scale the elastic and body-force terms
relative to the traction term, respectively; and $D_6$, $D_7$ and $D_8$ scale the penalty terms of the phase-field boundary, the displacement
boundary and the mechanical interface, respectively.
We take $T_0=L_0^4/(\lambda M_0)$, yielding $D_1=1$.
In the numerical examples, we use $L_0=M_0=1$ and constant mobility $M=M_0$.
The scaled potential is given as
\begin{equation}
\begin{aligned}
z^{\star}(\rho^{\star},\mathbf u^{\mathrm d,\star};\mathbf D)
={}&-D_3\frac{\partial\psi_E^{\star}}{\partial\rho^{\star}}
(\rho^{\star},\mathbf u^{\mathrm d,\star})
+D_2\frac{\mathrm{d}\psi_B^{\star}}{\mathrm{d}\rho^{\star}}
(\rho^{\star})
-\Delta^{\star}\rho^{\star}.
\end{aligned}
\label{eq:app_dimensionless_potential}
\end{equation}
The scaled energies are defined as
\begin{equation}
\begin{aligned}
J^{\star}&=J_E^{\star}+J_B^{\star}+J_I^{\star},
\\
J_E^{\star}
&=D_3\left[
\int_{\Omega^{\mathrm d,\star}}
\psi_E^{\star}(\rho^{\star},\mathbf u^{\mathrm d,\star})
\,\mathrm{d}\Omega^{\star}
+\int_{\Omega^{\mathrm{nd},\star}}
\psi_E^{\star}(\mathbf u^{\mathrm{nd},\star})
\,\mathrm{d}\Omega^{\star}
\right],
\\
J_B^{\star}
&=D_2\int_{\Omega^{\mathrm d,\star}}
\psi_B^{\star}(\rho^{\star})\,\mathrm{d}\Omega^{\star},
\\
J_I^{\star}
&=\int_{\Omega^{\mathrm d,\star}}
\psi_I^{\star}(\rho^{\star})\,\mathrm{d}\Omega^{\star}.
\end{aligned}
\label{eq:numerical_scaled_energies}
\end{equation}
Note that $J_E^{\star}$ includes both design and non-design domains,
whereas $J_B^{\star}$ and $J_I^{\star}$ only involve the design domain.

With $\phi^{\star}=\phi$, $\mathbf v^{\star}=\mathbf v/L_0$ and
$\mathbf g^{\star}=\mathbf g/L_0$, the dimensionless residuals are defined as
\begingroup
\allowdisplaybreaks[1]
\begin{subequations}
\label{eq:app_residual_scaling}
\begin{align}
&R_\rho^{\star}(\rho^{\star},\mathbf u^{\star};\mathbf D)
(\phi^{\star})
\nonumber :=D_1\int_{\Omega^{\mathrm d,\star}}
\frac{\partial\rho^{\star}}{\partial t^{\star}}
\phi^{\star}\,\mathrm{d}\Omega^{\star}
\nonumber-\int_{\Omega^{\mathrm d,\star}}
z^{\star}\nabla^{\star}\!\cdot
\bigl(M^{\star}(\rho^{\star})\nabla^{\star}\phi^{\star}\bigr)
\,\mathrm{d}\Omega^{\star}
\nonumber\\*
&\qquad+\int_{\partial\Omega^{\mathrm d,\star}}
z^{\star}M^{\star}(\rho^{\star})
(\nabla^{\star}\phi^{\star}\cdot\mathbf n^{\mathrm d})
\,\mathrm{d}\Gamma^{\star}
\nonumber\\*
&\qquad+D_6\int_{\partial\Omega^{\mathrm d,\star}}
\frac{1}{h^{\star}}
(\nabla^{\star}\rho^{\star}\cdot\mathbf n^{\mathrm d})
(\nabla^{\star}\phi^{\star}\cdot\mathbf n^{\mathrm d})
\,\mathrm{d}\Gamma^{\star},
\label{eq:app_phase_residual}
\\[6pt]
&R_{\mathbf u}^{\star}(\rho^{\star},\mathbf u^{\star};\mathbf D)
(\mathbf v^{\star})
\nonumber:=D_4\int_{\Omega^{\mathrm d,\star}}
\boldsymbol\sigma^{\star}(\rho^{\star},\mathbf u^{\mathrm d,\star})
:\boldsymbol\varepsilon^{\star}(\mathbf v^{\mathrm d,\star})
\,\mathrm{d}\Omega^{\star}
\nonumber\\*
&\qquad+D_4\int_{\Omega^{\mathrm{nd},\star}}
\boldsymbol\sigma^{\star}(\mathbf u^{\mathrm{nd},\star})
:\boldsymbol\varepsilon^{\star}(\mathbf v^{\mathrm{nd},\star})
\,\mathrm{d}\Omega^{\star}
-D_5\int_{\Omega^{\mathrm d,\star}}
\mathbf f^{\star}\cdot\mathbf v^{\mathrm d,\star}
\,\mathrm{d}\Omega^{\star}
\nonumber\\
&\qquad-D_5\int_{\Omega^{\mathrm{nd},\star}}
\mathbf f^{\star}\cdot\mathbf v^{\mathrm{nd},\star}
\,\mathrm{d}\Omega^{\star}
\nonumber-\int_{\Gamma_N^{\star}}
\mathbf h^{\star}\cdot\mathbf v^{\star}\,\mathrm{d}\Gamma^{\star}
\nonumber\\
&\qquad-D_4\int_{\Gamma_D^{\star}}
\boldsymbol\sigma^{\star}(\rho^{\star},\mathbf u^{\star})\mathbf n
\cdot\mathbf v^{\star}\,\mathrm{d}\Gamma^{\star}
\nonumber-D_4\int_{\Gamma_D^{\star}}
\boldsymbol\sigma^{\star}(\rho^{\star},\mathbf v^{\star})\mathbf n
\cdot(\mathbf u^{\star}-\mathbf g^{\star})
\,\mathrm{d}\Gamma^{\star}
\nonumber\\
&\qquad-D_4\int_{\Gamma_I^{\star}}
g(\rho^{\star})
\boldsymbol\sigma^{\star}(\mathbf u^{\mathrm{nd},\star})
\mathbf n^{\mathrm{nd}}\cdot[\mathbf v^{\star}]
\,\mathrm{d}\Gamma^{\star}
\nonumber\\
&\qquad-D_4\int_{\Gamma_I^{\star}}
g(\rho^{\star})
\boldsymbol\sigma^{\star}(\mathbf v^{\mathrm{nd},\star})
\mathbf n^{\mathrm{nd}}\cdot[\mathbf u^{\star}]
\,\mathrm{d}\Gamma^{\star}
\nonumber\\*
&\qquad+D_7\int_{\Gamma_D^{\star}}
\frac{1}{h^{\star}}
(\mathbf u^{\star}-\mathbf g^{\star})\cdot\mathbf v^{\star}
\,\mathrm{d}\Gamma^{\star}
\nonumber\\
&\qquad+D_8\int_{\Gamma_I^{\star}}
g(\rho^{\star})
\left(\frac{1}{h^{\mathrm d,\star}}
+\frac{1}{h^{\mathrm{nd},\star}}\right)
[\mathbf u^{\star}]\cdot[\mathbf v^{\star}]
\,\mathrm{d}\Gamma^{\star}.
\label{eq:app_displacement_residual}
\end{align}
\end{subequations}
\endgroup

With the above definitions, the dimensionless coupled problem takes
the following form:
\begin{equation}
\begin{aligned}
&\text{find }\rho^{\star}\in\mathcal H,
\quad\mathbf u^{\star}\in\mathcal V,
\\[2pt]
&\text{such that, for each }t^{\star}\in(0,T^{\star}],
\\[2pt]
&R_\rho^{\star}(\rho^{\star},\mathbf u^{\star};\mathbf D)
(\phi^{\star})=0,
\qquad\forall\phi^{\star}\in\mathcal H,
\\[2pt]
&R_{\mathbf u}^{\star}(\rho^{\star},\mathbf u^{\star};\mathbf D)
(\mathbf v^{\star})=0,
\qquad\forall\mathbf v^{\star}\in\mathcal V,
\\[2pt]
&\text{with }\rho^{\star}=\rho_0^{\star}\text{ in }\Omega,\ \text{when }t^{\star}=0.
\end{aligned}
\label{eq:app_dimensionless_coupled_problem}
\end{equation}
In the main text, the superscript $\star$ is omitted to simplify the notation, except for the energies $J^{\star}$, $J_E^{\star}$, $J_B^{\star}$ and $J_I^{\star}$.

\subsection{Selection of the phase-field parameters}
\label{app:barred_parameters}

The characteristic element size of the background mesh is defined as
\begin{equation}
h=\max_{K^0\in\mathcal{M}^0}\operatorname{diam}(K^0),
\label{eq:numerical_background_mesh_size}
\end{equation}
where $\mathcal{M}^0$ is the non-trimmed background mesh in the physical
domain and $K^0$ is a parent element, as defined in
Section~\ref{sec:geometric_representation}.
The phase-field parameters are chosen as
\begin{equation}
\lambda=\bar{\lambda}\,h^2,
\qquad
\gamma=\bar{\gamma}\,\gamma_E,
\label{eq:numerical_phasefield_parameters}
\end{equation}
where $\bar{\lambda}$ and $\bar{\gamma}$ are problem-dependent
dimensionless parameters. $\gamma_E$ is defined as the ratio between the
Ginzburg--Landau free energy and the strain energy,
\begin{equation}
\gamma_E=
\frac{
\displaystyle\int_{\Omega^{\mathrm d,\star}}
\bigl[\psi_B^{\star}(\rho_0^{\star})
+(h^{\star})^2\psi_I^{\star}(\rho_0^{\star})\bigr]
\,\mathrm{d}\Omega^{\star}
}{
\displaystyle\int_{\Omega^{\mathrm d,\star}}
\psi_E^{\star}(\rho_0^{\star},\mathbf u_0^{\mathrm d,\star})
\,\mathrm{d}\Omega^{\star}
+
\displaystyle\int_{\Omega^{\mathrm{nd},\star}}
\psi_E^{\star}(\mathbf u_0^{\mathrm{nd},\star})
\,\mathrm{d}\Omega^{\star}
},
\label{eq:numerical_initial_energy_balance}
\end{equation}
where $\mathbf u_0^{\star}$ is the equilibrium displacement for the initial
density $\rho_0^{\star}$. Accordingly, we have $D_2=L_0^2/(\bar{\lambda}h^2)$ and
$D_3=\bar{\gamma}\gamma_E D_2$ as a result of Eq.~\eqref{eq:app_dimensionless_parameters}.

\section{Temporal discretization and adaptive time stepping}
\label{app:time_integration}

We first discretize $\rho$ and $\mathbf u$ spatially as
\begin{equation}
\rho_h(t,\mathbf x)=\sum_{A=1}^{n_\rho}\rho_A(t)N_A(\mathbf x),
\qquad
\mathbf u_h(t,\mathbf x)=\sum_{A=1}^{n_u}\mathbf u_A(t)N_A(\mathbf x),
\end{equation}
where $N_A$ is a B-spline/NURBS basis function, $n_\rho$ and $n_u$ are
the numbers of basis functions for $\rho$ and $\mathbf u$, respectively,
and $\rho_A(t)$ and $\mathbf u_A(t)$ are the control variables at time $t$.

For temporal discretization, we adopt the adaptive generalized-$\alpha$
method to solve the coupled evolution efficiently \cite{Dede_PhaseFieldITO,Sheng_ThinWalledComplex}.
At each time step, a set of nonlinear equations are solved by
Newton--Raphson iterations. The control variables at time step $t_n$ are denoted as
$\mathbf P_n=\{\rho_A(t_n)\}_{A=1}^{n_\rho}$,
$\dot{\mathbf P}_n=\{\dot{\rho}_A(t_n)\}_{A=1}^{n_\rho}$ and
$\mathbf U_n=\{\mathbf u_A(t_n)\}_{A=1}^{n_u}$.
The detailed solution procedure is given in Algorithm~\ref{alg:appB_merged}.

\par\addvspace{\baselineskip}
\begingroup
\renewcommand{\baselinestretch}{1}\small
\algrenewcommand{\algorithmicindent}{1em}
\algrenewcommand{\algorithmicrequire}{\textbf{Given:}}
\captionsetup{
  type=algorithm,
  hypcap=false,
  labelfont=bf,
  font=small,
  justification=raggedright,
  singlelinecheck=false
}
\captionof{algorithm}
{Time step scheme and time adaptivity}
\label{alg:appB_merged}
\label{alg:appB_newton}
\label{alg:appB_adaptivity}
\nobreak\hrule\nobreak\vspace{3pt}

\begin{algorithmic}[1]
\Require $\mathbf U_n$, $\mathbf P_n$, $\dot{\mathbf P}_n$,
$\Delta t_{n,(0)}$ and the parameters $\alpha_m$, $\alpha_f$, $\delta$,
$\mathrm{tol}_R$, $i_{\max}$, $\chi_{NC}$, $\chi_A$, $\chi_{GR}$,
$\mathrm{tol}_A$ and $l_{\max}$.

\For{each $l=1,\ldots,l_{\max}$}
\label{line:appB_adaptive_loop}

  \State \textbf{Given:}
  \label{line:appB_newton_begin}
  $\displaystyle
  \mathbf U_{n+1,(0)}=\mathbf U_n,\quad
  \dot{\mathbf P}_{n+1,(0)}
  =\frac{\delta-1}{\delta}\dot{\mathbf P}_n,\quad
  \mathbf P_{n+1,(0)}=\mathbf P_n.$

  \State Assemble the residuals
  $\mathbf R_{u,(0)}$ and $\mathbf R_{\rho,(0)}$.

  \For{each $i=1,\ldots,i_{\max}$}

    \State Update control variables
    \nopagebreak[4]
    \Statex \hspace{2em}$\displaystyle
    \begin{aligned}
    \mathbf U_{n+1,(i)}
    &=\mathbf U_{n+1,(i-1)},
    \\
    \dot{\mathbf P}_{n+\alpha_m,(i)}
    &=\dot{\mathbf P}_n
    +\alpha_m\bigl(
    \dot{\mathbf P}_{n+1,(i-1)}-\dot{\mathbf P}_n
    \bigr),
    \\
    \mathbf P_{n+\alpha_f,(i)}
    &=\mathbf P_n
    +\alpha_f\bigl(
    \mathbf P_{n+1,(i-1)}-\mathbf P_n
    \bigr).
    \end{aligned}$

    \State Assemble the residuals
    \nopagebreak[4]
    \Statex \hspace{2em}$\displaystyle
    \begin{aligned}
    \mathbf R_{u,(i)}
    &=\mathbf R_u\bigl(
    \mathbf U_{n+1,(i)},
    \dot{\mathbf P}_{n+\alpha_m,(i)},
    \mathbf P_{n+\alpha_f,(i)}
    \bigr),
    \\
    \mathbf R_{\rho,(i)}
    &=\mathbf R_\rho\bigl(
    \mathbf U_{n+1,(i)},
    \dot{\mathbf P}_{n+\alpha_m,(i)},
    \mathbf P_{n+\alpha_f,(i)}
    \bigr).
    \end{aligned}$

    \If{$\displaystyle
    \frac{\|\mathbf R_{u,(i)}\|}
         {\|\mathbf R_{u,(0)}\|}
    <\mathrm{tol}_R$
    and $\displaystyle
    \frac{\|\mathbf R_{\rho,(i)}\|}
         {\|\mathbf R_{\rho,(0)}\|}
    <\mathrm{tol}_R$}

      \State Update control variables
      \nopagebreak[4]
      \Statex \hspace{3em}$\displaystyle
      \begin{aligned}
      \mathbf U_{n+1,(l-1)}
      &=\mathbf U_{n+1,(i-1)},
      \\
      \dot{\mathbf P}_{n+1,(l-1)}
      &=\dot{\mathbf P}_{n+1,(i-1)},
      \\
      \mathbf P_{n+1,(l-1)}
      &=\mathbf P_{n+1,(i-1)}.
      \end{aligned}$

      \State \textbf{break}.

    \Else

            \State Assemble the tangent matrices
      \nopagebreak[4]
      \Statex \hspace{3em}$\displaystyle
      \begin{alignedat}{2}
      \mathbf K_{uu,(i)}
      &=\frac{\partial\mathbf R_u}
             {\partial\mathbf U_{n+1}},
      \quad &
      \mathbf K_{u\rho,(i)}
      &=\alpha_f
      \frac{\partial\mathbf R_u}
           {\partial\mathbf P_{n+\alpha_f}},
      \\[6pt]
      \mathbf K_{\rho u,(i)}
      &=\frac{\partial\mathbf R_\rho}
             {\partial\mathbf U_{n+1}},
      \quad &
      \mathbf K_{\rho\rho,(i)}
      &=\frac{\alpha_m}{\delta\Delta t_{n,(l-1)}}
      \frac{\partial\mathbf R_\rho}
           {\partial\dot{\mathbf P}_{n+\alpha_m}}
      +\alpha_f
      \frac{\partial\mathbf R_\rho}
           {\partial\mathbf P_{n+\alpha_f}}.
      \end{alignedat}$

      \State Solve the linear system
      \nopagebreak[4]
      \Statex \hspace{3em}$\displaystyle
      \begin{aligned}
      \mathbf K_{uu,(i)}\Delta\mathbf U_{n+1,(i)}
      +\mathbf K_{u\rho,(i)}\Delta\mathbf P_{n+1,(i)}
      &=-\mathbf R_{u,(i)},
      \\
      \mathbf K_{\rho u,(i)}\Delta\mathbf U_{n+1,(i)}
      +\mathbf K_{\rho\rho,(i)}\Delta\mathbf P_{n+1,(i)}
      &=-\mathbf R_{\rho,(i)}.
      \end{aligned}$

      \State $\displaystyle
      \begin{aligned}[t]
      &\text{Update control variables}
      \\[2pt]
      &\mathbf U_{n+1,(i)}
      =\mathbf U_{n+1,(i-1)}
      +\Delta\mathbf U_{n+1,(i)},
      \\
      &\dot{\mathbf P}_{n+1,(i)}
      =\dot{\mathbf P}_{n+1,(i-1)}
      +\frac{1}{\delta\Delta t_{n,(l-1)}}
       \Delta\mathbf P_{n+1,(i)},
      \\
      &\mathbf P_{n+1,(i)}
      =\mathbf P_{n+1,(i-1)}
      +\Delta\mathbf P_{n+1,(i)}.
      \end{aligned}$

    \EndIf
  \EndFor
  \label{line:appB_newton_end}

  \State Compute control variables
  $\mathbf U_{n+1,(l-1)}^{\mathrm{BE}}$,
  $\dot{\mathbf P}_{n+1,(l-1)}^{\mathrm{BE}}$ and
  $\mathbf P_{n+1,(l-1)}^{\mathrm{BE}}$
  using backward Euler: repeat lines
  \ref{line:appB_newton_begin}--\ref{line:appB_newton_end}
  from $(\mathbf U_n,\mathbf P_n,\dot{\mathbf P}_n)$
  with $\alpha_m=\alpha_f=\delta=1$ and
  $\Delta t_{n,(l-1)}$.

  \If{generalized-$\alpha$ or backward Euler are not convergent}

    \State $\Delta t_{n,(l)}
    =\chi_{NC}\Delta t_{n,(l-1)}$.

    \State Go to line~\ref{line:appB_adaptive_loop}.

  \Else

    \State Evaluate the error
    \nopagebreak[4]
    \Statex \hspace{2em}$\displaystyle
    \begin{aligned}
    e_{n+1,(l-1)}
    ={}&
    \frac{
      \|\mathbf U_{n+1,(l-1)}
      -\mathbf U_{n+1,(l-1)}^{\mathrm{BE}}\|
    }{
      \|\mathbf U_{n+1,(l-1)}\|
    }
    +
    \frac{
      \|\mathbf P_{n+1,(l-1)}
      -\mathbf P_{n+1,(l-1)}^{\mathrm{BE}}\|
    }{
      \|\mathbf P_{n+1,(l-1)}\|
    }.
    \end{aligned}$

    \State Update time step
    \nopagebreak[4]
    \Statex \hspace{2em}$\displaystyle
    \Delta t_{n,(l)}
    =
    \min\!\left\{
      \chi_A
      \sqrt{\frac{\mathrm{tol}_A}{e_{n+1,(l-1)}}},
      1+\chi_{GR}
    \right\}
    \Delta t_{n,(l-1)}.$

    \If{$e_{n+1,(l-1)}<\mathrm{tol}_A$}

      \State Update control variables and time step
      \nopagebreak[4]
      \Statex \hspace{3em}$\displaystyle
      \begin{aligned}
      \mathbf U_{n+1}
      &=\mathbf U_{n+1,(l-1)},
      &
      \mathbf P_{n+1}
      &=\mathbf P_{n+1,(l-1)},
      \\
      \dot{\mathbf P}_{n+1}
      &=\dot{\mathbf P}_{n+1,(l-1)},
      &
      \Delta t_n
      &=\Delta t_{n,(l)}.
      \end{aligned}$

      \State \Return
      $\mathbf U_{n+1}$, $\mathbf P_{n+1}$,
      $\dot{\mathbf P}_{n+1}$ and $\Delta t_n$.

    \Else

      \State Go to line~\ref{line:appB_adaptive_loop}.

    \EndIf
  \EndIf
\EndFor

\end{algorithmic}
\nobreak\vspace{3pt}\hrule
\endgroup
\par\addvspace{\baselineskip}

For the generalized-$\alpha$ method, we set
$\alpha_m=\frac{5}{6}$, $\alpha_f=\delta=\frac{2}{3}$. For Newton's method, we set
$\mathrm{tol}_R=10^{-4}$ and $i_{\max}=10$.
For time adaptivity, we use
$\chi_{NC}=0.5$, $\chi_A=0.85$, $\chi_{GR}=0.2$,
$\mathrm{tol}_A=10^{-3}$ and $l_{\max}=30$, starting with
$\Delta t_0=10^{-12}$.

\section*{References}
\begingroup
\renewcommand{\section}[2]{}%

\endgroup

\end{document}